\documentclass[english]{article}
\usepackage[T1]{fontenc}
\usepackage[latin9]{inputenc}
\usepackage{textcomp}
\usepackage{mathrsfs}
\usepackage{amsmath}
\usepackage{amssymb}
\usepackage{cancel}
\usepackage{esint}
\usepackage{babel}
\begin{document}
\title{Conservation of Singularities in Functional Equations Associated to
Collatz-Type Dynamical Systems; or, Dreamcatchers for Hydra Maps}
\author{M.C. Siegel}

\maketitle
\tableofcontents{}

\vphantom{}\emph{Disclaimer: This mathematics in this paper were
raised free-range and organic, without the use of the Axiom of Choice.}

\pagebreak{}

\section*{Introduction}

\subsection{The Author's Idiosyncratic Notation}

The main results and arguments given in this paper tend to be computationally
intensive and/or require working with a significant amount of data
in mind. To that end, here are the notational conventions the reader
should be aware of prior to reading this paper. Additional notional
conventions will be introduced in the text as the need arises.

\textbullet{} $\overset{\textrm{def}}{=}$ (read ``is defined as'')
is used to denote the definition of a symbol.

\textbullet{} $\mathbb{P}$ is used to denote the set of prime numbers
($\mathbb{P}\overset{\textrm{def}}{=}\left\{ 2,3,5,\ldots\right\} $).

\textbullet{} Given an real number $r$, the symbol $\mathbb{N}_{r}$
denotes \textbf{the set of all integers} \textbf{greater than or equal
to }\emph{r}. $\mathbb{C}$ denotes the complex numbers.

\textbullet{} $\mathbb{D}$ denotes the open disk in $\mathbb{C}$
centered at $0$ of radius $1$ (the ``(open) unit disk''). $\partial\mathbb{D}$
is the unit circle in $\mathbb{C}$ centered at $0$.

\textbullet{} $\mathbb{H}_{+}$ denotes the half-plane of all $z\in\mathbb{C}$
with $\textrm{Re}\left(z\right)>0$ (the \textbf{right half-plane});
$\mathbb{H}_{+i}$ denotes the half-plane of all $z\in\mathbb{C}$
with $\textrm{Im}\left(z\right)>0$ (the \textbf{upper half-plane}).
Over-bars ($\overline{\mathbb{H}}_{+}$, $\overline{\mathbb{\mathbb{H}}}_{+i}$)
denote the topological closures of $\mathbb{H}_{+}$ and $\mathbb{\mathbb{H}}_{+i}$,
respectively.

\textbullet{} Given a ring $\mathcal{R}$ (typically, $\mathcal{R}$
will be $\mathbb{C}$, or $\overline{\mathbb{Z}}$, or $\mathbb{Z}_{p}$,
or some other nice, locally-compact topological ring), and given elements
$a,b\in\mathcal{R}$ and a subset $S\subseteq\mathcal{R}$, we write:
\[
a+bS\overset{\textrm{def}}{=}\left\{ a+bs:s\in S\right\} 
\]
Thus, for example, for $a\in\mathbb{C}$ and $r\in\mathbb{R}>0$,
$a+r\mathbb{D}$ is the open disk in $\mathbb{C}$ of radius $r$
centered at $a$; $a+r\partial\mathbb{D}$, meanwhile, is the circle
in $\mathbb{C}$ of radius $r$ centered at $a$. In a more number-theoretic
vein, the statement that real numbers $x,y\in\mathbb{R}$ are congruent
to one another modulo an integer $n$ (i.e., $x-y=mn$ for some $m\in\mathbb{Z}$)
can be written as $x\in y+n\mathbb{Z}$, $y\in x+n\mathbb{Z}$, or
$x-y\in n\mathbb{Z}$.

\textbullet{} Given a function or map (ex: $P$), $P^{\circ n}$ denotes
the $n$th iterate of $P$:
\[
P^{\circ n}=\underbrace{P\circ P\circ...\circ P}_{n\textrm{ times}}
\]

\textbullet{} Given two expressions $A$ and $B$ (or statements,
or equations, etc.), the notation $A\Leftrightarrow B$ or, equivalently:
\[
\begin{array}{c}
A\\
\Updownarrow\\
B
\end{array}
\]
denotes logical equivalence (i.e., $A$ is true if and only if $B$
is true). Single arrows $\Leftarrow,\Rightarrow,\Uparrow,\Downarrow$
denote (mono-directional) if-then implications.

\textbullet{} With regards to inequalities and set-inclusions, the
symbols $<$, $>$, $\subset$, and $\supset$ will \textbf{only be
used in the strict sense}. That is, 
\[
<\Leftrightarrow\lneqq
\]
\[
\subset\Leftrightarrow\subsetneqq
\]
and so on and so forth.

\textbullet{} Let $a\in\mathbb{N}_{1}$ and let $r\in\mathbb{Q}$.
Then:
\[
\xi_{a}^{r}\overset{\textrm{def}}{=}e^{\frac{2\pi ir}{a}}
\]

\textbullet{} Let $a\in\mathbb{N}_{1}$. The symbol $\overset{a}{\equiv}$
is used to denote \textbf{congruence of real numbers modulo }\textbf{\emph{a}}.

\textbullet{} Let $a\in\mathbb{N}_{1}$. For sake of concreteness,
quotient rings $\mathbb{Z}/a\mathbb{Z}$ (where $a\in\mathbb{N}_{2}$)
will always be viewed as the set $\left\{ 0,\ldots,a-1\right\} \subseteq\mathbb{N}_{0}$
equipped with special arithmetic laws. As such, given $a\in\mathbb{N}_{1}$
and any $x\in\mathbb{R}$, the expression $\left[x\right]_{a}$ is
written to denote the unique real number in the interval $\left[0,a\right)$
so that $x=an+\left[x\right]_{a}$ for some integer $n\in\mathbb{Z}$.
Given an integer $n\in\mathbb{Z}$ which is co-prime to $a$ (that
is, $\gcd\left(a,n\right)=1$), the expression $\left[n\right]_{a}^{-1}$
denotes the unique integer in $\left\{ 0,\ldots,a-1\right\} $ so
that:
\[
n\left[n\right]_{a}^{-1}\overset{a}{\equiv}1
\]
that is, so that the product of $n$ and $\left[n\right]_{a}^{-1}$
is an integer which is congruent to $1$ mod $a$.

\textbullet{} In the aim of keeping notation as independent of context
as pleasantly possible, different fonts will be used for Archimedean
(real or complex) and non-Archimedean quantities (elements of $\overline{\mathbb{Z}}$,
$\mathbb{Z}_{p}$, $\mathbb{Q}_{p}$, extensions thereof, etc.). All
non-Archimedean variables and constants are written using lower case
letters in the $\mathfrak{fraktur}$ font. This notation will also
be used for Haar measures for non-Archimedean spaces (for example,
writing $d\mathfrak{z}_{p}$ to denote the unit-normalized Haar measure
on $\mathbb{Z}_{p}$).

\textbullet{} Let $p\in\mathbb{P}$, let $\mathfrak{y}\in\mathbb{Z}_{p}$,
and let $n\in\mathbb{N}_{1}$. The expression $\left[\mathfrak{y}\right]_{p^{n}}$
denotes the unique integer in $\left\{ 0,\ldots,p^{n}-1\right\} $
for which $\mathfrak{y}\overset{p^{n}}{\equiv}\left[\mathfrak{y}\right]_{p^{n}}$

\textbullet{} Let $p\in\mathbb{P}$. The expression $\mathbb{Z}_{p}^{\times}$
denotes the set of all multiplicatively invertible elements of $\mathbb{Z}_{p}$.

\textbullet{} Let $G$ be a locally compact abelian group with Haar
measure $dx$, and let $r$ be a real number (or infinity) satisfying
$1\leq r\leq\infty$. $L^{r}\left(G\right)$ is written to denote
the $\mathbb{C}$-linear space of all complex valued functions $f:G\rightarrow\mathbb{C}$
so that:
\[
\int_{G}\left|f\left(x\right)\right|^{r}dx<\infty
\]

\textbullet{} Let $a\in\mathbb{R}$. Then, the expressions $\lim_{x\downarrow a}$
(resp. $\lim_{x\uparrow a}$) denote the limit of $x$ as $x$ decreases
to $a$ (i.e., approaches $a$ from above) (resp. the limit of $x$
as $x$ increases to $a$ (i.e., approaches $a$ from below)).

\textbullet{} The\textbf{ ring of profinite integers} is denoted by
$\widetilde{\mathbb{Z}}$; recall that there is an isomorphism of
topological rings given by:
\[
\widetilde{\mathbb{Z}}\cong\prod_{p\in\mathbb{P}}\mathbb{Z}_{p}
\]

\textbullet{} When working in a measure space $\left(X,\Sigma,\mu\right)$,
The symbol $\tilde{\forall}$ is used to denote ``for almost all'';
so, the statement ``$S$ is true $\tilde{\forall}x\in X$'' means
that $S$ is true for all $x\in X$ except possibly for those $x$
belonging in some subset of $X$ with zero $\mu$-measure.

\textbullet{} (\textbf{Iverson bracket notation}) Given a statement
$S$, the Iverson bracket $\left[S\right]$ denotes a quantity which
is $1$ if $S$ is true and which is $0$ if $S$ is false.

\emph{Remark}: Let $X$ be a set, let $S\left(x\right)$ be a statement
about elements $x\in X$, and let $f:X\rightarrow\mathbb{C}$ be a
function. Write $f\left(S\left(x\right)\right)$ to denote what $f\left(x\right)$
would be under the assumption that $S\left(x\right)$ were true. Then,
the \textbf{evaluation property }of Iverson brackets is that:
\[
f\left(x\right)\left[S\left(x\right)\right]=f\left(S\left(x\right)\right)\left[S\left(x\right)\right]
\]
Example:

\[
e^{2\pi ix}\left[x\overset{1}{\equiv}\frac{1}{a}\right]=\xi_{a}\left[x\overset{1}{\equiv}\frac{1}{a}\right]
\]
where, recall, $\xi_{a}=e^{2\pi i/a}$.

\textbullet{} Given a set $S$, the \textbf{indicator function }of
$S$ is denoted $\mathbf{1}_{S}$:
\[
\mathbf{1}_{S}\left(x\right)\overset{\textrm{def}}{=}\begin{cases}
1 & \textrm{if }x\in S\\
0 & \textrm{else}
\end{cases}=\left[x\in S\right]
\]
This notation will be modified in §3.1.2 to deal with co-sets of $\mathbb{Z}$
in $\mathbb{R}$.

\pagebreak{}

\subsection{Background}

The impetus for this dissertation came from the author's interest
(read ``obsession'') with the infamous \textbf{Collatz Conjecture},
also known as the $3n+1$ problem; defining the function:
\[
C\left(n\right)\overset{\textrm{def}}{=}\begin{cases}
\frac{n}{2} & \textrm{if }n\textrm{ is even}\\
3n+1 & \textrm{if }n\textrm{ is odd}
\end{cases}
\]
on the positive integers, does there exist an integer $n\geq1$ which
$C$ does \emph{not} eventually iterate to $1$? The orbits of integers
under $C$ rise and fall with an almost organic randomness, like the
movements of hailstones in the depths of storm clouds. Raw computational
force has shown that every positive integer up to some ungodly upper
bound gets mapped to $1$; and yet, the answer remains elusive. To
the author's overactive imagination, the conjecture is like the Hydra
of Greek myth. Try as we might to slay it, no matter how many of its
heads we lop off, more and more continue to grow in its place; endless
waves of complex simplicity. Barring a miracle, for the time being,
mathematics will have to content itself with striving for incremental
progress on understanding the phenomena which underlie the Collatz
map and its many relatives.

One such source of progress is the notion of a ``sufficient set'',
for which Collatz-relevant results have been proven by Andaloro (Cited
in (Monks. 2006)), and subsequently refined by K. Monks in \emph{The
Sufficiency of Arithmetic Progressions for the $3x+1$} \emph{Conjecture}
(Monks. 2006). Stated informally, given a map $H:\mathbb{N}_{0}\rightarrow\mathbb{N}_{0}$
and a conjecture \textbf{C} about the behavior of $H$ on $\mathbb{N}_{0}$,
a set $S\subseteq\mathbb{N}_{0}$ is ``sufficient'' if proving that
\textbf{C }holds for all $n\in S$ implies that \textbf{C }holds for
all $n\in\mathbb{N}_{0}$. The cited authors' work deals exclusively
with the Collatz Conjecture. Monks' conclusion---the stronger of
the two---is that infinite arithmetic progressions (IAPs) are sufficient
for the Collatz Conjecture; \emph{that is to say, if you can prove
there are integers $a\geq1$, $b\geq0$, and $\nu\geq0$ so that every
element of the IAP:
\[
\left\{ an+b:n\in\mathbb{N}_{\nu}\right\} \cap\mathbb{N}_{1}
\]
is eventually sent to $1$ by the Collatz map, then you will have
succeeded in proving the Collatz Conjecture}. Equivalently---and
this is how we shall speak of sufficiency results in this paper---\emph{if
an orbit class of the Collatz map in $\mathbb{N}_{0}$ contains an
IAP, then that orbit class must contain all positive integers}. For
reasons to be discussed later on, the author call results of this
form \textbf{Rationality Theorems}.

The genesis of this paper was an aspiration to generalize Monks' result
to other Collatz-type maps, using more broadly applicable analytic
methods, rather than Monks' intricate use of elementary arithmetic
and congruences in manipulating $2$s and $3$s. The present paper
only partially succeeds in realizing that aspiration. The result obtained
herein is, essentially, a weaker form of Monks' sufficiency result,
albeit one applicable to a relatively large family of Collatz-type
maps---the eponymous ``Hydra maps'' of the paper's title. The use
of the qualifier ``essentially'' in the previous sentence cannot
be dropped. The methods used to prove these results are twice-removed
from the context of the original problem (arithmetic dynamical systems
on the integers), and there are several apparently non-trivial analytic
technicalities that must be dealt with in order to fully justify applying
our findings to the problems that birthed them. These technicalities
will be discussed throughout the paper, as the need arises.

With that in mind, our generalization (a \textbf{Weak Rationality
Theorem}) can be roughly stated as: \emph{for any Hydra map $H$ satisfying
the appropriate conditions, if an orbit class $V\subseteq\mathbb{N}_{0}$
of $H$ is the union of a finite set and finitely many arithmetic
progressions of infinite length, then $\mathbb{N}_{0}\backslash V$
contains at most finitely many elements}. This and other results are
detailed in §4.2; however, the terminology used in Chapter 1 and §4.1
will be needed to fully understand them. Although not detailed in
this paper, the author has generalized this result to Hydra maps defined
on the integer rings of an arbitrary number field (a finite-degree
field extension of $\mathbb{Q}$). Additionally, the methods presented
here have been extended to deal with deeper features of the orbit
classes of Hydra maps, however, these extensions produce as many questions
as they do answers. Both of these generalizations (as well as their
combination) shall be dealt with in subsequent publications.

That being said, the author believes that the most important ``result''
of the present paper is not the \textbf{Weak Rationality Theorem},
but rather the tool used to obtain it: the\textbf{ Dreamcatcher}.
This apparently novel construction can be used tool to analyze the
properties of solutions of functional equations over the complex numbers.
It is hoped that Dreamcatchers and the Fourier/Harmonic-analytic methods
employed in conjunction with them will be of interest not only for
their application to future studies of Collatz and Collatz-type conjectures,
but also their potential use in any situations where functional equations
of holomorphic functions arise: transcendental number theory, analytic
combinatorics, the theory of automata, and the like.\pagebreak{}

\subsection{An Overview of the Paper}

Thematically speaking, the paper has three central concerns.

\vphantom{}I.\textbf{ Permutation Operators} \textbf{\& Set-series}:
These---the first part of the ``twice-removal''---arise by studying
the Collatz Map and its kin not as maps on $\mathbb{Z}$, but as linear
operators on spaces of holomorphic functions defined by convergent
infinite series of a specified type. Though the author arrived at
this idea independently in late March of 2017, this approach is not
new; its first appearance in the literature (to our knowledge) was
in 1995, in a paper by Meinardus and Berg (Berg \& Meinardus. 1995).

As is common, rather than working with the Collatz map as defined
at the start of the paper, consider the ``shortened'' version of
the Collatz map:
\[
H_{3}\left(n\right)\overset{\textrm{def}}{=}\begin{cases}
\frac{n}{2} & \textrm{if }n\overset{2}{\equiv}0\\
\frac{3n+1}{2} & \textrm{if }n\overset{2}{\equiv}1
\end{cases}
\]
where the use of the notation $H_{3}$ is part of the author's convention
to write $H_{p}$ to denote the map:
\[
H_{p}\left(n\right)\overset{\textrm{def}}{=}\begin{cases}
\frac{n}{2} & \textrm{if }n\overset{2}{\equiv}0\\
\frac{pn+1}{2} & \textrm{if }n\overset{2}{\equiv}1
\end{cases}
\]
where $p$ is an odd prime. The idea is to consider $H_{3}$ acting
on power series. Specifically, letting $f:\mathbb{D}\rightarrow\mathbb{C}$
be a holomorphic defined by a convergent power series: 
\[
f\left(z\right)=\sum_{n=0}^{\infty}c_{n}z^{n}
\]
we can define two linear operators:
\[
\mathcal{Q}_{3}\left\{ \sum_{n=0}^{\infty}c_{n}z^{n}\right\} \overset{\textrm{def}}{=}\sum_{n=0}^{\infty}c_{H_{3}\left(n\right)}z^{n}
\]
\[
\mathcal{P}_{3}\left\{ \sum_{n=0}^{\infty}c_{n}z^{n}\right\} \overset{\textrm{def}}{=}\sum_{n=0}^{\infty}c_{n}z^{H_{3}\left(n\right)}
\]
These operators do not appear to have any established name; the author
calls them the \textbf{Permutation Operators }induced by $H_{3}$
of the \textbf{first }($\mathcal{Q}_{3}$) and \textbf{second }($\mathcal{P}_{3}$\textbf{)
kinds}, respectively. The focus of this paper will be on $\mathcal{Q}$,
rather than $\mathcal{P}$. It remains to be seen whether studying
$\mathcal{P}$ will prove to be as productive as studying $\mathcal{Q}$.

As is discussed in Chapter 1, the Permutation Operators induced by
a dynamical system translate questions about the system's invariant
sets into questions about the fixed points of the Permutation Operators.
For example, the functions fixed by $\mathcal{Q}_{3}$ correspond
exactly to the different possible \textbf{$H_{3}$-invariant sets}
in $\mathbb{N}_{0}$; that is, sets $V\subseteq\mathbb{N}_{0}$ for
which $V=H_{3}\left(V\right)=H_{3}^{-1}\left(V\right)$ (this is the
\textbf{Permutation Operator Theorem}).\textbf{ }The fixed points
of permutation operators are characterized by functional equations;
the fixed points of $\mathcal{Q}_{3}$, for instance, are exactly
those functions $f\left(z\right)$ holomorphic in a neighborhood of
$0\in\mathbb{C}$ that solve the functional equation\footnote{This exact equation (albeit with a $\lambda$ instead of a $\xi_{3}$,
and with an $h$ instead of an $f$) can be found in (Berg \& Meinardus.
1995).}:
\[
f\left(z\right)=f\left(z^{2}\right)+\frac{z^{-\frac{1}{3}}}{3}\sum_{k=0}^{2}\xi_{3}^{2k}f\left(\xi_{3}^{k}z^{2/3}\right)
\]
The solutions of this equation will be linear combinations of a type
of ordinary generating function the author calls \textbf{set-series}.
Given $V\subseteq\mathbb{N}_{0}$, the (ordinary) set-series of $V$
is the function $\varsigma_{V}\left(z\right)$ defined by:
\[
\varsigma_{V}\left(z\right)\overset{\textrm{def}}{=}\sum_{v\in V}z^{v}
\]
When $V\subseteq\mathbb{N}_{0}$ is a set satisfying $V=H_{3}\left(V\right)=H_{3}^{-1}\left(V\right)$
(that is, when $V$ is $H_{3}$-invariant), $\varsigma_{V}$ will
be a fixed point of $\mathcal{Q}_{3}$. Changing perspectives from
the open unit disk in $\mathbb{C}$ to the upper half-lane by letting
$\psi\left(z\right)\overset{\textrm{def}}{=}f\left(e^{2\pi iz}\right)$
gives the equivalent equation:
\[
\psi\left(z\right)=\psi\left(2z\right)+\frac{e^{-\frac{2\pi iz}{3}}}{3}\sum_{k=0}^{2}\xi_{3}^{2k}\psi\left(\frac{2z+k}{3}\right)
\]
whose solutions are linear combinations of \textbf{Fourier set-series}:
\[
\psi_{V}\left(z\right)\overset{\textrm{def}}{=}\sum_{v\in V}e^{2\pi ivz}
\]
where, again, $V$ is an $H_{3}$-invariant set.

\vphantom{}II. \textbf{Dreamcatchers \& Singularity Conservation
Laws}: The \textbf{Hardy-Littlewood Tauberian Theorem} from analytic
number theory establishes a correspondence between the arithmetic
properties of $V$ (its density; its distribution modulo a given positive
integer, and so on) and the singularities\emph{ }of $V$'s set-series
($\varsigma_{V}$, $\psi_{V}$, etc.). We will show that this correspondence
is strong enough that, provided the appropriate technical details
can be borne out, \emph{solutions of the above functional equations
are determined entirely by their singularities on the boundary of
their domain of holomorphy}.

Rather than attacking the functional equations directly, we exploit
the fact that $\mathcal{Q}_{H}$ (the permutation operator of the
first kind induced by a hydra map $H$) satisfies a \textbf{Singularity
Conservation Law}; in simplest terms, this is the observation that
if $\mathcal{Q}_{H}\left\{ f\right\} \left(z\right)$ possesses a
singularity of a given type (degree, asymptotic growth rate, etc.),
then $f$ must possess a singularity of that type. As such, we can
simplify the problem of solving for fixed points of $\mathcal{Q}_{H}$
by instead looking for functions whose sets of singularities are fixed
by $\mathcal{Q}$$_{H}$. To do this, we generalize the classical
complex-analytic concept of a pole and its residue to construct a
\textbf{Dreamcatcher}: mathematical objects that record the locations
and ``virtual'' residues of a given series' singularities on the
boundary of its region of convergence.

For example, given an appropriate $1$-periodic function $\psi$ holomorphic
on the upper half plane $\mathbb{H}_{+i}$, the dreamcatcher of $\psi$
is the function:
\[
R\left(x\right)\overset{\textrm{def}}{=}\lim_{y\downarrow0}y\psi\left(x+iy\right)
\]
defined for all $x\in\mathbb{Q}$ for which the above limit exists.
If $x$ is a simple pole (in the classical sense), $R\left(x\right)$
will be the residue (again, in the classical sense) of $\psi$ at
$x$. However, as is demonstrated by example at the start of Chapter
2, the above definition has the advantage of being defined at singularities
which have the same asymptotic growth rate as a simple pole ($\psi\left(x+iy\right)$
growing like $\frac{1}{y}$ as $y$ decreases to $0$), but which
are not true poles, due to the presence of lower-order growth rates,
or in the case where $x$ is a non-isolated singularity of $\psi$.

The surprise is that, as defined, dreamcatchers enjoy a non-trivial
level of regularity. Fixing $\psi$ as described above, let $T$ denote
the $1$-periodic set of $x\in\mathbb{Q}$ for which $R$ exists (note
that $T$ depends on $\psi$). Then, as will be shown in Chapter 2
(the \textbf{Square Sum Lemma}),\textbf{ }$R$ belongs to the Hilbert
space $L^{2}\left(T/\mathbb{Z}\right)$; that is, space of all functions
$R:T\rightarrow\mathbb{C}$ satisfying:
\[
R\left(x+1\right)=R\left(x\right),\textrm{ }\forall x\in T
\]
\[
\sum_{x\in\left[0,1\right)\cap T}\left|R\left(x\right)\right|^{2}<\infty
\]

By applying the same limiting procedure to $\mathscr{Q}_{H}$---the
analogue of $\mathcal{Q}_{H}$ for Fourier series on the upper half
plane, rather than power series on the unit disk---one obtains the
\textbf{Dreamcatcher operator }$\mathfrak{Q}_{H}:L^{2}\left(\mathbb{Q}/\mathbb{Z}\right)\rightarrow L^{2}\left(\mathbb{Q}/\mathbb{Z}\right)$;
this is the second step of the ``twice-removal''. Formulated in
this manner, the \textbf{Singularity Conservation Law }can, roughly,
be realized as the statement ``$\psi\in\textrm{Ker}\left(1-\mathscr{Q}_{H}\right)$
implies $R\in\textrm{Ker}\left(1-\mathfrak{Q}_{H}\right)$''; that
is to say, the dreamcatcher of a fixed point of $\mathscr{Q}_{H}$
is itself a fixed point of $\mathfrak{Q}_{H}$. As such, by studying
the fixed points of $\mathfrak{Q}_{H}$ over $L^{2}\left(\mathbb{Q}/\mathbb{Z}\right)$,
we can (modulo the technicalities) establish control on the singular
behavior of $\psi$ (and thus, on the distribution in $\mathbb{N}_{0}$
of the $H$-invariant set represented by $\psi$) by solving the fixed
point functional equation of $\mathfrak{Q}_{H}$, one that will hopefully
be less resilient to analysis than those of $\mathscr{Q}_{H}$ or
$\mathcal{Q}_{H}$.

\vphantom{}III. \textbf{Profinite Harmonic Analysis}: Reformulating
our functional equations in terms of dreamcatchers on $L^{2}\left(\mathbb{Q}/\mathbb{Z}\right)$
allows us to take advantage of the tools and structural properties
that $L^{2}\left(\mathbb{Q}/\mathbb{Z}\right)$ enjoys as the Hilbert
space of square-integrable complex-valued functions on a locally compact
abelian group (LCAG), such as orthogonality, and a Fourier Transform.
Fortunately, for the case dealt with in this paper, the problem on
$L^{2}\left(\mathbb{Q}/\mathbb{Z}\right)$ is indeed somewhat simpler
than the original functional equation from whence it came. We will
show that the fixed points $R$ of $\mathfrak{Q}_{H}$ with finite
support in $\mathbb{Q}/\mathbb{Z}$ (i.e., for which $R\left(x\right)=0$
for all but finitely many $x\in\mathbb{Q}\cap\left[0,1\right)$) must
be of the form $R\left(x\right)=R\left(0\right)\mathbf{1}_{\mathbb{Z}}\left(x\right)$
for a large family of $H$ (the \textbf{Finite Dreamcatcher Theorem}.

Proving the \textbf{Finite Dreamcatcher Theorem} requires a mix of
elementary \textbf{$p$-adic analysis} and the theory of \textbf{Pontryagin
Duality }for $\mathbb{Q}/\mathbb{Z}$. If every branch of $H$ entails
dividing out by some integer $\varrho$, using $\varrho$-adic analysis,
we can show $R\in\textrm{Ker}\left(1-\mathfrak{Q}_{H}\right)$ with
finite support must vanish for all rational numbers whose denominators
are multiples of $\varrho$. In particular, any such $R\in\textrm{Ker}\left(1-\mathfrak{Q}_{H}\right)$
must then have a trivial component ``on $\varrho$'' with respect
to the tensor product decomposition:
\[
L^{2}\left(\mathbb{Q}/\mathbb{Z}\right)\cong L^{2}\left(\mathbb{Z}\left[\frac{1}{\varrho}\right]/\mathbb{Z}\right)\otimes L^{2}\left(\mathbb{Q}_{\bcancel{\varrho}}/\mathbb{Z}\right)
\]
where: 
\[
\mathbb{Q}_{\bcancel{\varrho}}\overset{\textrm{def}}{=}\mathbb{Q}\cap\mathbb{Z}_{\varrho}=\left\{ x\in\mathbb{Q}:\left|x\right|_{\varrho}\leq1\right\} 
\]
(read {[}Q off rho{]}) with:
\[
\mathbb{Q}/\mathbb{Z}\cong\left(\mathbb{Z}\left[\frac{1}{\varrho}\right]/\mathbb{Z}\right)\times\left(\mathbb{Q}_{\bcancel{\varrho}}/\mathbb{Z}\right)
\]
That is to say, writing $x\in\mathbb{Q}/\mathbb{Z}$ as $x=\left(t,y\right)\in\left(\mathbb{Z}\left[\frac{1}{\varrho}\right]/\mathbb{Z}\right)\times\left(\mathbb{Q}_{\bcancel{\varrho}}/\mathbb{Z}\right)$,
there is a function $F\in L^{2}\left(\mathbb{Q}_{\bcancel{\varrho}}/\mathbb{Z}\right)$
so that:
\[
R\left(x\right)=R\left(t,y\right)=\left(\mathbf{1}_{\mathbb{Z}}\otimes F\right)\left(t,y\right)=\mathbf{1}_{\mathbf{\mathbb{Z}}}\left(t\right)F\left(y\right),\textrm{ }\forall x=\left(t,y\right)\in\left(\mathbb{Z}\left[\frac{1}{\varrho}\right]/\mathbb{Z}\right)\times\left(\mathbb{Q}_{\bcancel{\varrho}}/\mathbb{Z}\right)
\]
We describe this situation by saying that $R$ is ``supported off
$\varrho$''.

The existence of this decomposition overcomes a central obstacle faced
when studying $R$ over $L^{2}\left(\mathbb{Q}/\mathbb{Z}\right)$,
or---even more palpably---when studying its Fourier Transform, $\check{R}\in L^{2}\left(\widetilde{\mathbb{Z}}\right)$.
Indeed:
\[
R\left(x\right)=\left(\mathbf{1}_{\mathbb{Z}}\otimes F\right)\left(t,y\right)
\]
implies:
\[
\check{R}\left(\mathfrak{z}\right)=\left(1\otimes\check{F}\right)\left(\mathfrak{t},\mathfrak{y}\right)=\check{F}\left(\mathfrak{y}\right),\textrm{ }\tilde{\forall}\mathfrak{z}=\left(\mathfrak{t},\mathfrak{y}\right)\in\mathbb{Z}_{\varrho}\times\widetilde{\mathbb{Z}}_{\bcancel{\varrho}}
\]
where $\mathbb{Z}_{\varrho}$ is the ring of $\varrho$-adic integers
(the Pontryagin dual of $\mathbb{Z}\left[\frac{1}{\varrho}\right]/\mathbb{Z}$),
and where: 
\[
\widetilde{\mathbb{Z}}_{\bcancel{\varrho}}\overset{\textrm{def}}{=}\prod_{\begin{array}{c}
p\in\mathbb{P}\\
p\nmid\varrho
\end{array}}\mathbb{Z}_{p}
\]
(read {[}Z tilde off rho{]}), with:
\[
\widetilde{\mathbb{Z}}\cong\mathbb{Z}_{\varrho}\times\widetilde{\mathbb{Z}}_{\bcancel{\varrho}}
\]
As such, we can think of $\check{R}$ as an element of $L^{2}\left(\widetilde{\mathbb{Z}}_{\bcancel{\varrho}}\right)$;
this removes the difficulty posed by the fact that, as a subring of
$\widetilde{\mathbb{Z}}$, $\widetilde{\mathbb{Z}}_{\bcancel{\varrho}}$
has zero measure.

With the aid of the off-$\varrho$ decomposition, we can show that
the support in $\mathbb{Q}_{\bcancel{\varrho}}/\mathbb{Z}$ of any
such $R$ must be invariant under multiplication by $\varrho$ mod
$1$ (a group automorphism of $\mathbb{Q}_{\bcancel{\varrho}}/\mathbb{Z}$).
This invariance allows us to show that the support of any $R\in\textrm{Ker}\left(1-\mathfrak{Q}_{H}\right)$
supported off $\varrho$ must be unbounded in $\mu$-adic magnitude,
where $\mu$ is at least one (but possibly more) of the integers so
that there is a $j\in\left\{ 0,\ldots,\varrho-1\right\} $ for which:
\[
H\left(n\right)=\frac{\mu}{\varrho}n+\frac{b}{d},\textrm{ }\forall n\in\varrho\mathbb{N}_{0}+j
\]
for some constants $b,d\in\mathbb{Z}$, with $d\neq0$. In particular,
this shows that if $R\in\textrm{Ker}\left(1-\mathfrak{Q}_{H}\right)$
is finitely supported on $\mathbb{Q}/\mathbb{Z}$, then $R$ must
be trivial:

\[
R\left(x\right)=R\left(0\right)\mathbf{1}_{\mathbb{Z}}\left(x\right),\textrm{ }\forall x\in\mathbb{Q}
\]
which is, of course, exactly what is asserted by the \textbf{Finite
Dreamcatcher Theorem}.

\pagebreak{}

\section{Set-Series, Hydra Maps, and Permutation Operators}

The principal \emph{leitmotif }of our study can be stated as: given
a map that generates a discrete (and, generally, arithmetical) dynamical
system of interest, encode that map through a transformation of sufficiently-well-behaved
analytic objects. In general, the objects will be functions (holomorphic,
square-integrable, etc.) in some linear space\footnote{``Linear space'' $\overset{\textrm{def}}{=}$ ``Vector space''}
and the transformations will be linear operators on those spaces.
This procedure can be applied in significantly greater generality
(generalized Laplace transforms of indicator functions of subsets
of locally compact abelian groups), in this paper, we will confine
ourselves to the archetypical case of what happens over $\mathbb{Z}$.
To that end, \textbf{Set-series }are the means by which we will represent
$\mathbb{Z}$ and its subsets as analytic objects, \textbf{Hydra maps
}are the particular class of Collatz-type maps we will consider in
this paper, and \textbf{Permutation Operators }are the linear transformations
that will encode the action of our Hydra maps as transformations of
set-series, and of holomorphic functions in general.

\subsection{Set-Series on $\mathbb{Z}$}

The title of this section is (quite literally) only \emph{half} true;
it says ``Set-Series on $\mathbb{Z}$'', but it means ``Set-Series
on $\mathbb{N}_{0}$''. The reason for this is one of analytic necessity:
the sets under consideration must be bounded from at least one direction
(below or above) in order for their set-series to have a positive
radius of convergence. Without loss of generality, we can (and shall)
restrict our attentions to subsets of $\mathbb{N}_{0}$. The reasons
for this will be examined in §1.2.2, after we have gone over the definitions
and terminology needed for that discussion.

\subsubsection{Definitions and Basic Operations}

This subsection is just a list of definitions.

\vphantom{}\textbf{Definition 1 (Set-Series on $\mathbb{Z}$)}: Let
$V\subseteq\mathbb{N}_{0}$.

I. The \textbf{ordinary set-series }of $V$ is the function $\varsigma_{V}\left(z\right)$
defined by:
\begin{equation}
\varsigma_{V}\left(z\right)\overset{\textrm{def}}{=}\sum_{v\in V}z^{v}\label{eq:Def of ordinary set-series}
\end{equation}

\emph{Remark}: Since $V\subseteq\mathbb{N}_{0}$:
\[
\left|\varsigma_{V}\left(z\right)\right|\leq\sum_{v\in V}\left|z\right|^{v}\overset{V\subseteq\mathbb{N}_{0}}{\leq}\sum_{n=0}^{\infty}\left|z\right|^{n}=\frac{1}{1-\left|z\right|},\textrm{ }\forall z\in\mathbb{D}
\]
which shows that ordinary set-series are holomorphic functions on
$\mathbb{D}$; that is, $\varsigma_{V}\in\mathcal{A}\left(\mathbb{D}\right)$
for all $V\subseteq\mathbb{N}_{0}$.

\vphantom{}II. The \textbf{exponential set-series }of $V$ is the
function $\sigma_{V}\left(z\right)$ defined by:
\begin{equation}
\sigma_{V}\left(z\right)\overset{\textrm{def}}{=}\varsigma_{V}\left(e^{-z}\right)=\sum_{v\in V}e^{-vz}\label{eq:Def of exponential set-series}
\end{equation}

\emph{Remark}: Since: 
\[
e^{-z}\in\mathbb{D}\Leftrightarrow z\in\mathbb{H}_{+}
\]
it follows that:
\[
\left|\sigma_{V}\left(z\right)\right|=\left|\varsigma_{V}\left(e^{-z}\right)\right|\overset{V\subseteq\mathbb{N}_{0}}{\leq}\frac{1}{1-\left|e^{-z}\right|}=\frac{1}{1-e^{-\textrm{Re}\left(z\right)}},\textrm{ }\forall z\in\mathbb{H}_{+}
\]
which shows that exponential set-series are holomorphic functions
on the right half-plane $\mathbb{H}_{+}$.

\vphantom{}III. The \textbf{Fourier set-series }of $V$ is the function
$\psi_{V}\left(z\right)$ defined by:
\begin{equation}
\psi_{V}\left(z\right)\overset{\textrm{def}}{=}\varsigma_{V}\left(e^{2\pi iz}\right)=\sum_{v\in V}e^{2\pi ivz}\label{eq:Def of fourier set-series}
\end{equation}

\emph{Remark}: Since: 
\[
e^{2\pi iz}\in\mathbb{D}\Leftrightarrow z\in\mathbb{H}_{+i}
\]
it follows that:
\[
\left|\psi_{V}\left(z\right)\right|=\left|\varsigma_{V}\left(e^{2\pi iz}\right)\right|\overset{V\subseteq\mathbb{N}_{0}}{\leq}\frac{1}{1-\left|e^{2\pi iz}\right|}=\frac{1}{1-e^{-2\pi\textrm{Im}\left(z\right)}},\textrm{ }\forall z\in\mathbb{H}_{+i}
\]
which shows that Fourier set-series are holomorphic functions on the
upper half-plane $\mathbb{H}_{+i}$.

\vphantom{}IV. The \textbf{Dirichlet-set-series / Zeta-set-series
}of $V$ is the function $\zeta_{V}\left(s\right)$ defined by:
\begin{equation}
\zeta_{V}\left(s\right)\overset{\textrm{def}}{=}\sum_{v\in V}\frac{1}{v^{s}}\label{eq:Def of fourier set-series-1}
\end{equation}
where $s$ is a complex variable; $\zeta_{V}\left(s\right)$ is always
uniformly convergent for $s\in1+\mathbb{H}_{+}$.

V. When working with a series of a given type, we will use the term
``\textbf{digital}'' (ex: digital polynomial, digital function,
digital coefficients) to indicate that the coefficients of the function
or linear combination in question take values in the set $\left\{ 0,1\right\} \subseteq\mathbb{C}$.

\vphantom{}Unsurprisingly, the exact formulae for the operators and
operations we will consider later on depend on the particular type
of set-series being used. In this paper, we will focus primarily on
\emph{Ordinary }and \emph{Fourier }set-series. Fourier set-series
will be the type used from Chapter 2 onward, owing to the greater
simplicity of their permutation operators. For the legwork in this
chapter, however, we will mostly rely on ordinary set-series, which
are better suited for that sort of busy work. As such, the reader
should keep in mind that the results stated for ordinary set-series
also apply to every other type of set-series \emph{mutatis mutandis}.
Fortunately, our definitions will, for the most part, be independent
of the particular type of set-series, and work the same way for all.

\vphantom{}Almost everything one can do to a set, one can do to a
set-series on $\mathbb{Z}$.

\vphantom{}\textbf{Definition} \textbf{2}: Let $V,W\subseteq\mathbb{N}_{0}$.

I. Regardless of type, the set-series for the empty set ($\varnothing$)
is defined to be the constant function $0$.

II. We say that $\varsigma_{V}$ is \textbf{present }in $\varsigma_{W}$,
written $\varsigma_{V}\subseteq\varsigma_{W}$, exactly when $V\subseteq W$.
Likewise, when $\varsigma_{V}$ is not present in $\varsigma_{W}$,
we write $\varsigma_{V}\nsubseteq\varsigma_{W}$.

\emph{Examples}: $z^{2}\subseteq\varsigma_{\mathbb{N}_{1}}\left(z\right)$;
$1+e^{2\pi iz}\nsubseteq\frac{e^{2\pi iz}}{1-e^{2iz}}$, etc.

III. We say that $\varsigma_{V}$ is \textbf{disjoint }from $\varsigma_{W}$
whenever $V\cap W=\varnothing$.

\emph{Examples}: $\varsigma_{2^{\mathbb{N}_{0}}}\left(z\right)=z+z^{2}+z^{4}+z^{8}+\cdots$
and $\varsigma_{3\mathbb{N}_{0}}\left(z\right)=1+z^{3}+z^{6}+z^{9}+\cdots$
are disjoint.

More generally, we will say two linear combinations of set-series
are \textbf{disjoint }from one-another whenever they are $\mathbb{C}$-linearly
independent of one another.

IV. The \textbf{intersection }of $\varsigma_{V}$ and $\varsigma_{W}$,
denoted $\varsigma_{V}\cap\varsigma_{W}$, is defined by:
\[
\left(\varsigma_{V}\cap\varsigma_{W}\right)\left(z\right)\overset{\textrm{def}}{=}\varsigma_{V\cap W}\left(z\right)
\]
The symbol $\cap$ can be realized as a commutative binary operator.
In the case of ordinary set-series, $\cap$ corresponds to the \textbf{Hadamard
Product }of ordinary generating functions:
\[
\left(\varsigma_{V}\cap\varsigma_{W}\right)\left(z\right)=\frac{1}{2\pi i}\oint_{\partial\mathbb{D}}\varsigma_{V}\left(s\right)\varsigma_{W}\left(\frac{z}{s}\right)\frac{ds}{s},\textrm{ }\forall\left|z\right|<1
\]
There are corresponding versions of $\cap$ and its integral formula
for exponential, Fourier, and dirichlet set-series.

V. The \textbf{union }of $\varsigma_{V}$ and $\varsigma_{W}$, denoted
$\varsigma_{V}\cup\varsigma_{W}$, is defined by:
\[
\left(\varsigma_{V}\cup\varsigma_{W}\right)\left(z\right)\overset{\textrm{def}}{=}\varsigma_{V\cup W}\left(z\right)
\]

\vphantom{}\textbf{Proposition} \textbf{1} \textbf{(Properties of
Digital Set-Series and Function Intersections)}: Let $V,W\subseteq\mathbb{N}_{0}$,
and let $\varsigma_{V}$ and $\varsigma_{W}$ be their respective
set-series.

I. (Inclusion-Exclusion formula) $\varsigma_{V}\cup\varsigma_{W}=\varsigma_{V}+\varsigma_{W}-\left(\varsigma_{V}\cap\varsigma_{W}\right)$

II. $\varsigma_{V}\subseteq\varsigma_{W}$ if and only if $V\subseteq W$.

III. $\varsigma_{W\backslash V}=\varsigma_{W}-\varsigma_{V}$ if and
only if $V\subseteq W$.

\vphantom{}Proof:

I. There are two cases.

First, suppose $V\cap W=\varnothing$. Then, $\varsigma_{V\cap W}\left(z\right)=\varsigma_{\varnothing}\left(z\right)=0$.
As such:
\[
\varsigma_{V}\left(z\right)+\varsigma_{W}\left(z\right)-\varsigma_{V\cap W}\left(z\right)=\varsigma_{V}\left(z\right)+\varsigma_{W}\left(z\right)=\sum_{v\in V}z^{v}+\sum_{w\in W}z^{w}=\sum_{n\in V\cup W}z^{n}=\varsigma_{V\cup W}\left(z\right)
\]
as desired.

Second, supposing that $V\cap W\neq\varnothing$, let $S=V\cap W$,
and let $V^{\prime}=V\backslash S$ and $W^{\prime}=W\backslash S$.
Then $V^{\prime}$ and $W^{\prime}$ are disjoint, as are $V^{\prime}\cup W^{\prime}$
and $S$. Consequently, by the first case:
\begin{eqnarray*}
\varsigma_{V\cup W} & = & \varsigma_{\left(V^{\prime}\cup W^{\prime}\right)\cup S}\\
\left(\left(V^{\prime}\cup W^{\prime}\right)\cap S=\varnothing\right); & = & \varsigma_{V^{\prime}\cup W^{\prime}}+\varsigma_{S}\\
\left(V^{\prime}\cap W^{\prime}=\varnothing\right); & = & \varsigma_{V^{\prime}}+\varsigma_{W^{\prime}}+\varsigma_{S}\\
 & = & \left(\varsigma_{V^{\prime}}+\varsigma_{S}\right)+\left(\varsigma_{W^{\prime}}+\varsigma_{S}\right)-\varsigma_{S}\\
\left(V^{\prime}\cap S=W^{\prime}\cap S=\varnothing\right); & = & \varsigma_{V^{\prime}\cap S}+\varsigma_{W^{\prime}\cap S}-\varsigma_{S}\\
 & = & \varsigma_{V}+\varsigma_{W}-\varsigma_{V\cap W}
\end{eqnarray*}
as desired. $\checkmark$

\vphantom{}II. Let $V\subseteq W$. Then: 
\[
\varsigma_{W}\left(z\right)=\sum_{w\in W}z^{w}=\sum_{w\in W\backslash V}z^{w}+\sum_{v\in V}z^{v}=\sum_{w\in W\backslash V}z^{w}+\varsigma_{V}\left(z\right)
\]
Since $\varsigma_{W}\left(z\right)-\varsigma_{V}\left(z\right)$ is
the digital function $\sum_{w\in W\backslash V}z^{w}$, $\varsigma_{V}\left(z\right)\subseteq\varsigma_{W}\left(z\right)$.
For the other direction, let $\varsigma_{V}\left(z\right)\subseteq\varsigma_{W}\left(z\right)$.
Then, $\varsigma_{W}\left(z\right)-\varsigma_{V}\left(z\right)$ is
a digital function. Consequently, each for each $v\in V$, $z^{v}\subseteq\varsigma_{V}\left(z\right)$
forces $z^{v}\subseteq\varsigma_{W}\left(z\right)$, in order to ensure
that there will be no negative coefficients in $\varsigma_{W}\left(z\right)-\varsigma_{V}\left(z\right)$.
Since $z^{v}\subseteq\varsigma_{W}\left(z\right)$ forces $v\in W$,
it then follows that $V\subseteq W$. $\checkmark$

\vphantom{}III. $V\subseteq W$ and (II) immediately imply $\varsigma_{W}-\varsigma_{V}=\varsigma_{W\backslash V}$.
On the other hand, if $\varsigma_{W}-\varsigma_{V}=\varsigma_{W\backslash V}$,
the digitality of $\varsigma_{W\backslash V}$ forces $V\subseteq W$,
else $\varsigma_{W}-\varsigma_{V}$ will have terms with negative
coefficients. $\checkmark$

Q.E.D.

\vphantom{}$\mathbb{Z}$ is more than just a set, however: it is
an abelian group under addition, as well as a commutative ring under
addition and multiplication. These operations manifest as operators
on set-series, just like the set-theoretic operations described above.
Of greatest significance will be the action on $\mathbb{Z}$ of its
quotient subgroup---the cyclic groups $\left(\mathbb{Z}/a\mathbb{Z},+\right)$,
where $a\in\mathbb{N}_{1}$. These operators are well known in combinatorics
(where they are referred to as ``series multisection''), as well
as in harmonic analysis, where they are instances of the action of
a compact abelian group (in this case, the circle group) on appropriate
formulated function spaces (Tao. 2009). Since this more general aspect
of series multisection will increase in complexity in future generalizations,
it seems reasonable to give them a name of their own; we call them
\textbf{action decomposition operators}.

\vphantom{}\textbf{Definition 3 }(Action Decomposition Operators
for $\mathbb{Z}$): Let $a,b\in\mathbb{Z}$ be integers satisfying
$0\leq b\leq a-1$.

I. For every $V\subseteq\mathbb{N}_{0}$, we write $V_{a,b}$ to denote
the subset:
\[
V_{a,b}\overset{\textrm{def}}{=}V\cap\left(a\mathbb{Z}+b\right)
\]
We call $a$ the \textbf{modulus }or \textbf{period }of $V_{a,b}$,
and $b$ the \textbf{residue }of $V_{a,b}$. We call $V_{a,b}$ the
\textbf{cross-section of $V$ congruent to $b$ mod $a$}.

Note that for each $a\in\mathbb{N}_{1}$, The cross-sections $\left\{ V_{a,j}\right\} _{j\in\left\{ 0,\ldots,a-1\right\} }$
form a partition of $V$; that is, they are pair-wise disjoint with
respect to $j$, and satisfy:
\[
V=\bigcup_{j=0}^{a-1}V_{a,j}
\]

II. Define the \textbf{ordinary }$\left(a,b\right)$\textbf{-decomposition
operator} $\varpi_{a,b}:\mathcal{A}\left(\mathbb{D}\right)\rightarrow\mathcal{A}\left(\mathbb{D}\right)$
(a.k.a ``series multisection operator'') by the formula:
\begin{equation}
\varpi_{a,b}\left\{ \varsigma\right\} \left(z\right)\overset{\textrm{def}}{=}\frac{1}{a}\sum_{k=0}^{a-1}\xi_{a}^{-bk}\varsigma\left(\xi_{a}^{k}z\right),\textrm{ }\forall\varsigma\in\mathcal{A}\left(\mathbb{D}\right)\label{eq:ordinary (a,b)-decomposition operator definition}
\end{equation}
Then:
\begin{equation}
\varpi_{a,b}\left\{ \sum_{n=0}^{\infty}c_{n}z^{n}\right\} \left(z\right)=\sum_{n=0}^{\infty}c_{an+b}z^{an+b},\textrm{ }\forall\sum_{n=0}^{\infty}c_{n}z^{n}\in\mathcal{A}\left(\mathbb{D}\right)\label{eq:ordinary (a,b)-decomposition operator alternative formula}
\end{equation}
Moreover:
\begin{equation}
\varsigma\left(z\right)=\sum_{j=0}^{a-1}\varpi_{a,j}\left\{ \varsigma\right\} \left(z\right),\textrm{ }\forall z\in\mathbb{D},\textrm{ }\forall\varsigma\in\mathcal{A}\left(\mathbb{D}\right)\label{eq:ordinary (a,b)-decomposition}
\end{equation}
where the $\varpi_{a,j}\left\{ \varsigma\right\} $s are pair-wise
disjoint with respect to $j$.

Furthermore:
\begin{equation}
\varpi_{a,j}\left\{ \varsigma\right\} \left(\xi_{a}z\right)=\xi_{a}^{j}\varpi_{a,j}\left\{ \varsigma\right\} \left(z\right),\textrm{ }\forall j\in\mathbb{Z},\textrm{ }\forall z\in\mathbb{D},\textrm{ }\forall\varsigma\in\mathcal{A}\left(\mathbb{D}\right)\label{eq:ordianary (a,b)-decomposition operator automorphy}
\end{equation}

III. Define the \textbf{Fourier }$\left(a,b\right)$\textbf{-decomposition
operator} $\tilde{\varpi}{}_{a,b}:\mathcal{A}\left(\mathbb{H}_{+i}\right)\rightarrow\mathcal{A}\left(\mathbb{H}_{+i}\right)$
by the formula:
\begin{equation}
\tilde{\varpi}_{a,b}\left\{ \psi\right\} \left(z\right)\overset{\textrm{def}}{=}\frac{1}{a}\sum_{k=0}^{a-1}\xi_{a}^{-bk}\psi\left(z+\frac{k}{a}\right),\textrm{ }\forall\psi\in\mathcal{A}\left(\mathbb{H}_{+i}\right)\label{eq:fourier (a,b)-decomposition operator definition}
\end{equation}
Then:
\begin{equation}
\tilde{\varpi}_{a,b}\left\{ \sum_{n=0}^{\infty}c_{n}e^{2n\pi iz}\right\} \left(z\right)=\sum_{n=0}^{\infty}c_{an+b}e^{2\left(an+b\right)\pi iz},\textrm{ }\forall\sum_{n=0}^{\infty}c_{n}e^{2n\pi iz}\in\mathcal{A}\left(\mathbb{H}_{+i}\right)\label{eq:fourier (a,b)-decomposition operator alternative formula}
\end{equation}
Moreover:
\begin{equation}
\psi\left(z\right)=\sum_{j=0}^{a-1}\tilde{\varpi}_{a,j}\left\{ \psi\right\} \left(z\right),\textrm{ }\forall z\in\mathbb{H}_{+i},\textrm{ }\forall\psi\in\mathcal{A}\left(\mathbb{D}\right)\label{eq:fourier (a,b)-decomposition}
\end{equation}
where the $\tilde{\varpi}_{a,j}\left\{ \psi\right\} $s are pair-wise
disjoint with respect to $j$.
\begin{equation}
\tilde{\varpi}_{a,j}\left\{ \psi\right\} \left(z+\frac{1}{a}\right)=\xi_{a}^{j}\tilde{\varpi}_{a,j}\left\{ \psi\right\} \left(z\right),\textrm{ }\forall j\in\mathbb{Z},\textrm{ }\forall z\in\mathbb{H}_{+i},\textrm{ }\forall\psi\in\mathcal{A}\left(\mathbb{H}_{+i}\right)\label{eq:fourier (a,b)-decomposition operator automorphy-1}
\end{equation}

\subsubsection{Analytic Aspects of Set-Series on $\mathbb{Z}$}

Given an $V\subseteq\mathbb{N}_{0}$, a set-series will encode the
properties of $V$ in various ways. For our purposes, one of the most
important of these is the behavior of the set-series as we approach
the boundary of the domain of region of the series' convergence.

\vphantom{}\textbf{Definition} \textbf{4}: Let $D$ be an open, connected,
non-empty subset of $\mathbb{C}$ (with boundary $\partial D$), and
let $f:D\rightarrow\mathbb{C}$ be a holomorphic function $f$ is
said to have a \textbf{natural boundary }on $\partial D$ (equivalently,
$D$ is said to be a \textbf{domain of holomorphy} of $f$) if there
is no open set $U\subseteq\mathbb{C}$, containing $D$ as a proper
subset, such that $f$ admits an analytic continuation to $U$.

Slightly modifying the definition given by Breuer \& Simon (Breuer
\& Simon. 2011), we say that $\partial D$ is a \textbf{strong natural
boundary} of $f$ whenever:
\[
\sup_{z\in U}\left|f\left(z\right)\right|=\infty
\]
for every open, connected set $U\subseteq\mathbb{C}$ for which $U\cap\partial D\neq\varnothing$.
(That is to say, $f$ blows up every point of $\partial D$).

\vphantom{}One of the reasons set-series are useful at all comes
from the fact that their digitality places places rather stringent
conditions on the possible forms that they can take, thanks to a beautiful
from a century ago, due to Gabor Szeg\H{o}.

\vphantom{}\textbf{Theorem 1: Szeg\H{o}'s Theorem}: (Remmert. 1991)
(Breuer \& Simon. 2011)): Let $C$ be a finite subset of $\mathbb{C}$.
and let $f\left(z\right)=\sum_{n=0}^{\infty}c_{n}z^{n}$ be a power
series so that $c_{n}\in C$ for all $n\in\mathbb{N}_{0}$. Then,
exactly one of the following occurs:

I. $f\left(z\right)$ is a polynomial.

II. $f\left(z\right)$ is a non-polynomial rational function $f\left(z\right)=\frac{P\left(z\right)}{1-z^{\alpha}}$,
where $\alpha$ is a positive integer, and $P\left(z\right)$ is a
polynomial in $z$ whose coefficients are in $C$.

III. $f\left(z\right)$ is a transcendental function with a strong
natural boundary on $\partial\mathbb{D}$.

\vphantom{}Using this theorem, we can classify set-series into three
types.

\vphantom{}\textbf{Corollary 1 (Set-Series Classification Theorem)}:\textbf{
}Let $V\subseteq\mathbb{N}_{0}$. Then, exactly one of the following
occurs:

I. $\varsigma_{V}\left(z\right)$ is a digital polynomial. Moreover,
this occurs if and only if $V$ is finite. In particular, $\varsigma_{V}\left(z\right)$
is identically zero if and only if $V$ is empty.

II. $\varsigma_{V}\left(z\right)$ is a non-polynomial digital rational
function, which can be written in the form:
\[
\varsigma_{V}\left(z\right)=Q\left(z\right)+\frac{z^{\alpha\nu}P\left(z\right)}{1-z^{\alpha}}
\]
where $\alpha\in\mathbb{N}_{1}$, $\nu\in\mathbb{N}_{0}$, and where
$P$ and $Q$ are digital polynomials with $\deg Q\leq z^{\alpha\nu-1}$
and $\deg P\leq\alpha-1$, and where $P\left(0\right)\geq1$ (also,
$\gcd\left(P\left(z\right),1-z^{\alpha}\right)=1$).

Moreover, this occurs if and only if $V$ is the union of a finite
set (whose set-series is $Q\left(z\right)$) and a finite number of
IAPs of period $\alpha$ (the set-series for these IAPs being the
function $\frac{z^{\alpha\nu}P\left(z\right)}{1-z^{\alpha}}$).

III. $\varsigma_{V}\left(z\right)$ is a transcendental digital function
for which the unit circle is a strong natural boundary.

\vphantom{}Parallel to this classification theorem, we have the following
terminology for sets in $\mathbb{N}_{0}$:

\vphantom{}\textbf{Definition 5}: Let $V\subseteq\mathbb{N}_{0}$
be an infinite set.

I. We say $V$ is \textbf{rational }whenever $\varsigma_{V}\left(z\right)$
is a rational function (i.e., if and only if $V$ is the union of
a finite set and a finite number of IAPs of period $\alpha$).

II. We say $V$ is \textbf{irrational }whenever $\varsigma_{V}\left(z\right)$
is not a rational function (i.e., if and only if $\varsigma_{V}\left(z\right)$
is a transcendental function with a natural boundary on $\partial\mathbb{D}$).

\vphantom{}As a rule of thumb, the boundary behavior of power series
lead to the following correlations:

\textbullet{} The \emph{location }of set-series' singularity to the
\emph{arithmetic }properties of the set (that is, its cross-sections);

\textbullet{} The \emph{strength }of the singularity (that is, the
asymptotic behavior of the set-series as its input tends toward the
singularity) and the \emph{distribution} of the set (its density).

\vphantom{}\textbf{Corollary 1 }highlights the first of these correspondences.
Understanding the second correspondence, on the other hand, is the
central concern of what, in analytic number theory, are known as \textbf{Tauberian
Theorems}---results relate the asymptotic properties of a function
to those of the coefficients in that function's defining series. For
our purposes, however, we only need one such result.

\vphantom{}\textbf{Theorem 2: The Hardy-Littlewood Tauberian Theorem
(HLTT)} (Hardy. \emph{Divergent Series}. 154-155\footnote{Hardy states a series of related results; I'm combining things from
several places in order to get the version of the theorem given here.}) Let $\left\{ c_{n}\right\} _{n\in\mathbb{N}_{0}}\subseteq\mathbb{C}$
be a bounded sequence of complex numbers. Then:

I. 
\[
\lim_{x\uparrow1}\left(1-x\right)\sum_{n=0}^{\infty}c_{n}x^{n}=C
\]
implies:
\[
\lim_{N\rightarrow\infty}\frac{1}{N}\sum_{n=0}^{N}c_{n}=C
\]

II. 
\[
\lim_{x\downarrow0}x\sum_{n=0}^{\infty}c_{n}e^{-nx}=C
\]
implies:
\[
\lim_{N\rightarrow\infty}\frac{1}{N}\sum_{n=0}^{N}c_{n}=C
\]

\vphantom{}There is much to be said about this theorem's applications
to set-series, both now, and in Chapter 2, where they will be used
to define the concept of a ``virtual pole''. First, however, some
definitions:

\vphantom{}\textbf{Definition 6: }Let $V\subseteq\mathbb{N}_{0}$.

I. Given any $N\in\mathbb{N}_{0}$, we write: 
\[
V\left(N\right)\overset{\textrm{def}}{=}\left\{ v\in V:v\leq N\right\} 
\]
In a similar vein, $\left|V\left(N\right)\right|$ ($V$'s \textbf{counting
function}) denotes the number of elements of $V\left(N\right)$. In
a slight abuse of terminology, when we say something like ``the growth
of $V$'' or ``the asymptotics of $V$'', we mean the growth and
asymptotics of $\left|V\left(N\right)\right|$ as $N\rightarrow\infty$,
respectively.

II. Define the \textbf{upper density }(resp. \textbf{lower density})
of $V$ by the limits:
\[
\overline{d}\left(V\right)\overset{\textrm{def}}{=}\limsup_{N\rightarrow\infty}\frac{\left|V\left(N\right)\right|}{N}
\]
and:
\[
\underline{d}\left(V\right)\overset{\textrm{def}}{=}\liminf_{N\rightarrow\infty}\frac{\left|V\left(N\right)\right|}{N}
\]
respectively.

\emph{Remark}s: Enumerating the elements of $V$ in increasing order
as $\left\{ v_{n}\right\} _{n\in\mathbb{N}_{1}}$, one can also write:
\[
\overline{d}\left(V\right)=\limsup_{n\rightarrow\infty}\frac{n}{v_{n}}
\]
\[
\underline{d}\left(V\right)=\liminf_{n\rightarrow\infty}\frac{n}{v_{n}}
\]

iii. $\overline{d}\left(\mathbb{N}_{0}\backslash V\right)=1-\overline{d}\left(V\right)$,
and $\underline{d}\left(\mathbb{N}_{0}\backslash V\right)=1-\underline{d}\left(V\right)$.

iv. $W\subseteq V$ implies $\overline{d}\left(W\right)\leq\overline{d}\left(V\right)$
and $\underline{d}\left(W\right)\leq\underline{d}\left(V\right)$.

III. The \textbf{(natural) density }of $V$ is defined by the limit:
\[
d\left(V\right)\overset{\textrm{def}}{=}\lim_{N\rightarrow\infty}\frac{\left|V\left(N\right)\right|}{N}
\]
provided that the limit exists.

\emph{Remarks}:

i. Not every $V\subseteq\mathbb{N}_{0}$ has a well-defined natural
density!

ii. $V$ has a well-defined natural density if and only if $\mathbb{N}_{0}\backslash V$
has a well-defined natural density, in which case $d\left(\mathbb{N}_{0}\backslash V\right)=1-d\left(V\right)$.

\vphantom{}Using this terminology, we can apply the \textbf{HLTT
}to get a few elementary conclusions about the relationships between
singularities of a set-series and the growth of the represented sets.
In the following \textbf{Corollary}, we work with Fourier set-series,
and for two reasons: (1) because we will need part (III) near the
end of the paper in Chapter 4; (2), because parts (I) and (II) will
serve to motivate our work with Fourier set-series in Chapters 2 and
3.

\vphantom{}\textbf{Corollary 2: Minor Set-Series Corollaries of the
HLTT}: Let $V\subseteq\mathbb{N}_{0}$.

I. The limit $\lim_{y\downarrow0}y\psi_{V}\left(iy\right)$ exists
if and only if $V$ has a well-defined natural density, in which case:
\[
\lim_{y\downarrow0}y\psi_{V}\left(iy\right)=d\left(V\right)
\]
More generally, given any rational number $t=\frac{\beta}{\alpha}$
written in irreducible form ($\gcd\left(\alpha,\beta\right)=1$):
\[
\lim_{y\downarrow0}y\psi_{V}\left(\frac{\beta}{\alpha}+iy\right)=\frac{1}{2\pi}\lim_{N\rightarrow\infty}\frac{1}{N}\sum_{k=0}^{\alpha-1}\xi_{\alpha}^{k\beta}\left|V_{\alpha,k}\left(N\right)\right|=\frac{1}{2\pi}\sum_{k=0}^{\alpha-1}d\left(V_{\alpha,k}\right)\xi_{\alpha}^{k\beta}
\]
provided that the limits/densities exist.

II. Given any rational number $t=\frac{\beta}{\alpha}$ written in
irreducible form, if the limit $\lim_{y\downarrow0}y\psi_{V}\left(\frac{\beta}{\alpha}+iy\right)$
fails to exist, then $y\psi_{V}\left(\frac{\beta}{\alpha}+iy\right)$
fails to be of bounded variation as a function of $y$ on any interval
of the form $\left[0,\epsilon\right]$, where $\epsilon\in\mathbb{R}>0$.

Proof:

I. Let $t=\frac{\beta}{\alpha}\in\mathbb{Q}$ be in irreducible form.
Then, since $V$ is partitioned by the sets $V_{\alpha,0},V_{\alpha,1},\ldots,V_{\alpha,\alpha-1}$,
we can write:
\begin{eqnarray*}
\lim_{y\downarrow0}y\psi_{V}\left(t+iy\right) & = & \lim_{y\downarrow0}y\sum_{v\in V}e^{2\pi iv\left(t+iy\right)}\\
 & = & \lim_{y\downarrow0}y\sum_{k=0}^{\alpha-1}\sum_{v\in V_{\alpha,k}}e^{2\pi iv\left(t+iy\right)}\\
\left(t=\frac{\beta}{\alpha}\right); & = & \lim_{y\downarrow0}y\sum_{k=0}^{\alpha-1}\sum_{v\in V_{\alpha,k}}\xi_{\alpha}^{v\beta}e^{-2\pi vy}\\
\left(v\overset{\alpha}{\equiv}k\Leftrightarrow v\in V_{\alpha,k}\right); & = & \sum_{k=0}^{\alpha-1}\xi_{\alpha}^{k\beta}\lim_{y\downarrow0}y\sum_{v\in V_{\alpha,k}}e^{-2\pi vy}\\
 & = & \frac{1}{2\pi}\sum_{k=0}^{\alpha-1}\xi_{\alpha}^{k\beta}\lim_{y\downarrow0}2\pi y\sum_{n=0}^{\infty}\mathbf{1}_{V_{\alpha,k}}\left(n\right)e^{-2\pi ny}\\
\left(\textrm{let }x=2\pi y\right); & = & \frac{1}{2\pi}\sum_{k=0}^{\alpha-1}\xi_{\alpha}^{k\beta}\underbrace{\lim_{x\downarrow0}x\sum_{n=0}^{\infty}\mathbf{1}_{V_{\alpha,k}}\left(n\right)e^{-nx}}_{\textrm{apply \textbf{HLTT}}}\\
 & = & \frac{1}{2\pi}\sum_{k=0}^{\alpha-1}\xi_{\alpha}^{k\beta}\lim_{N\rightarrow\infty}\frac{1}{N}\sum_{n=0}^{N-1}\mathbf{1}_{V_{\alpha,k}}\left(n\right)\\
 & = & \frac{1}{2\pi}\sum_{k=0}^{\alpha-1}\xi_{\alpha}^{k\beta}\lim_{N\rightarrow\infty}\frac{1}{N}\left|V_{\alpha,k}\left(N\right)\right|\\
 & = & \frac{1}{2\pi}\lim_{N\rightarrow\infty}\frac{1}{N}\sum_{k=0}^{\alpha-1}\xi_{\alpha}^{k\beta}\left|V_{\alpha,k}\left(N\right)\right|
\end{eqnarray*}
Additionally, if all of the cross-sections of $V$ mod $\alpha$ have
well-defined natural densities, we can write:
\[
\lim_{y\downarrow0}y\psi_{V}\left(t+iy\right)=\frac{1}{2\pi}\sum_{k=0}^{\alpha-1}\xi_{\alpha}^{k\beta}\underbrace{\lim_{N\rightarrow\infty}\frac{1}{N}\left|V_{\alpha,k}\left(N\right)\right|}_{d\left(V_{\alpha,k}\right)}=\frac{1}{2\pi}\sum_{k=0}^{\alpha-1}d\left(V_{\alpha,k}\right)\xi_{\alpha}^{k\beta}
\]
as desired.

II. Let $t=\frac{\beta}{\alpha}\in\mathbb{Q}$ be in irreducible form.
We prove the contrapositive: ``if there is an $\epsilon\in\mathbb{R}>0$
so that $y\psi_{V}\left(t+iy\right)$ is of bounded variation as a
function of $y\in\left[0,\epsilon\right]$, then $\lim_{y\downarrow0}y\psi_{V}\left(t+iy\right)$
exists''. The proof is a standard result: if there is such an $\epsilon$
so that $y\psi_{V}\left(t+iy\right)$ is a function of bounded variation
on $\left[0,\epsilon\right]$, the limit as $y$ decreases to $0$
\emph{must} exist.

Q.E.D.

\pagebreak{}

\subsection{Hydra Maps on $\mathbb{Z}$}

\subsubsection{Dynamical Systems - Terms \& Definitions}

To begin, let us introduce some terminology of dynamical systems.

\vphantom{}\textbf{Definition 7}:\textbf{ }For technical reasons
(to be discussed in detail later on), we will restrict our attention
to\textbf{ surjective (discrete) maps}: surjective functions $H:\mathbb{N}_{0}\rightarrow\mathbb{N}_{0}$
($H$ stands for \emph{hydra}). Given such a map, we will use the
following terminology; excepting (V), which is the author's coinage,
we follow (Devaney. 2003):

\vphantom{}I. We say a non-empty set $V\subseteq\mathbb{N}_{0}$
is an \textbf{invariant (set) }of $H$ (or is simply said to be ``invariant'')
whenever both $H\left(V\right)=V$ and $H^{-1}\left(V\right)=V$ (here,
$H^{-1}\left(V\right)$ is the pre-image of $V$ under $H$).

\vphantom{}II. We say that an integer $v$ is a \textbf{periodic
point of }\textbf{\emph{H}} (equivalently, that $v$ is \textbf{periodic
(under }\textbf{\emph{H}}\textbf{)})\textbf{\emph{ }}if there exists
a positive integer $N$ such that $H^{\circ N}\left(v\right)=v$.

\vphantom{}III. Given an integer $v$, we denote the (\textbf{forward)}
\textbf{orbit }of $v$ under $H$ by:

\[
O_{H}\left(v\right)\overset{\textrm{def}}{=}\left\{ H^{\circ k}\left(v\right)\right\} _{k\in\mathbb{N}_{0}}
\]
Note that, for all $v$, $H\left(O_{H}\left(v\right)\right)\subseteq O_{H}\left(v\right)$.

\vphantom{}IV. We say that an integer $v$ is \textbf{pre-periodic
(under }\textbf{\emph{H}}\textbf{)}\textbf{\emph{ }}\textbf{with period
}\textbf{\emph{N}} if $O_{H}\left(v\right)$ contains a periodic point
of $H$, but $v\notin O_{H}\left(v\right)$ (i.e., $v$ is not a periodic
point, but there is some $n\geq1$ such that $H^{\circ n}\left(v\right)$
\emph{is }a periodic point.)

\vphantom{}V. A \textbf{cycle of }\textbf{\emph{H }}is the orbit
of a periodic point of $H$. Moreover, in a slight abuse of terminology,
letting $\Omega$ denote a cycle of $H$, we say that $\Omega$ is
\textbf{closed }if $H^{-1}\left(\Omega\right)=\Omega$. Alternatively,
we say that $\Omega$ is \textbf{open }if $H^{-1}\left(\Omega\right)\neq\Omega$.

\vphantom{}\emph{Remark}: Note that open and closed cycles are the
discrete analogues of the notions of stable and unstable orbits/trajectories
for continuous dynamical systems.

\vphantom{}\textbf{Definition 8: }Let $H:\mathbb{N}_{0}\rightarrow\mathbb{N}_{0}$
be an arbitrary map. Then define the following relation $\sim$ on
$\mathbb{N}_{0}$: two integers $v,w$ satisfy $v\sim w$ if and only
if there exist integers $j,k\geq0$ such that either $H^{\circ j}\left(v\right)=H^{\circ k}\left(w\right)$
(that is to say, $H$ eventually iterates $v$ and $w$ to a common
point.) The relation $\sim$ is called the \textbf{(total) orbit equivalence
relation} \textbf{for }\textbf{\emph{H}}. The equivalence classes
of $\mathbb{N}_{0}$ under $\sim$ are the \textbf{(total) orbit classes
of }\textbf{\emph{H}}.

\emph{Remark:} \emph{We will denote the orbit classes of a map $H$
by $V_{0},V_{1},\ldots,V_{j},\ldots$ Throughout this paper, whenever
$V_{j}$ is written, it refers to the orbit class of a map under consideration.
}Additonally, we say an orbit class of $H$ is \textbf{divergent }if
it contains an element whose forward orbit under $H$ is unbounded.

\vphantom{}The reader should note that the definition of orbit classes
does not require the surjectivity condition we imposed upon our maps.
The following proposition shows why the surjectivity is important,
and why we restrict our attention to \emph{surjective maps }on $\mathbb{N}_{0}$

\vphantom{}\textbf{Proposition} \textbf{2}: Let $H:\mathbb{N}_{0}\rightarrow\mathbb{N}_{0}$
be an arbitrary map, let $J$ be the countable indexing set for the
orbit classes $V_{j}$ of $H$.

I. Any invariant $V$ of $H$ is either an individual orbit class
of $H$ (i.e., $V=V_{j}$ for some $j$) or a union thereof ($V=\bigcup_{j\in J^{\prime}}V_{j}$
for some $J^{\prime}\subseteq J$).

II. $V_{j}$ is an invariant of $H$ for all $j$ if and only if surjective.

\vphantom{}Proof:

\vphantom{}I. Let $V$ be an invariant of $H$. Since the $V_{j}$s
are equivalence classes, they form a partition of $\mathbb{N}_{0}$.
As such, $V=\bigcup_{j\in J}V\cap V_{j}$. Since $V$ is invariant,
$V$ is non-empty, and so, there is a $j$ such that $V\cap V_{j}\neq\varnothing$.
Letting $v\in V\cap V_{j}$, the invariance of $V$ under $H$ then
guarantees that $w\in V$ for all integers $w$ such that $w\sim v$.
Thus, $V$ containing a single element of $V_{j}$ forces $V$ to
contain all of that $V_{j}$. Thus:
\[
V=\bigcup_{j\in J^{\prime}}V_{j}
\]
where the $J^{\prime}$ is the set of all $j\in J$ for which $V\cap V_{j}\neq\varnothing$.
$\checkmark$

\vphantom{}II. First, let us see how things might go wrong; let $H$
be non-surjective. Then, there is an integer $w$ such that $H\left(v\right)\neq w$
for any integer $v$. Let $V_{j}$ be the (unique) orbit class of
$H$ containing $w$. Then, $H\left(V_{j}\right)$ does not contain
$w$; if it did, there would be an integer $v$ in $V_{j}$ that $H$
sent to $w$. As such, $V_{j}$ is an orbit class of $H$ which is
not invariant. $\checkmark$

For the other direction, assuming $H$ is surjective, suppose by way
of contradiction that there is an orbit class $V_{j}$ which is \emph{not
}invariant. Since $H\left(V_{j}\right)\subseteq V_{j}$, it must be
that there is a $v\in V_{j}$ and a $w\in H^{-1}\left(\left\{ v\right\} \right)$
such that $w\notin V_{j}$; we need the assumed surjectivity of $H$
to guarantee the existence of this $w$. Since $w\notin V_{j}$ and
$v\in V_{j}$, $w\nsim v$ (they belong to different equivalence classes).
However, since $H\left(w\right)=v$, by definition, $w\sim v$---this
is a contradiction. Thus, all the $V_{j}$s must be invariant under
$H$ when $H$ is surjective. $\checkmark$

Q.E.D.

\vphantom{}Next, we consider how the invariant sets relate to the
dynamical behaviors of $H$.

\vphantom{}\textbf{Proposition} \textbf{3} \textbf{(Classification
of Orbit Classes of Surjective Maps)}: Let $H$ be a surjective map.
Then, for any given orbit class $V_{j}$ of $H$, exactly one of the
following statements is true:

I. $V_{j}$ contains exactly one cycle of $H$ (and, necessarily,
that cycle is either open or closed)

II. $V_{j}$ contains an integer $\lambda$ such that $O_{H}\left(\lambda\right)$
is unbounded. Moreover, in this case, the restriction of $H$ to $O_{H}\left(\lambda\right)$
is (topologically) conjugate to the translation $n\mapsto n+1$ on
$\mathbb{N}_{0}$; that is to say, there is a bijection $B:\mathbb{N}_{0}\rightarrow\mathbb{N}_{0}$
such that: 
\[
\left(B\circ H\mid_{O_{H}\left(\lambda\right)}\circ B^{-1}\right)\left(n\right)=n+1,\textrm{ }\forall n\in\mathbb{N}_{0}
\]

\vphantom{}Proof: Let $V_{j}$ be an orbit class of $H$. We prove
the proposition by showing first that $V_{j}$ must contain either
a cycle or an unbounded orbit, but not both. Then, we show that $V_{j}$
cannot contain two or more distinct (i.e., disjoint) cycles.

\vphantom{}\textbullet{} ($V_{j}$ contains either a cycle or an
unbounded orbit) If $V_{j}$ just so happens to contain an unbounded
orbit, we are done. So, suppose that $V_{j}$ does \emph{not} contain
an unbounded orbit; that is, suppose the orbit of \emph{every} element
of $V_{j}$ is bounded. Then, letting $v\in V$, $O_{H}\left(v\right)$
is therefore a bounded infinite sequence of integers, as such, it
must contain a convergent subsequence. 

Since the set of integers is discrete, the convergence of this subsequence
implies that there is a term of the subsequence, $H^{\circ m}\left(v\right)$,
such that $H^{\circ n}\left(H^{\circ m}\left(v\right)\right)=H^{\circ m}\left(v\right)$
for infinitely many $n\geq1$. Thus, $H^{\circ m}\left(v\right)$
is a periodic point of $H$ contained in $V_{j}$ and hence, $O_{H}\left(H^{\circ m}\left(v\right)\right)$
is a cycle contained in $V_{j}$. $\checkmark$

\vphantom{}\textbullet{} ($V_{j}$ cannot contain both an unbounded
orbit and a cycle) Suppose that $V_{j}$ contains both an unbounded
orbit and a cycle. Since $V_{j}$ contains a cycle, it contains a
periodic point $v$; likewise, since $V_{j}$ contains an unbounded
orbit, there is a $\lambda\in V_{j}$ for which $O_{H}\left(\lambda\right)$
is unbounded. Since $V_{j}$ is an orbit class, there are integers
$m,n\geq0$ such that $H^{\circ m}\left(v\right)=H^{\circ n}\left(\lambda\right)$.
Letting $w$ denote $H^{\circ m}\left(v\right)$, we then see that
$O_{H}\left(w\right)$ must be unbounded (since $w\in O_{H}\left(\lambda\right)$).
At the same time, though, since $w\in O_{H}\left(v\right)$, and since
$O_{H}\left(v\right)$ is a cycle, $O_{H}\left(w\right)$ is therefore
a cycle, and hence, contains only finitely many distinct integers.
Thus, $O_{H}\left(w\right)$ is an unbounded set containing only finitely
many distinct integers---but this is impossible. $\checkmark$

\vphantom{}\textbullet{} ($V_{j}$ can contain at most one cycle).
Suppose $V_{j}$ contains two disjoint cycles $\Omega$ and $\Omega^{\prime}$.
Then, there exist integers $v,v^{\prime}\in V_{j}$ such that $\Omega=O_{H}\left(v\right)$
and $\Omega^{\prime}=O_{H}\left(v^{\prime}\right)$. Since $v$ and
$v^{\prime}$ belong to the same orbit class of $H$, there are non-negative
integers $m,n$ such that $H^{\circ m}\left(v\right)=H^{\circ n}\left(v^{\prime}\right)$.
Let $w=H^{\circ m}\left(v\right)$. Then, $w$ is an element of both
$\Omega$ and $\Omega^{\prime}$. But this is impossible, we assumed
$\Omega\cap\Omega^{\prime}=\varnothing$. Consequently, $V_{j}$ can
contain at most one cycle. $\checkmark$

\vphantom{}Finally, to see the topological conjugacy, enumerate $O_{H}\left(\lambda\right)$
as $\lambda_{0},\lambda_{1},\ldots$ where $\lambda_{n}=H^{\circ n}\left(\lambda\right)$,
with $\lambda_{0}=\lambda$. Then, by construction: 
\[
H\left(\lambda_{n}\right)=\lambda_{n+1}
\]
Hence, the ``index-giving'' map: $B:\mathbb{N}_{0}\rightarrow\mathbb{N}_{0}$
defined by $B:\lambda_{n}\mapsto n$ is a bijection and:
\[
\left(B\circ H_{\mid O_{H}\left(\lambda\right)}\circ B^{-1}\right)\left(n\right)=n+1
\]
as desired. $\checkmark$

Q.E.D.

\vphantom{}Having thus covered all of the essential terminology for
our study of the dynamics of maps on $\mathbb{N}_{0}$, we can now
begin the task of translating this terminology into the language of
linear spaces of holomorphic functions.

\subsubsection{Hydra Maps on $\mathbb{Z}$ - Definitions, Examples, \& Conjectures}

As mentioned previously, the Collatz map is an example of what are
sometimes cumbersomely referred to as ``Residue-class-wise-affine
maps''. Studying generalizations of the Collatz map, in the hopes
of being able to discern broader trends that might be of use in attacking
the Collatz Conjecture, is an obvious and natural way to proceed.
Jeffery Lagarias, in his explanatory article (Lagarias. 2010) and
annotated bibliographies (Lagarias. 2011, 2012) on the Conjecture
and K.R. Matthews on his website (<http://www.numbertheory.org/3x+1/>)
both give wonderful surveys on these intrepid efforts. Of course,
different authors have used their own terminology in dealing with
and classifying Collatz-type maps. Given that many examples ($H_{3}$,
$H_{5}$, and so on) are often studied independently, this can complicate
matters, especially when---as is the case in this paper---attempts
are made to generalize useful arguments. For these reasons---and
for the admittedly irresistable allure of getting to set a new terminological
trend---the author has decided to give it a go of his own.

In this subsection, we will record these definitions, as well discuss
several examples, and the kinds of conjectures that tend to crop up
around them.

\vphantom{}\textbf{Definition 9 }(\textbf{Hydra Maps on $\mathbb{Z}$)}:\textbf{
}Let $\varrho\in\mathbb{N}_{2}$. A \textbf{hydra map }(on $\mathbb{Z}$)
of \textbf{order }\textbf{\emph{$\varrho$ }}(also called a \textbf{$\varrho$-hydra}
\textbf{map})\textbf{\emph{ }}is a \emph{surjective }map $H:\mathbb{N}_{0}\rightarrow\mathbb{N}_{0}$
of the form:
\begin{equation}
H\left(n\right)=\begin{cases}
\frac{a_{0}n+b_{0}}{d_{0}} & \textrm{if }n\overset{\varrho}{\equiv}0\\
\frac{a_{1}n+b_{1}}{d_{1}} & \textrm{if }n\overset{\varrho}{\equiv}1\\
\vdots & \vdots\\
\frac{a_{\varrho-1}n+b_{\varrho-1}}{d_{\varrho}} & \textrm{if }n\overset{\varrho}{\equiv}\varrho-1
\end{cases}\label{eq:Def of a Hydra Map on Z}
\end{equation}
where for each $j\in\left\{ 0,1,\ldots,\varrho-1\right\} $, $a_{j}$,
$b_{j}$, and $d_{j}$ are non-negative integers (with $a_{j},d_{j}>0$
and $\gcd\left(a_{j},d_{j}\right)=1$ for all $j$) such that $\frac{a_{j}n+b_{j}}{d_{j}}\in\mathbb{N}_{0}$
for all non-negative integers $n\overset{\varrho}{\equiv}j$. The
maps $n\mapsto\frac{a_{j}n+b_{j}}{d_{j}}$ are called the \textbf{branches
}of the hydra map.

\emph{Remark}:\emph{ }As a consequence of the definitions, observe
that for all $m\in\mathbb{N}_{0}$ and all $j\in\left\{ 0,\ldots,\varrho-1\right\} $:
\[
\underbrace{H\left(\varrho m+j\right)}_{\in\mathbb{N}_{0}}=\frac{a_{j}\left(\varrho m+j\right)+b_{j}}{d_{j}}=\frac{\varrho a_{j}m}{d_{j}}+\frac{ja_{j}+b_{j}}{d_{j}}=\frac{\varrho a_{j}}{d_{j}}m+\underbrace{H\left(j\right)}_{\in\mathbb{N}_{0}}
\]
Setting $m=1$, this then forces the quantity $\frac{\varrho a_{j}}{d_{j}}=H\left(\varrho m+j\right)-H\left(j\right)$
to be an integer for all $j$. This number occurs often enough in
computations that it is worth having its own symbol:
\begin{equation}
\mu_{j}\overset{\textrm{def}}{=}\frac{\varrho a_{j}}{d_{j}}\in\mathbb{N}_{1},\textrm{ }\forall j\in\left\{ 0,\ldots,\varrho-1\right\} \label{eq:Def of mu_j}
\end{equation}
In particular, since $a_{j}$ and $d_{j}$ are always co-prime, note
that $d_{j}$ must be a divisor of $\varrho$.

\vphantom{}Next, there are some flavors of hydra maps to be considered.

\vphantom{}\textbf{Definition 10 (Varieties of Hydra Maps)}: Let
$H:\mathbb{N}_{0}\rightarrow\mathbb{N}_{0}$ be a $\varrho$-hydra
map.

I. We say $H$ is \textbf{prime} if $\varrho\in\mathbb{P}$ and if
$a_{j}\in\left\{ 1\right\} \cup\mathbb{P}$ for all $j\in\left\{ 0,\ldots,\varrho-1\right\} $.

\emph{Remark}: If $H$ is prime, then, the fact that $d_{j}\mid\varrho$
for each $j$ forces $d_{j}\in\left\{ 1,\varrho\right\} $ for each
$j$.

II. We say $H$ is \textbf{regulated }if there exists a $j\in\left\{ 0,\ldots,\varrho-1\right\} $
(called a \textbf{regulated index})\textbf{ }for which $\mu_{j}=1$.

\emph{Remark}: Originally, the author worked with a notion of ``attraction'':
for all $n\overset{\varrho}{\equiv}0$, $H\left(n\right)=\frac{n}{d}$,
where $d\in\mathbb{N}_{2}$ is a divisor of $\varrho$, believing
that this was the condition necessary in order for the \textbf{Dreamcatcher
Theorem }to hold true. However, it turned out it sufficed that $H$
be regulated. The regulation condition can be shown to be equivalent
to the statement that the inverse of the $j$th branch of $H$ (where
$j$ is a regulated index) is an a measure-preserving bijection on:
\[
\widetilde{\mathbb{Z}}_{\bcancel{\varrho}}=\prod_{\begin{array}{c}
p\in\mathbb{P}\\
p\nmid\varrho
\end{array}}\mathbb{Z}_{p}
\]
where the measure (as will be discussed in Chapter 3); or, equivalently,
to a bijection of $\mathbb{Q}_{\bcancel{\varrho}}/\mathbb{Z}$.

\vphantom{}Before continuing to examples and conjectures, as promised,
there are several subtleties in the definition of a $\varrho$-hydra
map $H$ that must be discussed, as promised.

\vphantom{}\textbullet{} \emph{Why $\mathbb{N}_{0}$?}

As was seen in §1.2.1, for any $V\subseteq\mathbb{N}_{0}$, the region
of convergence of $\varsigma_{V}\left(z\right)$ set-series is easy
to handle; $V$ containing only non-negative integers guarantees the
holomorphy of $\varsigma_{V}\left(z\right)$ on $\mathbb{D}$, seeing
as negative integer elements of $V$ would produce negative integer
powers of $z$, which are singular at $0$. As such, we need for our
$H$s to have ``one-sided'' invariant sets; this feature extends
to the generalizations of set-series to discrete abelian groups othere
than $\mathbb{Z}$.

This condition is not as restrictive as it might seem; we can extend
it by \textbf{topological conjugation}, writing:
\[
H\left(n\right)\overset{\textrm{def}}{=}\textrm{sgn}\left(a\right)\left(T\left(\textrm{sgn}\left(a\right)n+a\right)-a\right),\textrm{ }\forall n\in\mathbb{N}_{0}
\]
where $T:\mathbb{Z}\rightarrow\mathbb{Z}$ is any map and $a$ is
any integer so that $H\left(n\right)$, as defined, becomes a $\varrho$-hydra
map for some $\varrho\in\mathbb{N}_{0}$. Here:
\[
\textrm{sgn}\left(a\right)\overset{\textrm{def}}{=}\begin{cases}
1 & \textrm{if }a\geq0\\
-1 & \textrm{if }a<0
\end{cases}
\]
is the \textbf{signum/sign function}.

\vphantom{}\textbullet{} \emph{Why surjectivity?}

As was seen in §1.2.1, surjectivity guarantees that the orbit classes
of $H$ on $\mathbb{N}_{0}$ form a partition of $\mathbb{N}_{0}$.
This partition property is essential, for it guarantees that we can
decompose $\varsigma_{\mathbb{N}_{0}}\left(z\right)=\frac{1}{1-z}$
as:
\[
\varsigma_{\mathbb{N}_{0}}\left(z\right)=\sum_{\begin{array}{c}
V\subseteq\mathbb{N}_{0}\\
V\textrm{ is an orbit class of }H
\end{array}}\varsigma_{V}\left(z\right)
\]
where the $\varsigma_{V}$s are pair-wise disjoint with respect to
$V$. This property is the very reason why the \textbf{Permutation
Operator Theorem }of the next section has anything useful to say at
all.

\vphantom{}(Lagarias. 2010) gives many examples of the kinds of questions
that have been asked Collatz-type maps. Most of these are qualitative:
how many distinct orbit classes does the map have (in particular,
are there infinitely many)?; what are the orbit classes' densities?;
are there any divergent orbit classes, and if so, how many?; are there
any open cycles, and if so, how many? Other paths of inquiry are more
quantiative: what upper or lower bounds can be achieved for the number
of elements in a cycle?; exploring connections between cycles and
problems in transcendence theory and exponential diophantine equations
(Tao. 2011), and the like. As discussed in the introduction, the motivation
for this paper was the desire to determine those families $\mathcal{H}$
of hydra maps for which the possession of a sufficient amount \emph{arithmetic
}regularity in a given orbit class $V$ of any $H\in\mathcal{H}$
(for example, $V$ being rational; $V$ containing a rational subset,
etc.) would then be sufficient to give monopoly on the dynamics of
$H$. Ideally, this would mean that the complement $\mathbb{N}_{0}\backslash V$
would be finite; a weaker (but still non-trivial) conclusion would
be that the complement $\mathbb{N}_{0}\backslash V$ would be a set
of zero density. To this end, we introduce two types of regularity
properties for a given family $\mathcal{H}$ of hydra maps.

\vphantom{}\textbf{Definition} \textbf{11}:\textbf{ }Let $\mathcal{H}$
be a collection of hydra maps on $\mathbb{N}_{0}$. For any $H\in\mathcal{H}$,
let $\textrm{W}\left(H\right)$ denote the set of orbit classes of
$H$ in $\mathbb{N}_{0}$.

I. A \textbf{Weak Rationality Theorem for $\mathcal{H}$ (WRT-$\mathcal{H}$)
}is the statement ``For any $H\in\mathcal{H}$, if there is a $V\in\textrm{W}\left(H\right)$
which \emph{is rational}, then $\mathbb{N}_{0}\backslash V$ is finite''.
$\mathcal{H}$ is said to \textbf{satisfy a Weak Rationality Theorem
}whenever the\textbf{ WRT-$\mathcal{H}$} is true.

II. A \textbf{Strong Rationality Theorem for $\mathcal{H}$ (SRT-$\mathcal{H}$)
}is the statement ``For any $H\in\mathcal{H}$, if there is a $V\in\textrm{W}\left(H\right)$
which \emph{contains a rational} \emph{subset}, then $\mathbb{N}_{0}\backslash V$
is finite''. $\mathcal{H}$ is said to \textbf{satisfy a Strong Rationality
Theorem }whenever the\textbf{ SRT-$\mathcal{H}$} is true.

\vphantom{}One particularly noteworthy consequence of a \textbf{SRT-$\mathcal{H}$
}is that, for any $H\in\mathcal{H}$, the existence of an IAP in some
$V\in\textrm{W}\left(H\right)$ then forces $\mathbb{N}_{0}\backslash V$
to be finite. The appellations of ``weak'' and ``strong'' refer
not to the minimal amount of structure required of $V$ in each case,
but to the ``strength'' (that is, generality) of the results obtained.
Rationality Theorems---particularly \textbf{SRT}s---can provide
a wedge by which one can approach certain Collatz-type conjectures.
As an example, Professor K.R. Matthews of the University of Queensland,
Australia, makes the following conjecture:

\vphantom{}\textbf{Matthews' Conjectures} (<http://www.numbertheory.org/php/markov.html>)):
consider the maps $T_{1}:\mathbb{Z}\rightarrow\mathbb{Z}$ and $T_{2}:\mathbb{Z}\rightarrow\mathbb{Z}$
defined by:
\[
T_{1}\left(n\right)\overset{\textrm{def}}{=}\begin{cases}
2n & \textrm{if }n\overset{3}{\equiv}0\\
\frac{7n+2}{3} & \textrm{if }n\overset{3}{\equiv}1\\
\frac{n-2}{3} & \textrm{if }n\overset{3}{\equiv}2
\end{cases}
\]
\[
T_{2}\left(n\right)\overset{\textrm{def}}{=}\begin{cases}
2n & \textrm{if }n\overset{3}{\equiv}0\\
\frac{5n-2}{3} & \textrm{if }n\overset{3}{\equiv}1\\
\frac{n-2}{3} & \textrm{if }n\overset{3}{\equiv}2
\end{cases}
\]
Then, he conjectures that:

I. Every trajectory {[}of $T_{1}${]} starting from $n\geq1$ will
eventually enter the zero residue class (mod $3$).

II. Every trajectory {[}of $T_{1}${]} starting from $n\leq-1$ will
eventually enter the zero residue class (mod $3)$, or reach one of
the cycles $-1,-1$ or $-2,-4,-2$.

III. Every trajectory {[}of $T_{2}${]} starting from $n\geq1$ will
eventually enter the zero residue class (mod $3$), or reach the cycle
$1,1$.

IV. Every trajectory {[}of $T_{2}${]} starting from $n\leq-1$ will
eventually enter the zero residue class (mod $3$), or reach one of
the cycles $-1,-1$ or $-2,-4,-2$.

\vphantom{}For the sake of the present discussion we will consider
only $T_{1}$. First, note that $\mathbb{N}_{0}$, and $-\mathbb{N}_{1}$
are orbit classes of $T_{1}$. As such, we shall split $T_{1}$ into
two $3$-hydra maps: $T_{+1},T_{-1}:\mathbb{N}_{0}\rightarrow\mathbb{N}_{0}$.
The former is simply the restriction of $T_{1}$ to $\mathbb{N}_{0}$:

\[
T_{+1}\left(n\right)\overset{\textrm{def}}{=}\begin{cases}
2n & \textrm{if }n\overset{3}{\equiv}0\\
\frac{7n+2}{3} & \textrm{if }n\overset{3}{\equiv}1\\
\frac{n-2}{3} & \textrm{if }n\overset{3}{\equiv}2
\end{cases}
\]
while the latter is the restriction of $T_{1}$ to $-\mathbb{N}_{0}$,
albeit conjugated by multiplication by $-1$:
\[
T_{-1}\left(n\right)=-T_{+1}\left(-n\right)
\]
This gives:
\[
T_{-1}\left(n\right)=\begin{cases}
2n & \textrm{if }-n\overset{3}{\equiv}0\\
-\frac{7\left(-n\right)+2}{3} & \textrm{if }-n\overset{3}{\equiv}1\\
-\frac{-n-2}{3} & \textrm{if }-n\overset{3}{\equiv}2
\end{cases}
\]
Since $-n\overset{3}{\equiv}0$ is the same as $n\overset{3}{\equiv}-0\overset{3}{\equiv}0$,
since $-n\overset{3}{\equiv}1$ is the same as $n\overset{3}{\equiv}-1\overset{3}{\equiv}2$,
and since $-n\overset{3}{\equiv}2$ is the same as $n\overset{3}{\equiv}-2\overset{3}{\equiv}1$,
this can be written as:
\[
T_{-1}\left(n\right)\overset{\textrm{def}}{=}\begin{cases}
2n & \textrm{if }n\overset{3}{\equiv}0\\
\frac{n+2}{3} & \textrm{if }n\overset{3}{\equiv}1\\
\frac{7n-2}{3} & \textrm{if }n\overset{3}{\equiv}2
\end{cases}
\]

Note that $\left\{ 0\right\} $ is then an orbit class of $T_{-1}$
(as opposed to the case for $T_{+1}$, where both $0$ and $2$ belong
to the same orbit class (since $T_{+1}\left(2\right)=\frac{2-2}{3}=0$)).

\vphantom{}\textbf{Proposition 4}: Let $V_{+1}$ (resp. $V_{-1}$)
denote the orbit class of $T_{+1}$ (resp. $T_{-1}$) in $\mathbb{N}_{0}$
containing $3\mathbb{N}_{1}$ (that is, containing all positive integer
multiples of $3$). Then $d\left(V_{+1}\right)=d\left(V_{-1}\right)=1$.

Proof: On page 24 of <https://arxiv.org/pdf/math/0608208.pdf>, Lagarias
comments that ``almost all trajectories contain an element $n\equiv0\textrm{ (mod 3)}$.''
Technically, Lagarias says this for a map which has the rule $7n+3$
for $n\overset{3}{\equiv}0$ instead of $2n$, but, since both $7n+3$
and $2n$ map $0$ mod $3$ to $0$ mod $3$, this ``almost all''
property works for both the version that Lagarias deals with and the
one given here. Moreover, the map Lagarias gives is defined on $\mathbb{Z}$,
rather than on $\mathbb{N}_{0}$. As such, assuming that when Lagarias
says ``almost all trajectories contain an element $n\equiv0\textrm{ (mod 3)}$'',
he \emph{means} ``the set of trajectories in $\mathbb{N}_{1}$ (resp.
$-\mathbb{N}_{1}$) that contain an element $n\equiv0\textrm{ (mod 3)}$
has a well-defined natural density, and that natural density is equal
to $1$'', this then implies that $d\left(V_{+}\right)=d\left(V_{-}\right)=1$,
as desired.

Q.E.D.

\vphantom{}\emph{Remark}:\emph{ }We will write $T_{\pm1}$ (likewise,
$V_{\pm1}$) when referring to both $T_{+1}$ and $T_{-1}$ (likewise,
$V_{+1}$ and $V_{-1}$) simultaneously.

\vphantom{}\textbf{Proposition} \textbf{5}: $T_{\pm1}:\mathbb{N}_{0}\rightarrow\mathbb{N}_{0}$
are surjective.

Proof:

\textbullet{} $T_{+1}$: Letting $n\in\mathbb{N}_{0}$ be of the form
$3m+k$ for $m\in\mathbb{N}_{0}$ and $k\in\left\{ 0,1,2\right\} $,
observe that:
\[
T_{+1}\left(3m+k\right)=\begin{cases}
2\left(3m+0\right) & \textrm{if }k\overset{3}{\equiv}0\\
\frac{7\left(3m+1\right)+2}{3} & \textrm{if }k\overset{3}{\equiv}1\\
\frac{\left(3m+2\right)-2}{3} & \textrm{if }k\overset{3}{\equiv}2
\end{cases}=\begin{cases}
6m & \textrm{if }k\overset{3}{\equiv}0\\
7m+3 & \textrm{if }k\overset{3}{\equiv}1\\
m & \textrm{if }k\overset{3}{\equiv}2
\end{cases}
\]
Thus $\left\{ T_{+1}\left(3m+2\right):m\in\mathbb{N}_{0}\right\} =\left\{ m:m\in\mathbb{N}_{0}\right\} =\mathbb{N}_{0}$,
and $T_{+}$ is therefore surjective. $\checkmark$

\textbullet{} $T_{-1}$: Letting $n\in\mathbb{N}_{0}$ be of the form
$3m+k$ for $m\in\mathbb{N}_{0}$ and $k\in\left\{ 0,1,2\right\} $,
observe that:
\[
T_{-1}\left(3m+k\right)=\begin{cases}
2\left(3m\right) & \textrm{if }k\overset{3}{\equiv}0\\
\frac{\left(3m+1\right)+2}{3} & \textrm{if }k\overset{3}{\equiv}1\\
\frac{7\left(3m+2\right)-2}{3} & \textrm{if }k\overset{3}{\equiv}2
\end{cases}=\begin{cases}
6m & \textrm{if }k\overset{3}{\equiv}0\\
m+1 & \textrm{if }k\overset{3}{\equiv}1\\
7m+4 & \textrm{if }k\overset{3}{\equiv}2
\end{cases}
\]
Thus:
\[
\left\{ T_{-1}\left(3m+1\right):m\in\mathbb{N}_{0}\right\} \cup\left\{ T\left(0\right)\right\} =\left\{ m+1:m\in\mathbb{N}_{0}\right\} \cup\left\{ 0\right\} =\mathbb{N}_{0}
\]
and $T_{-}$ is therefore surjective. $\checkmark$

Q.E.D.

\vphantom{}Thus, both of $T_{\pm1}$ are $3$-hydra maps on $\mathbb{N}_{0}$.
In particular, they are both prime and regulated. Note that if $T_{+1}$
satisfies a \textbf{SRT}, then Matthews' Conjecture (I) for $T_{+}$
will automatically hold true. Specifically, if it satisfies an \textbf{SRT},
then the fact that $V_{+1}$ contains the IAP $3\mathbb{N}_{1}$ would
force $\mathbb{N}_{0}\backslash V_{+1}$ to be finite, which then
forces $\mathbb{N}_{0}\backslash V_{+1}$ to either be empty, or to
be a union of finitely many closed cycles. As will be done as demonstrative
example in the next section, \textbf{Permutation Operators} can be
used to easily prove that $T_{+1}$ cannot have any closed cycles,
and thus, that $T_{+1}$ satisfying an \textbf{SRT }would be sufficient
to prove (I).

While tantalizing, establishing an \textbf{SRT }for $T_{+1}$ is likely
to be non-trivial. Indeed, while $T_{+1}$ might yet satisfy an \textbf{SRT},
the negative branch of $T_{1}$, $T_{-1}$, demonstrably \emph{fails
}to satisfy an \textbf{SRT}. Indeed, turning to $T_{-1}$:
\[
T_{-1}\left(n\right)=\begin{cases}
2n & \textrm{if }n\overset{3}{\equiv}0\\
\frac{n+2}{3} & \textrm{if }n\overset{3}{\equiv}1\\
\frac{7n-2}{3} & \textrm{if }n\overset{3}{\equiv}2
\end{cases}
\]
we note that, just as $V_{+1}$ contains $3\mathbb{N}_{1}$, so too
does $V_{-1}$. Thus, if $T_{-1}$ satisfied an \textbf{SRT}, this
would force $\mathbb{N}_{0}\backslash V_{-1}$ to be empty or finite.
However, this is clearly not the case. Observe that:

\[
T_{-1}\left(3m+k\right)=\begin{cases}
6m & \textrm{if }k\overset{3}{\equiv}0\\
m+1 & \textrm{if }k\overset{3}{\equiv}1\\
7m+4 & \textrm{if }k\overset{3}{\equiv}2
\end{cases},\textrm{ }\forall m\in\mathbb{N}_{0}
\]
As such:
\[
T_{-1}\left(3^{n}+1\right)=T_{-1}\left(3\cdot3^{n-1}+1\right)=3^{n-1}+1
\]
and so, every number in the set $\left\{ 3^{n}+1:n\in\mathbb{N}_{0}\right\} $
is iterated by $T_{-1}$ to the cycle $\left\{ 2,4\right\} $ (in
fact, $\left\{ 3^{n}+1:n\in\mathbb{N}_{0}\right\} $ is then an orbit
class of $T_{-1}$ in $\mathbb{N}_{0}$ of zero density, but with
infinitely many elements, which shows that $\left|\mathbb{N}_{0}\backslash V_{-1}\right|=\infty$,
and hence, $T_{-1}$ cannot satisfy an \textbf{SRT}. That being said,
one would think it significant that we can produce an explicit description
of this non-rational orbit class of $T_{-1}$; we can explicitly compute
the asymptotic expansion of $\left\{ 3^{n}+1:n\in\mathbb{N}_{0}\right\} $'s
counting function. As such, it leads the author to wonder if it might
be possible for an intermediate rationality theorem, wherein the hypothesis
of the \textbf{SRT }(an orbit class containing an IAP, or more generally,
a rational subset) at least guarantees that the complement of that
orbit class must possess some sort of predictable structure.

\pagebreak{}

\subsection{Permutation Operators on $\mathbb{Z}$}

As per our leitmotif, permutation operators encode the behavior of
dynamical systems as transformations of functions, enabling certain
dynamical questions to be tackled indirectly, by way of an analysis
of the operators' fixed points. The main takeaway (the \textbf{Permutation
Operator Theorem})\textbf{ }is that the fixed points are linear combinations
of set-series for orbit classes (and, more generally, invariant sets)
of the dynamical system in question.

\subsubsection{Definitions \& Formulae}

\textbf{Definition 12 (Permutation Operators} on $\mathbb{Z}$): Let
$H:\mathbb{N}_{0}\rightarrow\mathbb{N}_{0}$ be a map.

I. $\mathcal{Q}_{H}$, the \textbf{(Ordinary) Permutation Operator
}induced by $H$, is the linear operator on $\mathcal{Q}_{H}:\mathcal{A}\left(\mathbb{D}\right)\rightarrow\mathcal{A}\left(\mathbb{D}\right)$
defined by:
\[
\mathcal{Q}_{H}\left\{ \sum_{n=0}^{\infty}c_{n}z^{n}\right\} \overset{\textrm{def}}{=}\sum_{n=0}^{\infty}c_{H\left(n\right)}z^{n}
\]

II. $\mathscr{Q}_{H}$, the \textbf{(Fourier) Permutation Operator}
induced by $H$, is the linear operator $\mathscr{Q}_{H}:\mathcal{A}\left(\mathbb{H}_{+i}\right)\rightarrow\mathcal{A}\left(\mathbb{H}_{+i}\right)$
defined by:
\[
\mathscr{Q}_{H}\left\{ \sum_{n=0}^{\infty}c_{n}e^{2\pi inz}\right\} \overset{\textrm{def}}{=}\sum_{n=0}^{\infty}c_{H\left(n\right)}e^{2\pi inz}
\]

\vphantom{}In keeping with the §1.1's convention, we will prove statements
using permutation operators for ordinary set-series. Most of the analogous
results for Fourier permutation operators can be obtained \emph{mutatis
mutandis}; when the analogue is not as direct, however, we will provide
the formulae for the reader's benefit.

\emph{\vphantom{}}\textbf{Proposition} \textbf{6}: Let $H:\mathbb{N}_{0}\rightarrow\mathbb{N}_{0}$
be a map, and let $\left\{ c_{n}\right\} _{n\in\mathbb{N}_{0}}\subseteq\mathbb{C}$.
Then:
\[
\mathcal{Q}_{H}\left\{ \sum_{n=0}^{\infty}c_{n}z^{n}\right\} =\sum_{n=0}^{\infty}c_{n}\varsigma_{H^{-1}\left(\left\{ n\right\} \right)}\left(z\right)
\]

Additionally, for any $V\subseteq\mathbb{N}_{0}$:
\[
\mathcal{Q}_{H}\left\{ \varsigma_{V}\right\} \left(z\right)=\varsigma_{H^{-1}\left(V\right)}\left(z\right)
\]
where $H^{-1}\left(V\right)$ is, recall, the pre-image of $V$ under
$H$.

Proof: Consider the case of a set-series $\varsigma_{V}\left(z\right)$
for a set $V\subseteq\mathbb{N}_{0}$. Then:
\begin{eqnarray*}
\mathcal{Q}_{H}\left\{ \varsigma_{V}\right\} \left(z\right) & = & \mathcal{Q}_{H}\left\{ \sum_{v\in V}z^{v}\right\} \\
 & = & \mathcal{Q}_{H}\left\{ \sum_{n=0}^{\infty}\mathbf{1}_{V}\left(n\right)z^{n}\right\} \\
 & = & \sum_{n=0}^{\infty}\mathbf{1}_{V}\left(H\left(n\right)\right)z^{n}\\
\left(H\left(n\right)\in V\Leftrightarrow n\in H^{-1}\left(V\right)\right); & = & \sum_{n=0}^{\infty}\mathbf{1}_{H^{-1}\left(V\right)}\left(n\right)z^{n}\\
 & = & \varsigma_{H^{-1}\left(V\right)}\left(z\right)
\end{eqnarray*}
Now, letting $f\left(z\right)\overset{\textrm{def}}{=}\sum_{n=0}^{\infty}c_{n}z^{n}$
be arbitrary, let: 
\[
C_{f}\overset{\textrm{def}}{=}\left\{ c_{n}:n\in\mathbb{N}_{0}\right\} \backslash\left\{ 0\right\} 
\]
denote the set of non-zero values that occur as coefficients of $f$.
Then, for every complex number $c$, let: 
\[
I_{f}\left(c\right)\overset{\textrm{def}}{=}\left\{ n\in\mathbb{N}_{0}:c_{n}=c\right\} 
\]
Using these sets, observe that:
\[
f\left(z\right)=\sum_{n=0}^{\infty}c_{n}z^{n}=\sum_{c\in C_{f}}c\sum_{n\in I_{f}\left(c\right)}z^{n}=\sum_{c\in C_{f}}c\varsigma_{I_{f}\left(c\right)}\left(z\right)
\]
Consequently, by the linearity of $\mathcal{Q}_{H}$ and the set-series
case:

\begin{eqnarray*}
\mathcal{Q}_{H}\left\{ f\right\} \left(z\right) & = & \mathcal{Q}_{H}\left\{ \sum_{c\in C_{f}}c\varsigma_{I_{f}\left(c\right)}\right\} \left(z\right)\\
 & = & \sum_{c\in C_{f}}c\mathcal{Q}_{H}\left\{ \varsigma_{I_{f}\left(c\right)}\right\} \left(z\right)\\
 & = & \sum_{c\in C_{f}}c\mathcal{Q}_{H}\left\{ \varsigma_{I_{f}\left(c\right)}\right\} \left(z\right)\\
 & = & \sum_{c\in C_{f}}c\varsigma_{H^{-1}\left(I_{f}\left(c\right)\right)}\left(z\right)
\end{eqnarray*}
Next, for any disjoint sets $A,B\subseteq\mathbb{N}_{0}$, note that:
\[
H^{-1}\left(A\cup B\right)=H^{-1}\left(A\right)\cup H^{-1}\left(B\right)
\]
and that $H^{-1}\left(A\right)\cap H^{-1}\left(B\right)=\varnothing$.
Consequently, for every $c\in C_{f}$:
\[
\varsigma_{H^{-1}\left(I_{f}\left(c\right)\right)}=\varsigma_{H^{-1}\left(\left\{ n\in\mathbb{N}_{0}:c_{n}=c\right\} \right)}=\varsigma_{\bigcup_{n\in\mathbb{N}_{0}:c_{n}=c}H^{-1}\left(\left\{ n\right\} \right)}=\sum_{\begin{array}{c}
n\in\mathbb{N}_{0}\\
c_{n}=c
\end{array}}\varsigma_{H^{-1}\left(\left\{ n\right\} \right)}
\]
and so:
\begin{eqnarray*}
\mathcal{Q}_{H}\left\{ f\right\} \left(z\right) & = & \sum_{c\in C_{f}}c\varsigma_{H^{-1}\left(I_{f}\left(c\right)\right)}\left(z\right)\\
 & = & \sum_{c\in C_{f}}c\sum_{\begin{array}{c}
n\in\mathbb{N}_{0}\\
c_{n}=c
\end{array}}\varsigma_{H^{-1}\left(\left\{ n\right\} \right)}\left(z\right)\\
 & = & \sum_{n=0}^{\infty}c_{n}\varsigma_{H^{-1}\left(\left\{ n\right\} \right)}\left(z\right)
\end{eqnarray*}
which is the desired formula.

Q.E.D.

\vphantom{}\textbf{Proposition 7 (Permutation Operator Formulae}):
Let $H:\mathbb{N}_{0}\rightarrow\mathbb{N}_{0}$ be a $\varrho$-hydra
map. Then:
\begin{equation}
\mathcal{Q}_{H}\left\{ f\right\} \left(z\right)=\sum_{j=0}^{\varrho-1}\frac{1}{\mu_{j}}\sum_{k=0}^{\mu_{j}-1}\xi_{\mu_{j}}^{-kH\left(j\right)}f\left(\xi_{\mu_{j}}^{k}z^{\varrho/\mu_{j}}\right),\textrm{ }\forall f\in\mathcal{A}\left(\mathbb{D}\right)\label{eq:Q_h ordinary}
\end{equation}

\begin{equation}
\mathscr{Q}_{H}\left\{ \psi\right\} \left(z\right)=\sum_{j=0}^{\varrho-1}\frac{1}{\mu_{j}}\sum_{k=0}^{\mu_{j}-1}\xi_{\mu_{j}}^{-kH\left(j\right)}\psi\left(\frac{\varrho z+k}{\mu_{j}}\right),\textrm{ }\forall\psi\in\mathcal{A}\left(\mathbb{H}_{+i}\right)\label{eq:Q_H fourier}
\end{equation}
That is:
\begin{equation}
\mathcal{Q}_{H}\left\{ f\right\} \left(z\right)=\sum_{j=0}^{\varrho-1}z^{-b_{j}/a_{j}}\varpi_{\mu_{j},H\left(j\right)}\left\{ f\right\} \left(z^{\varrho/\mu_{j}}\right),\textrm{ }\forall f\in\mathcal{A}\left(\mathbb{D}\right)\label{eq:Q_h ordinary-1}
\end{equation}

\begin{equation}
\mathscr{Q}_{H}\left\{ \psi\right\} \left(z\right)=\sum_{j=0}^{\varrho-1}e^{-2\pi i\frac{b_{j}}{a_{j}}z}\tilde{\varpi}_{\mu_{j},H\left(j\right)}\left\{ \psi\right\} \left(\frac{\varrho z}{\mu_{j}}\right),\textrm{ }\forall\psi\in\mathcal{A}\left(\mathbb{H}_{+i}\right)\label{eq:Q_H fourier-1}
\end{equation}

Proof: Since we will actually be using \ref{eq:Q_H fourier}, it is
the one we shall prove. The proof is but a matter of compuation. To
begin, let:
\[
\psi\left(z\right)=\sum_{n=0}^{\infty}c_{n}e^{2\pi inz}\in\mathcal{A}\left(\mathbb{H}_{+i}\right)
\]
be arbitrary. Note that:
\[
H\left(\varrho n+j\right)=\frac{a_{j}\left(\varrho n+j\right)+b_{j}}{d_{j}}=\frac{\varrho a_{j}}{d_{j}}n+\frac{ja_{j}+b_{j}}{d_{j}}=\mu_{j}n+H\left(j\right)
\]
Now, by the definition of $\mathscr{Q}_{H}$:
\[
\mathscr{Q}_{H}\left\{ \psi\right\} \left(z\right)=\sum_{n=0}^{\infty}c_{H\left(n\right)}e^{2\pi inz}
\]
Splitting $n$ mod $\varrho$ gives:
\[
\mathscr{Q}_{H}\left\{ \psi\right\} \left(z\right)=\sum_{j=0}^{\varrho-1}\sum_{n=0}^{\infty}c_{H\left(\varrho n+j\right)}e^{2\pi i\left(\varrho n+j\right)z}=\sum_{j=0}^{\varrho-1}\sum_{n=0}^{\infty}c_{\mu_{j}n+H\left(j\right)}e^{2\pi i\left(\varrho n+j\right)z}
\]
To finish, note that we can invoke the $\tilde{\varpi}$ operators
if we re-write $\varrho n+j$ as $\mu_{j}n+H\left(j\right)$, like
so:
\begin{eqnarray*}
\left(\varrho n+j\right)z & = & \left(\varrho n+j\right)\left(\frac{a_{j}}{d_{j}}\times\frac{d_{j}}{a_{j}}z\right)\\
 & = & \left(\frac{\varrho a_{j}}{d_{j}}n+\frac{ja_{j}}{d_{j}}\right)\frac{d_{j}z}{a_{j}}\\
 & = & \left(\frac{\varrho a_{j}}{d_{j}}n+\frac{ja_{j}+b_{j}}{d_{j}}-\frac{b_{j}}{d_{j}}\right)\frac{d_{j}z}{a_{j}}\\
 & = & \left(\frac{\varrho a_{j}}{d_{j}}n+\frac{ja_{j}+b_{j}}{d_{j}}\right)\frac{d_{j}z}{a_{j}}-\left(\frac{b_{j}}{d_{j}}\times\frac{d_{j}z}{a_{j}}\right)\\
 & = & \left(\mu_{j}n+H\left(j\right)\right)\frac{d_{j}z}{a_{j}}-\frac{b_{j}z}{a_{j}}\\
 & = & \left(\mu_{j}n+H\left(j\right)\right)\frac{\varrho z}{\frac{\varrho a_{j}}{d_{j}}}-\frac{b_{j}z}{a_{j}}\\
 & = & \left(\mu_{j}n+H\left(j\right)\right)\frac{\varrho z}{\mu_{j}}-\frac{b_{j}z}{a_{j}}
\end{eqnarray*}
and so:
\begin{eqnarray*}
\mathscr{Q}_{H}\left\{ \psi\right\} \left(z\right) & = & \sum_{j=0}^{\varrho-1}\sum_{n=0}^{\infty}c_{\mu_{j}n+H\left(j\right)}e^{2\pi i\left(\varrho n+j\right)z}\\
 & = & \sum_{j=0}^{\varrho-1}\sum_{n=0}^{\infty}c_{\mu_{j}n+H\left(j\right)}e^{2\pi i\left(\left(\mu_{j}n+H\left(j\right)\right)\frac{\varrho z}{\mu_{j}}-\frac{b_{j}z}{a_{j}}\right)}\\
 & = & \sum_{j=0}^{\varrho-1}e^{-2\pi i\frac{b_{j}}{a_{j}}z}\sum_{n=0}^{\infty}c_{\mu_{j}n+H\left(j\right)}e^{2\pi i\left(\mu_{j}n+H\left(j\right)\right)\frac{\varrho z}{\mu_{j}}}
\end{eqnarray*}
Since 
\[
\tilde{\varpi}_{\mu_{j},H\left(j\right)}\left\{ \underbrace{\sum_{n=0}^{\infty}c_{n}e^{2\pi inz}}_{\psi}\right\} \left(z\right)=\sum_{n=0}^{\infty}c_{\mu_{j}n+H\left(j\right)}e^{2\pi i\left(\mu_{j}n+H\left(j\right)\right)z}
\]
it follows that:
\[
\sum_{n=0}^{\infty}c_{\mu_{j}n+H\left(j\right)}e^{2\pi i\left(\mu_{j}n+H\left(j\right)\right)\frac{\varrho z}{\mu_{j}}}=\tilde{\varpi}_{\mu_{j},H\left(j\right)}\left\{ \psi\right\} \left(\frac{\varrho z}{\mu_{j}}\right)
\]
and hence, that:
\[
\mathscr{Q}_{H}\left\{ \psi\right\} \left(z\right)=\sum_{j=0}^{\varrho-1}e^{-2\pi i\frac{b_{j}}{a_{j}}z}\tilde{\varpi}_{\mu_{j},H\left(j\right)}\left\{ \psi\right\} \left(\frac{\varrho z}{\mu_{j}}\right)
\]
Which is \ref{eq:Q_H fourier-1}. Using \ref{eq:fourier (a,b)-decomposition operator definition}
then gives \ref{eq:Q_H fourier}. Change of coordinates from $\mathbb{H}_{+i}$
to $\mathbb{D}$ then yields the other two formulae.

Q.E.D.

\subsubsection{The Permutation Operator Theorem}

In preparation for the shift to Fourier set-series, we state and prove
the \textbf{Permutation Operator Theorem }using Fourier set-series.

\textbf{Theorem 3: The Permutation Operator Theorem} (\textbf{POT}):
Let $H:\mathbb{N}_{0}\rightarrow\mathbb{N}_{0}$ be a surjective map.
Let $J$ be the (countable) index set so that $\left\{ V_{j}\right\} _{j\in J}$
is the collection of the distinct orbit classes of $H$ in $\mathbb{N}_{0}$.
Then, the Fourier set-series $\left\{ \psi_{V_{j}}\right\} _{j\in J}$
form a basis for $\textrm{Ker}\left(1-\mathscr{Q}_{H}\right)$; that
is to say, every fixed point of the permutation operator induced by
$H$ is a $\mathbb{C}$-linear combination of set-series of orbit
classes of $H$.

Proof: Let $H$ be as given. Let $J$ and the $V_{j}$s be as given.
By \textbf{Proposition 6}, it follows that: 
\begin{eqnarray*}
\mathscr{Q}_{H}\left\{ \psi_{V_{j}}\right\} \left(z\right) & = & \psi_{H^{-1}\left(V_{j}\right)}\left(z\right)\\
\left(V\textrm{ is an orbit class of }H\right); & = & \psi_{V_{j}}\left(z\right)
\end{eqnarray*}
and thus, that each $\psi_{V_{j}}$ is fixed by $\mathscr{Q}_{H}$.
Since the $V_{j}$s are pair-wise disjoint in $j$, the $\psi_{V_{j}}$s
are then linearly independent over $\mathbb{C}$.

To finish, let $\varphi\left(z\right)$ be an arbitrary fixed point
of $\mathscr{Q}_{H}$, and let $\left\{ c_{n}\right\} _{n\in\mathbb{N}_{0}}$
be the coefficients of $\varphi$'s Fourier series. Then, using the
notation from the proof of \textbf{Proposition 6}:
\[
\varphi\left(z\right)=\sum_{c\in C_{\varphi}}c\varphi_{I_{\varphi}\left(c\right)}\left(z\right)
\]
\[
\mathscr{Q}_{H}\left\{ \varphi\right\} \left(z\right)=\sum_{c\in C_{\varphi}}c\psi_{H^{-1}\left(I_{\varphi}\left(c\right)\right)}\left(z\right)
\]
where $\psi_{X}$ is a Fourier set-series for $X$. By the uniqueness
of series expansions, it must be that: 
\[
\psi_{I_{\varphi}\left(c\right)}=\psi_{H^{-1}\left(I_{\varphi}\left(c\right)\right)},\textrm{ }\forall c\in C_{\varphi}
\]
Since this is an equality of generic set-series, it occurs if and
only if $I_{\varphi}\left(c\right)=H^{-1}\left(I_{\varphi}\left(c\right)\right)$.
Thus, $I_{\varphi}\left(c\right)$ is an $H$-invariant set, and as
such, $I_{\varphi}\left(c\right)$ is necessarily a union of orbit
classes of $H$; for each $c\in C_{\varphi}$, there is an index set
$J_{c}\subseteq J$ so that:
\[
I_{\varphi}\left(c\right)=\bigcup_{j\in J_{c}}V_{j}
\]
Consequently:
\begin{eqnarray*}
\varphi & = & \sum_{c\in C_{\varphi}}c\psi_{I_{\varphi}\left(c\right)}\\
 & = & \sum_{c\in C_{\varphi}}c\psi_{\bigcup_{j\in J_{c}}V_{j}}\\
\left(\textrm{pair-wise disjointness}\right); & = & \sum_{c\in C_{\varphi}}\sum_{j\in J_{c}}c\psi_{V_{j}}
\end{eqnarray*}
Thus, the fixed point $\varphi$ is a $\mathbb{C}$-linear combination
of the $\psi_{V_{j}}$, which shows that the $\psi_{V_{j}}$s span
$\textrm{Ker}\left(1-\mathscr{Q}_{H}\right)$, as desired. $\checkmark$

Q.E.D.

\vphantom{}\emph{Example}: The Collatz Conjecture is equivalent to
the statement that $\textrm{Ker}\left(1-\mathscr{Q}_{3}\right)$ is
spanned by the functions $1$ and $\frac{e^{2\pi iz}}{1-e^{2\pi iz}}$;
equivalently, $\textrm{Ker}\left(1-\mathcal{Q}_{3}\right)$ by $1$
and $\frac{z}{1-z}$. Though the author discovered this result independently,
it was first obtained by (Berg \& Meinardus. 1995).

\vphantom{}Finally, as promised, here is a demonstrative example
of how to use permutation operators to prove that Matthews' Map $T_{+1}$
has no finite, non-empty invariant sets---and hence, no finite, non-empty
orbit classes.

\vphantom{}\textbf{Lemma 1}: Let $V$ be a $T_{+}$-invariant subset
of $\mathbb{N}_{0}$. If $\left|V\right|<\infty$, then $V=\varnothing$.

Proof: Let $V$ be a $T_{+}$-invariant subset of $\mathbb{N}_{0}$,
and suppose that $V$ is finite. By the \textbf{Permutation Operator
Theorem}, since $V$ is $T_{+}$-invariant, its ordinary set-series
$\varsigma_{V}\left(z\right)$ (which, in this case, is a digital
polynomial) is fixed by $\mathcal{Q}_{T_{+1}}$. Since $\varsigma_{\varnothing}\left(z\right)$
is the constant function $0$, to show that $V=\varnothing$, it suffices
to prove that the only polynomial fixed by is the zero polynomial;
that will force $\varsigma_{V}\left(z\right)$ to be identically $0$,
which forces $V$ to be empty.

Now, the formula for $\mathcal{Q}_{T_{+1}}$ is:
\begin{equation}
\mathcal{Q}_{T_{+1}}\left\{ f\right\} \left(z\right)=\frac{1}{6}\sum_{j=0}^{5}f\left(\xi_{6}^{j}z^{1/2}\right)+\frac{z^{-2/7}}{7}\sum_{k=0}^{6}\xi_{7}^{-3k}f\left(\xi_{7}^{k}z^{3/7}\right)+z^{2}f\left(z^{3}\right)\label{eq:Formula for Q_T_+1}
\end{equation}
So, suppose $f\left(z\right)$ is a polynomial which is fixed by $\mathcal{Q}_{T_{+1}}$.
This gives: 
\begin{equation}
f\left(z\right)=\frac{1}{6}\sum_{j=0}^{5}f\left(\xi_{6}^{j}z^{1/2}\right)+\frac{z^{-2/7}}{7}\sum_{k=0}^{6}\xi_{7}^{-3k}f\left(\xi_{7}^{k}z^{3/7}\right)+z^{2}f\left(z^{3}\right)\label{eq:Fixed point equation for Q_T_+1-1}
\end{equation}
Replacing $z$ with $z^{14}$ and moving things around gives:
\begin{equation}
f\left(z^{14}\right)-z^{28}f\left(z^{42}\right)=\frac{1}{6}\sum_{j=0}^{5}f\left(\xi_{6}^{j}z^{7}\right)+\frac{1}{7z^{4}}\sum_{k=0}^{6}\xi_{7}^{-3k}f\left(\xi_{7}^{k}z^{6}\right)\label{eq:T+1 functional equation (1)}
\end{equation}
Replacing $z$ with $\xi_{6}z$ and re-indexing the $j$-sum gives:
\begin{equation}
f\left(\xi_{3}z^{14}\right)-\xi_{3}^{2}z^{28}f\left(z^{42}\right)=\frac{1}{6}\sum_{j=0}^{5}f\left(\xi_{6}^{j}z^{7}\right)+\frac{\xi_{3}}{7z^{4}}\sum_{k=0}^{6}\xi_{7}^{-3k}f\left(\xi_{7}^{k}z^{6}\right)\label{eq:T+1 functional equation (2)}
\end{equation}
Subtracting \ref{eq:T+1 functional equation (2)} from \ref{eq:T+1 functional equation (1)}
cancels out the $j$-sum and leaves us with:
\begin{equation}
f\left(z^{14}\right)-f\left(\xi_{3}z^{14}\right)+\left(\xi_{3}^{2}-1\right)z^{28}f\left(z^{42}\right)=\frac{1-\xi_{3}}{7z^{4}}\sum_{k=0}^{6}\xi_{7}^{-3k}f\left(\xi_{7}^{k}z^{6}\right)\label{eq:T+1 functional equation (3)}
\end{equation}
Replacing $z$ with $z^{1/2}$ gives:
\begin{equation}
f\left(z^{7}\right)-f\left(\xi_{3}z^{7}\right)+\left(\xi_{3}^{2}-1\right)z^{14}f\left(z^{21}\right)=\frac{1-\xi_{3}}{7z^{2}}\sum_{k=0}^{6}\xi_{7}^{-3k}f\left(\xi_{7}^{k}z^{3}\right)\label{eq:T+1 functional equation (4)}
\end{equation}
Replacing $z$ with $\xi_{3}z$ gives:
\[
f\left(\xi_{3}z^{7}\right)-f\left(\xi_{3}^{2}z^{7}\right)+\xi_{3}^{2}\left(\xi_{3}^{2}-1\right)z^{14}f\left(z^{21}\right)=\frac{\xi_{3}-\xi_{3}^{2}}{7z^{2}}\sum_{k=0}^{6}\xi_{7}^{-3k}f\left(\xi_{7}^{k}z^{3}\right)
\]
Multiplying by $\xi_{3}^{2}$ on both sides gives:
\begin{equation}
\xi_{3}^{2}\left(f\left(\xi_{3}z^{7}\right)-f\left(\xi_{3}^{2}z^{7}\right)\right)+\left(1-\xi_{3}\right)z^{14}f\left(z^{21}\right)=\frac{1-\xi_{3}}{7z^{2}}\sum_{k=0}^{6}\xi_{7}^{-3k}f\left(\xi_{7}^{k}z^{3}\right)\label{eq:T+1 functional equation (5)}
\end{equation}
Subtracting \ref{eq:T+1 functional equation (5)} from \ref{eq:T+1 functional equation (4)}
cancels out the $k$-sum, leaving:
\[
-\left(1+\xi_{3}^{2}\right)f\left(\xi_{3}z^{7}\right)+\xi_{3}^{2}f\left(\xi_{3}^{2}z^{7}\right)+f\left(z^{7}\right)+\left(\xi_{3}^{2}+\xi_{3}-2\right)z^{14}f\left(z^{21}\right)=0
\]
Since $\xi_{3}^{2}+\xi_{3}=-1$ and $1+\xi_{3}^{2}=-\xi_{3}$, this
can be written as:
\begin{eqnarray*}
\xi_{3}f\left(\xi_{3}z^{7}\right)+\xi_{3}^{2}f\left(\xi_{3}^{2}z^{7}\right)+f\left(z^{7}\right)-3z^{14}f\left(z^{21}\right) & = & 0\\
\left(+3z^{14}f\left(z^{21}\right),\div3\right); & \Updownarrow\\
\frac{1}{3}\left(f\left(z^{7}\right)+\xi_{3}f\left(\xi_{3}z^{7}\right)+\xi_{3}^{2}f\left(\xi_{3}^{2}z^{7}\right)\right) & = & z^{14}f\left(z^{21}\right)\\
\left(z\rightarrow z^{1/7}\right); & \Updownarrow\\
\frac{1}{3}\sum_{k=0}^{3-1}\xi_{3}^{k}f\left(\xi_{3}^{k}z\right) & = & z^{2}f\left(z^{3}\right)\\
\left(\xi_{3}^{k}=\xi_{3}^{-2k}\right); & \Updownarrow\\
\underbrace{\frac{1}{3}\sum_{k=0}^{3-1}\xi_{3}^{-2k}f\left(\xi_{3}^{k}z\right)}_{\varpi_{3,2}\left\{ f\right\} \left(z\right)} & = & z^{2}f\left(z^{3}\right)
\end{eqnarray*}
And so:
\begin{equation}
\varpi_{3,2}\left\{ f\right\} \left(z\right)=z^{2}f\left(z^{3}\right)\label{eq:Final Identity for T+1}
\end{equation}
Since $f\left(z\right)$ is a polynomial, so are $\varpi_{3,2}\left\{ f\right\} \left(z\right)$
and $z^{2}f\left(z^{3}\right)$. Moreover:
\[
\deg f\geq\deg\left(\varpi_{3,2}\left\{ f\right\} \right)
\]
\[
\deg\left(z^{2}f\left(z^{3}\right)\right)=2+3\deg f
\]
and so: 
\begin{eqnarray*}
\varpi_{3,2}\left\{ f\right\} \left(z\right) & = & z^{2}f\left(z^{3}\right)\\
\left(\deg\right); & \Downarrow\\
\deg f\geq\deg\left(\varpi_{3,2}\left\{ f\right\} \right) & = & 2+3\deg f\\
 & \Downarrow\\
-2 & \geq & 2\deg f\\
 & \Updownarrow\\
-1 & \geq & \deg f
\end{eqnarray*}
Thus, $f$ must have negative degree. This forces $f$ to be identically
zero, since the zero polynomial is the only polynomial with a negative
degree (namely, $\deg\left(0\right)\overset{\textrm{def}}{=}-\infty$).
Since $f$ was an arbitrary polynomial fixed by $\mathcal{Q}_{T_{+1}}$,
we conclude that the only polynomial fixed point of $\mathcal{Q}_{T_{+1}}$
is the zero polynomial, and hence, that $V$ must be empty, as desired.

Q.E.D.

\pagebreak{}

\section{Dreamcatchers and Singularity Conservation Laws for Hydra Maps}

As was seen in Chapter 2, in passing from a set $V$ to its set-series,
the arithmetic and asymptotic properties of $V$ are encoded in the
set-series' boundary behavior. Given appropriate assumptions on $V$,
the ``dominant'' components of this boundary behavior can be captured
and thereby fashioned into well-defined mathematical objects that
the author calls \textbf{Dreamcatchers}. By a minor miracle, these
can be shown to be elements of the space $L^{2}\left(\mathbb{Q}/\mathbb{Z}\right)$.
Crucially, the dreamcatchers generated by a set-series fixed by permutation
operators will ``inherit'' a form of the functional equation that
characterized their progenitor. This the inherited equation---the\textbf{
Singularity Conservation Law} (\textbf{SCL}) in the paper's title---is
a statement of the fact that the permutation operators cannot change
a given type of singularity in its input function to a different type
of singularity in the output function. As such, we can approach the
task of finding fixed points of permutation operators by examining
what those operators do to singularities of a given growth rate.

We shalll begin this chapter by defining the notions of \textbf{virtual
poles }and \textbf{virtual residues}. These can then be used to define
the \textbf{dreamcatcher }of a set-series (or, more generally, of
any analytic function given by some infinite series which is only
known to converge on some region in $\mathbb{C}$). In the case where
the set-series represent invariant sets of some $\varrho$-hydra map,
it will be shown that the associated dreamcatchers satisfy a \textbf{SCL}.

\subsection{Origins}

Fix a $\varrho$-hydra map $H:\mathbb{N}_{0}\rightarrow\mathbb{N}_{0}$.
Considering an arbitrary $\psi\in\mathcal{A}\left(\mathbb{H}_{+i}\right)$,
the formula \ref{eq:Q_H fourier} for the effect of $\mathscr{Q}_{H}$
on $\psi$ is given by:
\[
\mathscr{Q}_{H}\left\{ \psi\right\} \left(z\right)=\sum_{j=0}^{\varrho-1}\frac{1}{\mu_{j}}\sum_{k=0}^{\mu_{j}-1}\xi_{\mu_{j}}^{-kH\left(j\right)}\psi\left(\frac{\varrho z+k}{\mu_{j}}\right)
\]
Rather than choosing to die upon the hill of Fourier series with which
we originally defined $\mathscr{Q}_{H}$, let us consider what happens
when we apply the permutation operator to functions defined by forms
other than Fourier series. (Besides, that hill is waiting for us in
Chapter 3). Testing out the effects of $\mathscr{Q}_{H}$ on various
every-day functions ($e^{\alpha z},\frac{1}{\left(z-c\right)^{n}},\ln z,\sqrt{z-c},c_{0}+c_{1}z+c_{2}z^{2}$,
etc.) for a given $H$---reveals some noteworthy, if simple, phenomena.

First and most importantly, since the right-hand side of \ref{eq:Q_H fourier}
is a sum of products of compositions of entire functions (exponentials,
the Möbius transformations $\frac{\varrho z+k}{\mu_{j}}$, etc.),
note that $\mathscr{Q}_{H}\left\{ \psi\right\} \left(z\right)$ will
be an entire function whenever $\psi\left(z\right)$ is entire. An
immediate corollary of this is that if $\mathscr{Q}_{H}\left\{ \psi\right\} \left(z\right)$
has any singularities, then those singularities are necessarily the
progeny of singularities of $\psi\left(z\right)$. Given a function
of the form:
\[
\psi\left(z\right)=\sum_{m=0}^{\infty}\sum_{n=1}^{d_{m}}\frac{c_{m,n}}{\left(z-s_{m}\right)^{n}}
\]
where $S$ is some discrete subset of $\mathbb{C}$, the $d_{m}$s
are positive integers, and the $c_{m,n}$s are arbitrary complex constants
(where, for each $m$, there is at least one value of $n$ for which
$c_{m,n}\neq0$) it can be shown that any singularities of $\mathscr{Q}_{H}\left\{ \psi\right\} \left(z\right)$
will be of the form $z=\frac{\varrho s_{m}+k}{\mu_{j}}$ for some
$m,j,k$, and that it is possible for the images of distinct singularities
$s_{m}$ and $s_{m^{\prime}}$ under $\mathscr{Q}_{H}$ to end up
cancelling each other out if the $c_{m,n}$s and $c_{m^{\prime},n}$s
satisfy certain systems of linear equations involving the roots of
unity present in $\mathscr{Q}_{H}$.

If we consider $\psi\left(z\right)$ as representing the ``meromorphic
part'', as it were, of an $H$-invariant Fourier set-series $\psi_{V}\left(z\right)$,
tools of analytic number theory such as the \textbf{HLTT }tell us
that the $s_{m}$s, the $n$s, and the $c_{m,n}$s are intimately
connected with the arithmetic and asymptotic properties of the set
$V$. Moreover, being an $H$-invariant Fourier set-series, $\psi_{V}\left(z\right)$
is necessarily fixed by $\mathscr{Q}_{H}$, and as such, those singularities
of $\psi_{V}\left(z\right)$ which are represented in $\psi\left(z\right)$
must exist in a delicate balance with the workings of $\mathscr{Q}_{H}$:
any singularities $s_{1},s_{2}$ of $\psi\left(z\right)$ which end
up cancelling one another in $\mathscr{Q}_{H}\left\{ \psi\right\} \left(z\right)$
must be replaced by the singularities generated by some singularities
$s_{3}$ and $s_{4}$ of $\psi\left(z\right)$. That is, to say these
singularities obey a \textbf{Singularity Conservation Law}.

\subsection{Edgepoints, Virtual Poles, and Virtual Residues}

The study of the singular behavior of functions defined by infinite
series (power series, dirichlet series, exponential power series (a.k.a.
``generalized dirichlet series''), Fourier series) is one of the
perennial topics of analysis. The classical complex-analytic treatment
of isolated singularities, however, is insufficient for the task at
hand. As\textbf{ }a motivational example, let us define:
\[
\phi_{d}\left(z\right)\overset{\textrm{def}}{=}\sum_{n=0}^{\infty}z^{d^{n}},\textrm{ }\forall d\in\mathbb{N}_{2}
\]
and then consider the set-series:

\[
\varsigma_{V}\left(z\right)\overset{\textrm{def}}{=}\phi_{2}\left(z\right)+\frac{1}{1-z^{3}}
\]
where, here: 
\[
V=2^{\mathbb{N}_{0}}\cup\left(3\mathbb{N}_{0}\right)
\]
and, as such, $\varsigma_{V}\left(z\right)$ will then inherit $\phi_{2}\left(z\right)$'s
(strong) natural boundary on $\partial\mathbb{D}$. As much as the
presence of the rational function $\frac{1}{1-z^{3}}$ might tempt
us to say that $\varsigma_{V}\left(z\right)$ has a pole at $1$,
$\xi_{3}$, and $\xi_{3}^{2}$, this would be a falsehood; the classical
theory breaks down in the face of the fact that these three singularities
are non-isolated, with the singularities that $\varsigma_{V}\left(z\right)$
inherits from $\phi_{2}\left(z\right)$ (namely $\bigcup_{n=1}^{\infty}\left\{ \xi_{2^{n}}^{m}:0\leq m\leq2^{n}-1\right\} $)
constitutes a dense subset of the unit circle.

Fortunately, there is a way to overcome this pathological situation.
An application of the \textbf{HLTT }shows that:
\[
\frac{1}{1-x^{3}}\sim\frac{1/3}{1-x}\textrm{ as }x\uparrow1
\]
Moreover, a Mellin transform analysis of $\phi_{2}\left(z\right)$
yields the (exact) asymptotic:
\[
\phi_{2}\left(x\right)\sim\frac{1}{2}-\frac{\gamma+\ln\left(-\ln x\right)}{\ln2}+\underbrace{\frac{1}{\ln2}\sum_{k\in\mathbb{Z}\backslash\left\{ 0\right\} }\Gamma\left(\frac{2k\pi i}{\ln2}\right)e^{-2k\pi i\log_{2}\left(-\ln x\right)}}_{\textrm{bounded, but oscillates wildly}}\textrm{ as }x\uparrow1
\]
giving $\varsigma_{V}\left(x\right)$ a singular expansion of:
\[
\varsigma_{V}\left(x\right)\asymp\frac{1/3}{1-x}+\frac{1}{2}-\frac{\gamma+\ln\left(-\ln x\right)}{\ln2}+\frac{1}{\ln2}\sum_{k\in\mathbb{Z}\backslash\left\{ 0\right\} }\Gamma\left(\frac{2k\pi i}{\ln2}\right)e^{-2k\pi i\log_{2}\left(-\ln x\right)}\textrm{ as }x\uparrow1
\]
Since the asymptotics of $\frac{1}{x}$ as $x\downarrow0$ are identical
to those of $\frac{1}{1-x}$ as $x\uparrow1$, we can write this as:
\[
\varsigma_{V}\left(e^{-x}\right)\asymp\frac{1/3}{x}+\frac{1}{2}-\frac{\gamma+\ln x}{\ln2}+\frac{1}{\ln2}\sum_{k\in\mathbb{Z}\backslash\left\{ 0\right\} }\Gamma\left(\frac{2k\pi i}{\ln2}\right)e^{-2k\pi i\log_{2}x}\textrm{ as }x\downarrow0
\]
with the asymptotic expansion (the part that actually grows as $x\downarrow0$)
being given by:
\[
\varsigma_{V}\left(e^{-x}\right)=\frac{1}{3x}-\log_{2}x+O\left(1\right)\textrm{ as }x\downarrow0
\]
As such, the simple pole at $x=0$ is indeed present in $\varsigma_{V}\left(e^{-x}\right)$,
but it is lost among the weeds generated by the fractal oscillations
of the infinite $k$-sum and, of course, the sluggish growth of the
logarithmic term . Nevertheless, if we agree to disregard any growth
rates in the asymptotic expansion of $\varsigma_{V}\left(e^{-x}\right)$
of order \emph{less than} $O\left(\frac{1}{x}\right)$, we can say
that the dominant singular behavior of the function at $x=0$ is $\frac{1}{3x}$.
It is this precisely this idea of disregarding the lower-order singular
terms that leads to the concept of a \emph{virtual pole}.

\vphantom{}\textbf{Definition 13}

\textbf{(Edgepoints, Virtual Residues \& Virtual Poles, and Related
Terminology)}: Consider a function $\psi\in\mathcal{A}\left(\mathbb{H}_{+i}\right)$.

I. $x\in\mathbb{R}$ is said to be an \textbf{edgepoint }of $\psi$
whenever $\lim_{y\downarrow0}y\psi\left(x+iy\right)$ exists and is
finite. The set of all edgepoints of $\psi$ is denoted $\textrm{P}\left(\psi\right)$:
\[
\textrm{P}\left(\psi\right)\overset{\textrm{def}}{=}\left\{ x\in\mathbb{R}:\lim_{y\downarrow0}y\psi\left(x+iy\right)\textrm{ exists and is finite}\right\} 
\]
If $\psi$ is a set-series for a set $V\subseteq\mathbb{N}_{0}$,
the author will abuse terminology slightly and speak of the \textbf{edgepoints
of $V$} to refer to the edgepoints of the set-series of $V$ under
consideration.

II. For each $x\in\textrm{P}\left(\psi\right)$, the \textbf{virtual
residue }of $\psi$ at $x$ is denoted $\textrm{VRes}_{x}\left[\psi\right]$,
and is defined by:
\[
\textrm{VRes}_{x}\left[\psi\right]\overset{\textrm{def}}{=}\lim_{y\downarrow0}y\psi\left(x+iy\right),\textrm{ }\forall x\in\textrm{P}\left(\psi\right)
\]

III. An edgepoint $x$ of $\psi$ is said to be a \textbf{virtual
pole }of $\psi$ whenever $\textrm{VRes}_{x}\left[\psi\right]\neq0$.

IV. A virtual pole is said to be \textbf{classical }if $x$ is a pole
of $\psi\left(z\right)$ in the classical sense---i.e., if $\psi\left(z\right)-\frac{\textrm{VRes}_{x}\left[\psi\right]}{z-x}$
is holomorphic $z=x$.

\vphantom{}It is not an accident that the definitions of virtual
poles and virtual residues only take into account the case where $x$
is of ``degree $1$''; when working with set-series, these are the
only kinds of virtual poles that can ever appear.

\vphantom{}\textbf{Proposition 8}: Let $V\subseteq\mathbb{N}_{0}$
be arbitrary. Then, $V$'s Fourier set-series satisfies:
\[
\lim_{y\downarrow0}y^{1+\delta}\psi\left(x+iy\right)=0,\textrm{ }\forall\delta\in\mathbb{R}>0,\textrm{ }\forall x\in\mathbb{R}
\]

Proof: Let $x\in\mathbb{R}$ and let $\delta\in\mathbb{R}>0$ be arbitrary.
Then:
\begin{eqnarray*}
\left|y^{1+\delta}\psi_{V}\left(x+iy\right)\right| & = & \left|y^{1+\delta}\sum_{n=0}^{\infty}\mathbf{1}_{V}\left(n\right)e^{2n\pi i\left(x+iy\right)}\right|\\
 & \leq & y^{1+\delta}\sum_{n=0}^{\infty}\mathbf{1}_{V}\left(n\right)e^{-2n\pi y}\\
\left(V\subseteq\mathbb{N}_{0}\right); & \leq & y^{1+\delta}\sum_{n=0}^{\infty}e^{-2n\pi y}\\
\left(\forall x+iy\in\mathbb{H}_{+i}\right); & = & \frac{y^{1+\delta}}{1-e^{-2\pi y}}
\end{eqnarray*}
Since:
\begin{eqnarray*}
\lim_{y\downarrow0}\frac{y^{1+\delta}}{1-e^{-2\pi y}} & = & \lim_{y\downarrow0}\frac{y^{1+\delta}}{y}\frac{y}{1-e^{-2\pi y}}\\
 & = & \lim_{y\downarrow0}y^{\delta}\frac{1}{2\pi e^{-2\pi y}}\\
 & = & \frac{1}{2\pi}\lim_{y\downarrow0}y^{\delta}\\
\left(\delta>0\right); & = & 0
\end{eqnarray*}
it follows that:
\[
\lim_{y\downarrow0}\left|y^{1+\delta}\psi_{V}\left(x+iy\right)\right|\leq\lim_{y\downarrow0}\frac{y^{1+\delta}}{1-e^{-2\pi y}}=\frac{1}{2\pi}\lim_{y\downarrow0}y^{\delta}=0
\]
and so:
\[
\lim_{y\downarrow0}y^{1+\delta}\psi_{V}\left(x+iy\right)=0,\textrm{ }\forall x\in\mathbb{R}
\]
as desired.

Q.E.D.

\vphantom{}It is because of \textbf{Proposition 8 }that we shall
speak of virtual poles as encoding the ``dominant'' asymptotic behavior
of set-series at the boundaries of their regions of convergence. Moreover,
by definition of $\textrm{P}\left(\psi_{V}\right)$, note that this
shows that: 
\[
\lim_{y\downarrow0}y\psi_{V}\left(x+iy\right)\textrm{ does not exist}\Leftrightarrow x\in\mathbb{R}\backslash\textrm{P}\left(\psi_{V}\right)
\]
seeing as $x\in\mathbb{R}\backslash\textrm{P}\left(\psi_{V}\right)$
forces the limit on the left to be either infinite or non-existent,
whereas \textbf{Proposition 8} shows that $\psi_{V}\left(x+iy\right)=O\left(\frac{1}{y}\right)$
for all $x\in\mathbb{R}$.

As something of a brief interlude, before we proceed on to construct
dreamcatchers, the function $\phi_{d}\left(z\right)$ mentioned at
the start of the section is noteworthy enough for it to be worth our
while to discuss it a bit before proceeding to the example. $\phi_{d}\left(z\right)$
is probably the best known example of a power series with a (strong)
natural boundary on the unit circle. Not only is it of interest to
pure analysis as a prototypical example of a \textbf{lacunary series
}(Mandelbrojt. 1927), it is also the genesis for an entire branch
of transcendental number theory---Mahler Theory, named after Kurt
Mahler---which utilizes the functional equations that characterize
such functions (for $\phi_{d}$, the equation is $\phi_{d}\left(z^{d}\right)=\phi_{d}\left(z\right)-z$)
to establish the transcendence of those functions at various algebraic
values. For instance, for all $d\in\mathbb{N}_{2}$, Mahler himself
established the transcendence of $\phi_{d}\left(\alpha\right)$, where
$\alpha\in\mathbb{C}$ is any algebraic number satisfying $0<\left|\alpha\right|<1$.
(Nishioka. 1996).

Mahler Theory has deep connections with automata theory, via the concept
of an ``automatic sequence''. In terms of power series, it follows
from a result of Eilenberg's (cited in (B. Adamczewski and J. Bell.
2017)) that the ``automatic'' property of a sequence is equivalent
to said sequence being the sequence of coefficients of a power series
(an ``automatica series'') with the property that the linear space
generated by applying iterates of certain linear operators (called
\textbf{Cartier Operators}) is finite dimensional. It is a simple
computation to show that the Cartier Operators and the action-decomposition
operators $\varpi_{a,b}$ considered in Chapter 1 can be expressed
in terms of one another. Though the similarities between these areas
of transcendence theory and this paper's work are striking, there
appear to be fundamental technical limitations that keep the Mahler
theory methods---fascinating though they are---from being applied
to analyze the functional equations for fixed points of the permutation
operators induced by the Collatz map and its kin. These can be seen
merely by glacing at, say, the fixed point equation for $\mathcal{Q}_{H_{3}}$:
\[
f\left(z\right)=f\left(z^{2}\right)+\frac{z^{-\frac{1}{3}}}{3}\sum_{k=0}^{2}\xi_{3}^{2k}f\left(\xi_{3}^{k}z^{2/3}\right)
\]
where we note the presence of both $3$ and $2$ in the powers of
$z$ within $f$, as well as the precomposition with multiplication
by $\xi_{3}^{k}$. That being said, it is the author's hope that the
methods of dreamcatcher analysis presented in this paper will be of
use in Mahler theory and related fields.

\pagebreak{}

\subsection{Constructing Dreamcatchers}

A dreamcatcher is, in part, an object that keeps track of the two
pieces of data associated to virtual poles: their \emph{location }and
their \emph{residues}. However, it is not enough simply to list this
data in the form of a set of location-residue pairs. We need to give
it some sort of regularity, so that our record of this data will be
in a form amenable to analytical techniques.

\vphantom{}\textbf{Definition} \textbf{14}: When speaking of a dreamcatcher,
we will use the term ``\textbf{support}'' (as in, ``the support
of the dreamcatcher''; ``the dreamcatcher is supported on {[}insert
set here{]}'') to denote the set of edgepoints of the function that
generated the dreamcatcher in question.

\vphantom{}The simplest case is where a dreamcatcher's support is
a finite---or, more generally, a \emph{discrete}---subset of $\mathbb{C}$.
In that case, dreamcatchers are essentially a spruced-up version of
what, in the asymptotic analysis of meromorphic functions, is known
as a ``singular expansion''. While meromorphic dreamcatchers have
their uses, for the present analysis, the $L^{2}$ case predominates.
Nevertheless, it is instructive to at least delinate the essentials
of the meromorphic case, as they provide a conceptual bedrock for
the generalizations to come.

In the types of asymptotic analysis common to complex analysis and
analytic number theory, one often works with ``singular expansions''
of meromorphic functions (ex: the ``principal part'' of a function's
Laurent series representation). For example, at $s=1$, the Riemann
Zeta function has the Laurent series:
\[
\zeta\left(s\right)=\frac{1}{s-1}+\sum_{n=0}^{\infty}\gamma_{n}\frac{\left(-1\right)^{n}}{n!}\left(s-1\right)^{n}
\]
where the real numbers $\left\{ \gamma_{n}\right\} _{n\in\mathbb{N}_{0}}$
are known as the \textbf{Stieltjes Constants}, with $\gamma_{0}=\gamma$,
the \textbf{Euler-Mascheroni Constant}. Since $\zeta\left(s\right)$
is a meromorphic function with exactly one singularity---the simple
pole at $s=1$---the function $\zeta\left(s\right)-\frac{1}{s-1}$
defiend by the Stieltjes constants' power series is entire. Consequently,
everything we need to know about $\zeta\left(s\right)$'s singular
behavior is contained in the $\frac{1}{s-1}$. This, the principal
part of $\zeta\left(s\right)$'s Laurent series at $s=1$, is the
meromorphic dreamcatcher of $\zeta\left(s\right)$. Our next proposition
generalizes this procedure to Fourier set-series.

\vphantom{}\textbf{Proposition 9 (Meromorphic Dreamcatchers}):\textbf{
}Let $V\subseteq\mathbb{N}_{0}$. $\psi_{V}\left(z\right)$ is said
to \textbf{admit a meromorphic dreamcatcher }whenever the set $T_{V}\subseteq\mathbb{R}$
of $\psi_{V}$'s virtual poles is a discrete\emph{ }subset of $\mathbb{R}$.
In that case, the \textbf{meromorphic dreamcatcher }of $\psi_{V}$
is the function: 
\[
F_{V}\left(z\right)\overset{\textrm{def}}{=}\sum_{t\in T_{V}}\frac{\textrm{VRes}_{t}\left[\psi_{V}\right]}{z-t}
\]
which is $1$-periodic and convergent for all $z\in\mathbb{C}\backslash T_{V}$.

Proof: Let everything be as given. Since $\psi_{V}\left(z\right)$
is $1$-periodic, so is the function $t\mapsto\textrm{VRes}_{t}\left[\psi_{V}\right]$
for $t\in\textrm{P}\left(\psi_{V}\right)$. Since $T_{V}$ is $1$-periodic,
its discreteness then forces it to be of the form:
\[
T_{V}=\bigcup_{m=1}^{M}\left(t_{m}+\mathbb{Z}\right)
\]
for some $M\in\mathbb{N}_{1}$ and some real numbers $t_{1},\ldots,t_{M}$
with $0\leq t_{1}<\cdots<t_{M}<1$. As such:
\begin{eqnarray*}
F_{V}\left(z\right) & = & \sum_{t\in T_{V}}\frac{\textrm{VRes}_{t}\left[\psi_{V}\right]}{z-t}\\
 & = & \sum_{m=1}^{M}\sum_{n\in\mathbb{Z}}\frac{\textrm{VRes}_{t_{m}+n}\left[\psi_{V}\right]}{z-\left(t_{m}+n\right)}\\
\left(\textrm{VRes}_{t_{m}+n}\left[\psi_{V}\right]=\textrm{VRes}_{t_{m}}\left[\psi_{V}\right]\right); & = & \sum_{m=1}^{M}\textrm{VRes}_{t_{m}}\left[\psi_{V}\right]\sum_{n\in\mathbb{Z}}\frac{1}{\left(z-t_{m}\right)-n}\\
\left(\sum_{n\in\mathbb{Z}}\frac{1}{z-n}=\pi\cot\left(\pi z\right),\textrm{ }\forall z\in\mathbb{C}\backslash\mathbb{Z}\right); & = & \sum_{m=1}^{M}\pi\cot\left(\pi\left(z-t_{m}\right)\right)\textrm{VRes}_{t_{m}}\left[\psi_{V}\right]
\end{eqnarray*}
and thus, being a linear combination of finitely many translates of
$\cot\left(\pi z\right)$, $F_{V}\left(z\right)$ is then meromorphic
on $\mathbb{C}$, with its series being compactly convergent for $z\in\mathbb{C}\backslash T_{V}$.

Q.E.D.

\vphantom{}Continuing in this vein, essentially everything we will
deal with for the remainder of this chapter can be done in the case
of meromorphic dreamcatchers. The common thread between the two approaches
(aside from taking virtual residues every five minutes) is that we
restrict our attention to the \emph{singular }part of the set-series
under consideration so as to obtain an object (the dreamcatcher) which
is simpler---and thus \emph{better behaved}---than the set-series
from whence it came. For example, the Fourier analogue $\psi_{2^{\mathbb{N}_{0}}\cup\left(3\mathbb{N}_{0}\right)}\left(z\right)$
of the example considered in §2.1.1 has a natural boundary on the
real axis, and, since the set of powers of $2$ has density $0$,
it can be shown that $\psi_{2^{\mathbb{N}_{0}}\cup\left(3\mathbb{N}_{0}\right)}\left(z\right)$
has the same set of virtual poles---the set $\frac{1}{3}\mathbb{Z}$---as
the function: 
\[
\psi_{3\mathbb{N}_{0}}\left(z\right)=\frac{1}{1-\left(e^{2\pi iz}\right)^{3}}
\]
which is rational in $e^{2\pi iz}$. But the dreamcatcher for $\psi_{2^{\mathbb{N}_{0}}\cup\left(3\mathbb{N}_{0}\right)}\left(z\right)$
cares not for the natural boundary induced by $2^{\mathbb{N}_{0}}$;
indeed, the dreamcatchers of $\psi_{2^{\mathbb{N}_{0}}\cup\left(3\mathbb{N}_{0}\right)}\left(z\right)$
and $\psi_{3\mathbb{N}_{0}}\left(z\right)$ are one and the same function,
and are holomorphic on $\mathbb{C}\backslash\frac{1}{3}\mathbb{Z}$.
Obviously, excising the pathological behavior of the Fourier set-series
by passing to their dreamcatchers significantly simplifies the analytical
issues.

Unfortunately---as can be shown by \textbf{Mittag-Leffler's Theorem
(ML'sT)} on the construction of meromorphic functions with specified
poles---$\psi_{V}$ will only admit useful meromorphic dreamcatchers
when $T_{V}$---the set of $\psi_{V}$'s virtual poles---is discrete.
Although \textbf{ML'sT} can be generalized to the case where the function
in question has a \emph{non-discrete}\footnote{This was quite a contentious issue in the late nineteenth century;
Laura Turner gives a fascinating account of the controversy in her
master's thesis (Turner. 2007).}\emph{ }set of poles, in doing so, the non-discrete poles form a natural
boundary for the meromorphic function thus constructed, and as a result,
cannot be studied in isolation like in the discrete case; this distinction
makes all the difference (Turner. 2007).

But now, let us proceed to the advertised $L^{2}$ case. Credit must
be given to the MathOverflow user \textbf{fedja} for unintentially
pointing out the phenomenon (the \textbf{Square Sum Lemma} (\textbf{SSL}))\textbf{
}that makes the contents of this chapter possible. The link to the
MathOverflow question asked by the author and the answer fedja provided
can be found in the Bibliography. While the \emph{proof} of the \textbf{SSL}
is due to fedja, the\emph{ statement} \emph{of the result} is the
author's; fedja's purpose in giving the proof was to show that the
set of edgepoints of a set-series had to be countable. The author,
however, soon realized this result could be re-interpreted as a statement
involving $L^{2}$ norms.

\vphantom{}\textbf{Lemma 2: The Square Sum Lemma (SSL)}: Let $V\subseteq\mathbb{N}_{0}$.
Recalling that $\textrm{P}\left(\psi_{V}\right)$ denotes the set
of edgepoints of $\psi_{V}$ in $\mathbb{R}$, and that $\textrm{P}\left(\psi_{V}\right)$
is a $1$-periodic set, let $\textrm{P}^{\prime}\left(\psi_{V}\right)$
denote:
\[
\textrm{P}^{\prime}\left(\psi_{V}\right)\overset{\textrm{def}}{=}\textrm{P}\left(\psi_{V}\right)\cap\left[0,1\right)
\]
Then, the function $R_{V}\left(x\right)\overset{\textrm{def}}{=}\textrm{VRes}_{x}\left[\psi_{V}\right]$
is in $L^{2}\left(\textrm{P}^{\prime}\left(\psi_{V}\right)\right)$;
that is, $R_{V}$ is a $1$-periodic complex-valued function on $\textrm{P}\left(\psi_{V}\right)$
which is square-integrable with respect to the counting measure on
$\textrm{P}^{\prime}\left(\psi_{V}\right)$:
\[
\sum_{x\in\textrm{P}^{\prime}\left(\psi_{V}\right)}\left|R_{V}\left(x\right)\right|^{2}=\sum_{x\in\textrm{P}^{\prime}\left(\psi_{V}\right)}\left|\textrm{VRes}_{x}\left[\psi_{V}\right]\right|^{2}<\infty
\]

Proof: For brevity, write:
\[
R\left(t\right)\overset{\textrm{def}}{=}\textrm{VRes}_{t}\left[\psi_{V}\right],\textrm{ }\forall t\in\textrm{P}^{\prime}\left(\psi_{V}\right)
\]
Note that:
\begin{eqnarray*}
R\left(t\right) & = & \lim_{y\downarrow0}y\psi_{V}\left(iy+t\right)\\
 & = & \lim_{y\downarrow0}y\sum_{v\in V}e^{2\pi iv\left(iy+t\right)}\\
 & = & \lim_{y\downarrow0}y\sum_{n=0}^{\infty}\mathbf{1}_{V}\left(n\right)e^{2\pi int}e^{-2\pi ny}
\end{eqnarray*}
Now, let:
\[
L\overset{\textrm{def}}{=}\sum_{t\in\textrm{P}^{\prime}\left(\psi_{V}\right)}\left|R\left(t\right)\right|^{2}
\]
Then, applying the \textbf{Cauchy-Schwarz Inequality} for infinite
series, we obtain:
\begin{eqnarray*}
L^{2} & = & \left(\sum_{t\in\textrm{P}^{\prime}\left(\psi_{V}\right)}\left|R\left(t\right)\right|^{2}\right)^{2}\\
 & = & \left(\sum_{t\in\textrm{P}^{\prime}\left(\psi_{V}\right)}\limsup_{y\downarrow0}y\sum_{n=0}^{\infty}\mathbf{1}_{V}\left(n\right)e^{2\pi int}e^{-2\pi ny}\overline{R}\left(t\right)\right)^{2}\\
 & = & \left(\limsup_{y\downarrow0}\sum_{n=0}^{\infty}\left(\sqrt{y}\mathbf{1}_{V}\left(n\right)e^{-\pi ny}\right)\left(\sqrt{y}\sum_{t\in\textrm{P}^{\prime}\left(\psi_{V}\right)}e^{2\pi int}\overline{R}\left(t\right)e^{-\pi ny}\right)\right)^{2}\\
\left(\textrm{CSI}\right); & \leq & \left(\limsup_{y\downarrow0}y\sum_{n=0}^{\infty}\left(\mathbf{1}_{V}\left(n\right)\right)^{2}e^{-2\pi ny}\right)\left(\limsup_{x\downarrow0}x\sum_{n=0}^{\infty}\left|\sum_{t\in\textrm{P}^{\prime}\left(\psi_{V}\right)}\overline{R\left(t\right)}e^{2\pi int}\right|^{2}e^{-2\pi nx}\right)\\
 & = & \frac{\overline{d}\left(V\right)}{2\pi}\limsup_{x\downarrow0}x\sum_{n=0}^{\infty}\left|\sum_{t\in\textrm{P}^{\prime}\left(\psi_{V}\right)}\overline{R\left(t\right)}e^{2\pi int}\right|^{2}e^{-2\pi nx}
\end{eqnarray*}
Next, writing:
\[
\left|\sum_{t\in\textrm{P}^{\prime}\left(\psi_{V}\right)}\overline{R\left(t\right)}e^{2\pi int}\right|^{2}=\sum_{t,\tau\in\textrm{P}^{\prime}\left(\psi_{V}\right)}\overline{R\left(t\right)}R\left(\tau\right)e^{2\pi in\left(t-\tau\right)}
\]
note that:
\begin{eqnarray*}
 & \sum_{n=0}^{\infty}\left|\sum_{t\in\textrm{P}^{\prime}\left(\psi_{V}\right)}\overline{R\left(t\right)}e^{2\pi int}\right|^{2}e^{-2\pi nx}\\
 & \shortparallel\\
 & \sum_{n=0}^{\infty}\left(\sum_{t,\tau\in\textrm{P}^{\prime}\left(\psi_{V}\right)}\overline{R\left(t\right)}R\left(\tau\right)e^{2\pi in\left(t-\tau\right)}\right)e^{-2\pi nx}\\
 & \shortparallel\\
 & \sum_{t,\tau\in\textrm{P}^{\prime}\left(\psi_{V}\right)}\overline{R\left(t\right)}R\left(\tau\right)\sum_{n=0}^{\infty}e^{2\pi in\left(t-\tau\right)}e^{-2\pi nx}
\end{eqnarray*}

For $t\neq\tau$, we have:
\[
\sum_{n=0}^{\infty}e^{2\pi in\left(t-\tau\right)}e^{-2\pi nx}=\frac{1}{1-e^{2\pi i\left(t-\tau\right)}e^{-2\pi x}},\textrm{ }\forall x>0
\]
and thus:
\[
\limsup_{x\downarrow0}x\sum_{n=0}^{\infty}e^{2\pi in\left(t-\tau\right)}e^{-2\pi nx}=\limsup_{x\downarrow0}\frac{x}{1-e^{2\pi i\left(t-\tau\right)}e^{-2\pi x}}=\frac{0}{1-e^{2\pi i\left(t-\tau\right)}}=0
\]
for all $t\neq\tau$. Hence:
\begin{eqnarray*}
 & \limsup_{x\downarrow0}x\sum_{t,\tau\in\textrm{P}^{\prime}\left(\psi_{V}\right)}\overline{R\left(t\right)}R\left(\tau\right)\sum_{n=0}^{\infty}e^{2\pi in\left(t-\tau\right)}e^{-2\pi nx}\\
 & \shortparallel\\
 & \limsup_{x\downarrow0}x\sum_{t\in\textrm{P}^{\prime}\left(\psi_{V}\right)}\left|R\left(t\right)\right|^{2}\sum_{n=0}^{\infty}e^{-2\pi nx}\\
 & \shortparallel\\
 & \sum_{t\in\textrm{P}^{\prime}\left(\psi_{V}\right)}\left|R\left(t\right)\right|^{2}\limsup_{x\downarrow0}\frac{x}{1-e^{-2\pi x}}\\
 & \shortparallel\\
 & \frac{1}{2\pi}\sum_{t\in\textrm{P}^{\prime}\left(\psi_{V}\right)}\left|R\left(t\right)\right|^{2}\\
 & \shortparallel\\
 & \frac{L}{2\pi}
\end{eqnarray*}
As such:
\begin{eqnarray*}
L^{2} & \leq & \frac{\overline{d}\left(V\right)}{2\pi}\limsup_{x\downarrow0}x\sum_{n=0}^{\infty}\left|\sum_{t\in\textrm{P}^{\prime}\left(\psi_{V}\right)}\overline{R\left(t\right)}e^{2\pi int}\right|^{2}e^{-2\pi nx}\\
 & = & \frac{\overline{d}\left(V\right)}{2\pi}\limsup_{x\downarrow0}x\sum_{t,\tau\in\textrm{P}^{\prime}\left(\psi_{V}\right)}\overline{R\left(t\right)}R\left(\tau\right)\sum_{n=0}^{\infty}e^{2\pi in\left(t-\tau\right)}e^{-2\pi nx}\\
 & = & \frac{\overline{d}\left(V\right)}{\left(2\pi\right)^{2}}L
\end{eqnarray*}
which gives:
\[
L\leq\frac{\overline{d}\left(V\right)}{\left(2\pi\right)^{2}}
\]
Thus: 
\[
L=\sum_{t\in\textrm{P}^{\prime}\left(\psi_{V}\right)}\left|R\left(t\right)\right|^{2}
\]
is finite, which shows that $R\in L^{2}\left(\textrm{P}^{\prime}\left(\psi_{V}\right)\right)$,
as desired.

Q.E.D.

\vphantom{}At this point, we can now give the general definition
of a dreamcatcher. Before doing so, however, we must discuss the white
elephant in the room. In its abstract, this paper advertised dreamcatchers
as objects in the Hilbert space $L^{2}\left(\mathbb{Q}/\mathbb{Z}\right)$.
On the other hand, the \textbf{SSL} works with $\textrm{P}^{\prime}\left(\psi_{V}\right)$,
the set of representatives for the equivalence class modulo $1$ of
$\textrm{P}\left(\psi_{V}\right)$, the set of $\psi_{V}$'s edgepoints.
As discussed in §0.2, the author's method of twice-removing the dynamics
of hydra maps from the original context is complicated by certain
technical difficulties. The above discrepency between $L^{2}\left(\mathbb{Q}/\mathbb{Z}\right)$
and $L^{2}\left(\textrm{P}^{\prime}\left(\psi_{V}\right)\right)$
is the first of those difficulties.

The convergence of the square-magnitude sum in the \textbf{SSL }shows
that the subset of $\textrm{P}^{\prime}\left(\psi_{V}\right)$ (and
thus, of $\textrm{P}\left(\psi_{V}\right)$) on which $R$ is non-vanishing
must be countable; that is, $R$'s support in $\textrm{P}\left(\psi_{V}\right)$
is countable. However, it says nothing about the \emph{location} of
this support within $\textrm{P}\left(\psi_{V}\right)$, nor about
the location of $\textrm{P}\left(\psi_{V}\right)$ within $\mathbb{R}$.
Obviously, this creates problems. The arguments to be made in Chapter
3 depend upon the algebraic and arithmetical properties of $\mathbb{Q}/\mathbb{Z}$;
as a discrete abelian group, it possesses a \textbf{Pontryagin Dual},
and the theory \textbf{Pontryagin Duality }can (and shall) be used
to formulate the Fourier Transform on $L^{2}\left(\mathbb{Q}/\mathbb{Z}\right)$.
Thus, on the one hand, if $\textrm{P}\left(\psi_{V}\right)$ contains
all of $\mathbb{Q}$, we can then realize $R$ as an element of $L^{2}\left(\mathbb{Q}/\mathbb{Z}\right)$
simply by restricting its domain of definition to $\mathbb{Q}$. On
the other hand, if there are points in $\mathbb{Q}$ that lie outside
of $\textrm{P}\left(\psi_{V}\right)$, the only way we can obtain
arithmetical structure from $L^{2}\left(\textrm{P}^{\prime}\left(\psi_{V}\right)\right)$
would be if, by some miracle, $\textrm{P}\left(\psi_{V}\right)$ turns
out to contain a subset $T$ which, when equipped with addition, becomes
a subgroup of $\mathbb{Q}/\mathbb{Z}$ (for example, $T=\mathbb{Z}\left[\frac{1}{a}\right]$,
where $a\in\mathbb{N}_{2}$) (or, more generally, of the topological
group obtained by equipping $\left(\mathbb{R}/\mathbb{Z},+\right)$
with the discrete topology), the Pontryagin duality arguments will
still apply. And---the reader should keep this in mind---this is
only the \emph{first} wave of technical complications; the second
wave lies in wait just beyond the next page break.

The take-away from all this is that, at least for the remainder of
this chapter, we will need to proceed somewhat gingerly, slowly modulating
our language to enable us to speak of dreamcatchers without reference
to the functions from which they came, and the burden of the accompanying
technicalities. To that end, in giving the rigorous definition of
a dreamcatcher, we will shift our weight from $\psi$ to $\textrm{P}\left(\psi\right)$---or,
rather, to subsets of $\textrm{P}\left(\psi\right)$.

\vphantom{}\textbf{Definition 15}:

I. Given a function $\psi\in\mathcal{A}\left(\mathbb{H}_{+i}\right)$
and a $1$-periodic set\footnote{We will write $T\subseteq S/\mathbb{Z}$ to denote that $T$ is a
$1$-periodic subset of $S$, where $S$ is any $1$-periodic set.} $T\subseteq\textrm{P}\left(\psi\right)$, $\psi$ is said to \textbf{admit
a dreamcatcher over $T$} whenever the function:
\[
x\in T\mapsto\textrm{VRes}_{x}\left[\psi\right]\overset{\textrm{def}}{=}\lim_{y\downarrow0}y\psi\left(x+iy\right)\in\mathbb{C}
\]
is an element of $L^{2}\left(T/\mathbb{Z}\right)$. In that case,
$\textrm{VRes}_{x}\left[\psi\right]$ is called the \textbf{dreamcatcher
of $\psi$ over $T$}, and is denoted $R\left(x\right)$---or $R_{\psi}\left(x\right)$,
or $R_{T}\left(x\right)$, or some other notational modification,
depending on situational need; for example, if $\psi$ is the Fourier
set-series $\psi=\psi_{V}$ of a set $V\subseteq\mathbb{N}_{0}$,
the dreamcatcher might also be denoted as $R_{V}\left(x\right)$,
and we will often speak of the \textbf{dreamcatcher of $V$} when
referring to the dreamcatcher of $\psi_{V}$.

II. Given an arbitrary \emph{subgroup }$T$ of $\mathbb{R}/\mathbb{Z}$,
we write $\mathcal{D}_{1}\left(T\right)$ to denote the set of all
$\psi\in\mathcal{A}\left(\mathbb{H}_{+i}\right)$ which admit a dreamcatcher
over $T$; that is:
\[
\mathcal{D}_{1}\left(T\right)\overset{\textrm{def}}{=}\left\{ \psi\in\mathcal{A}\left(\mathbb{H}_{+i}\right):\left(x\in T\mapsto\textrm{VRes}_{x}\left[\psi\right]\in\mathbb{C}\right)\in L^{2}\left(T/\mathbb{Z}\right)\right\} 
\]
Note that $T\subseteq\textrm{P}\left(\psi\right)$ for all $\psi\in\mathcal{D}_{1}\left(T\right)$.

\emph{Remark}: The subscript $1$ is the author's convention for distinguishing
the dreamcatchers defined here from the generalized form of dreamcatchers
(``subdominant dreamcatchers'') corresponding to boundary behavior
of order strictly less than $O\left(\frac{1}{y}\right)$.

\pagebreak{}

\subsection{Singularity Conservation Laws for Hydra Maps}

With the terminology of the previous section in hand, we can now give
rigorous footing to the idea of a \textbf{Singularity Conservation
Law}. In keeping with this Chapter's theme of taking virtual residues
to pass from complicated objects to hopefully-less-complicated objects,
having done so with set-series in §2.3, the next order of business
is to do the same to the (Fourier) permutation operator induced by
$H$. Doing so brings us to the second of our two main technical obstacles,
and this time, $H$ itself is directly implicated.

\vphantom{}\textbf{Definition} \textbf{16}: Let $H:\mathbb{N}_{0}\rightarrow\mathbb{N}_{0}$
be a $\varrho$-hydra map. A set $T\subseteq\mathbb{R}$ is said to
be \textbf{$H$-branch invariant} whenever:

\[
\frac{\varrho T+k}{\mu_{j}}\subseteq T,\textrm{ }\forall j\in\left\{ 0,\ldots,\varrho-1\right\} ,\forall k\in\left\{ 0,\ldots,\mu_{j}-1\right\} 
\]
meaning that $\frac{\varrho t+k}{\mu_{j}}\subseteq T$ for all $t\in T$
and all values of $j$ and $k$ as indicated above.

\vphantom{}The issues arising from $H$-branch invariance are two-fold.
Firstly, without the hypothesis that $T$ is $H$-branch invariant,
the interchange of limit and finite sum required obtain the right-hand
side of \ref{eq:VRF} in \textbf{Lemma 3} fails to be justified. The
second issue, on the other hand foreshadows our imminent change of
setting from $\mathcal{A}\left(\mathbb{H}_{+i}\right)$ to $L^{2}\left(T/\mathbb{Z}\right)$,
and from there up to $L^{2}\left(\mathbb{Q}/\mathbb{Z}\right)$; but
discussing that is best left until after our next lemma.

\vphantom{}\textbf{Lemma 3 (Virtual Residue Formula (VRF))}: Let
$H:\mathbb{N}_{0}\rightarrow\mathbb{N}_{0}$ be a $\varrho$-hydra
map, and let $T\subseteq\mathbb{R}/\mathbb{Z}$ be $H$-branch invariant.
Then, for every $\psi\in\mathcal{D}_{2}\left(T\right)$:

\begin{equation}
\textrm{VRes}_{t}\left[\mathscr{Q}_{H}\left\{ \psi\right\} \right]=\frac{1}{\varrho}\sum_{j=0}^{\varrho-1}e^{-2\pi i\frac{b_{j}}{a_{j}}t}\sum_{k=0}^{\mu_{j}-1}\xi_{\mu_{j}}^{-kH\left(j\right)}\textrm{VRes}_{\frac{\varrho t+k}{\mu_{j}}}\left[\psi\right],\textrm{ }\forall t\in T\label{eq:VRF}
\end{equation}

Proof: Let $H$, $T$ and $\psi$ be as given. The standard branch
formula for $\mathscr{Q}_{H}\left\{ \psi\right\} \left(z\right)$
is:
\begin{equation}
\mathscr{Q}_{H}\left\{ \psi\right\} \left(z\right)=\sum_{j=0}^{\varrho-1}\frac{e^{-2\pi i\frac{b_{j}}{a_{j}}z}}{\mu_{j}}\sum_{k=0}^{\mu_{j}-1}e^{-\frac{2k\pi iH\left(j\right)}{\mu_{j}}}\psi\left(\frac{\varrho z+k}{\mu_{j}}\right)
\end{equation}
Now, fix an arbitrary $t\in T$. Then, the virtual residue of $\mathscr{Q}_{H}\left\{ \psi\right\} \left(z\right)$
at $t$ is given by:

\begin{eqnarray*}
 & \textrm{VRes}_{t}\left[\mathscr{Q}_{H}\left\{ \psi\right\} \right]\\
 & \shortparallel\\
 & \textrm{VRes}_{t}\left[\sum_{j=0}^{\varrho-1}\frac{e^{-2\pi i\frac{b_{j}}{a_{j}}z}}{\mu_{j}}\sum_{k=0}^{\mu_{j}-1}e^{-\frac{2k\pi iH\left(j\right)}{\mu_{j}}}\psi\left(\frac{\varrho z+k}{\mu_{j}}\right)\right]\\
 & \shortparallel\\
 & \sum_{j=0}^{\varrho-1}\frac{e^{-2\pi i\frac{b_{j}}{a_{j}}t}}{\mu_{j}}\sum_{k=0}^{\mu_{j}-1}e^{-\frac{2k\pi iH\left(j\right)}{\mu_{j}}}\lim_{z\rightarrow t}\left(z-t\right)\psi\left(\frac{\varrho z+k}{\mu_{j}}\right)\\
 & \shortparallel\\
 & \sum_{j=0}^{\varrho-1}\frac{e^{-2\pi i\frac{b_{j}}{a_{j}}t}}{\mu_{j}}\sum_{k=0}^{\mu_{j}-1}e^{-\frac{2k\pi iH\left(j\right)}{\mu_{j}}}\textrm{VRes}_{t}\left[\psi\left(\frac{\varrho z+k}{\mu_{j}}\right)\right]
\end{eqnarray*}
That being done, all that remains is to compute the virtual residue
of $\psi\left(\frac{\varrho z+k}{\mu_{j}}\right)$ at $z=t$. For
brevity, write:
\begin{equation}
A_{j,k}\left(z\right)\overset{\textrm{def}}{=}\frac{\varrho z+k}{\mu_{j}},\textrm{ }\forall j\in\left\{ 0,\ldots,\varrho-1\right\} ,\forall k\in\left\{ 0,\ldots,\mu_{j}-1\right\} 
\end{equation}
\begin{equation}
A_{j,k}^{-1}\left(z\right)\overset{\textrm{def}}{=}\frac{\mu_{j}z-k}{\varrho},\textrm{ }\forall j\in\left\{ 0,\ldots,\varrho-1\right\} ,\forall k\in\left\{ 0,\ldots,\mu_{j}-1\right\} 
\end{equation}
(where $A_{j,k}^{-1}$ is the inverse of $A_{j,k}$) to denote the
branches of $H$ and their inverses. Then, make a change of variables
in the virtual-residue-defining limit, and multiply and divide by
an auxiliary factor in preparation for splitting the limit into two
parts:
\begin{eqnarray*}
\textrm{VRes}_{t}\left[\psi\left(\frac{\varrho z+k}{\mu_{j}}\right)\right] & = & \lim_{z\rightarrow t}\left(z-t\right)\psi\left(A_{j,k}\left(z\right)\right)\\
\left(\textrm{let }w=A_{j,k}\left(z\right)\right); & = & \lim_{w\rightarrow A_{j,k}\left(t\right)}\left(A_{j,k}^{-1}\left(w\right)-t\right)\psi\left(w\right)\\
\left(\times\textrm{ \& }\div\textrm{ by }w-A_{j,k}\left(t\right)\right); & = & \lim_{w\rightarrow A_{j,k}\left(t\right)}\frac{A_{j,k}^{-1}\left(w\right)-t}{w-A_{j,k}\left(t\right)}\left(w-A_{j,k}\left(t\right)\right)\psi\left(w\right)
\end{eqnarray*}

Since:
\[
\left(A_{j,k}^{-1}\right)^{\prime}\left(z\right)=\frac{\mu_{j}}{\varrho}
\]
an application of L'Hôpital's Rule gives:
\begin{equation}
\lim_{w\rightarrow A_{j,k}\left(t\right)}\frac{A_{j,k}^{-1}\left(w\right)-t}{w-A_{j,k}\left(t\right)}=\lim_{w\rightarrow A_{j,k}\left(t\right)}\frac{\frac{\mu_{j}}{\varrho}-0}{1-0}=\frac{\mu_{j}}{\varrho}
\end{equation}
Meanwhile:

\begin{equation}
\lim_{w\rightarrow A_{j,k}\left(t\right)}\left(w-A_{j,k}\left(t\right)\right)F\left(w\right)=\textrm{VRes}_{A_{j,k}\left(t\right)}\left[\psi\right]=\textrm{VRes}_{\frac{\varrho t+k}{\mu_{j}}}\left[\psi\right]
\end{equation}
Note that the $H$-branch invariance of $T$ guarantees that $\textrm{VRes}_{\frac{\varrho t+k}{\mu_{j}}}\left[\psi\right]$
exists for all $t\in T$, since $\frac{\varrho t+k}{\mu_{j}}\in T\subseteq\textrm{P}\left(\psi\right)$
for all $\psi\in\mathcal{D}_{2}\left(T\right)$.

Thus, our limit becomes:
\begin{eqnarray*}
\lim_{z\rightarrow t}\left(z-t\right)\psi\left(A_{j,k}\left(z\right)\right) & = & \lim_{x\rightarrow A_{j,k}\left(t\right)}\frac{A_{j,k}^{-1}\left(x\right)-t}{x-A_{j,k}\left(t\right)}\times\lim_{y\rightarrow A_{j,k}\left(t\right)}\left(y-A_{j,k}\left(t\right)\right)\psi\left(y\right)\\
 & = & \frac{\mu_{j}}{\varrho}\textrm{VRes}_{A_{j,k}\left(t\right)}\left[\psi\right]
\end{eqnarray*}
So:
\begin{equation}
\textrm{VRes}_{t}\left[\psi\left(\frac{\varrho z+k}{\mu_{j}}\right)\right]=\frac{\mu_{j}}{\varrho}\textrm{VRes}_{A_{j,k}\left(t\right)}\left[\psi\right]=\frac{\mu_{j}}{\varrho}\textrm{VRes}_{\frac{\varrho t+k}{\mu_{j}}}\left[\psi\right]
\end{equation}
Plugging these virtual residue computations into the formula for $\textrm{VRes}_{t}\left[\mathscr{Q}_{H}\left\{ \psi\right\} \right]$
produces:
\[
\textrm{VRes}_{t}\left[\mathscr{Q}_{H}\left\{ \psi\right\} \right]=\frac{1}{\varrho}\sum_{j=0}^{\varrho-1}e^{-2\pi i\frac{b_{j}}{a_{j}}t}\sum_{k=0}^{\mu_{j}-1}\xi_{\mu_{j}}^{-kH\left(j\right)}\textrm{VRes}_{\frac{\varrho t+k}{\mu_{j}}}\left[\psi\right],\textrm{ }\forall t\in T
\]
as desired.

Q.E.D.

\vphantom{}Note, in terms of $R_{\psi}\left(t\right)=\textrm{VRes}_{t}\left[\psi\right]$,
this can be written as:
\[
\textrm{VRes}_{t}\left[\mathscr{Q}_{H}\left\{ \psi\right\} \right]=\frac{1}{\varrho}\sum_{j=0}^{\varrho-1}e^{-2\pi i\frac{b_{j}}{a_{j}}t}\sum_{k=0}^{\mu_{j}-1}\xi_{\mu_{j}}^{-kH\left(j\right)}R_{\psi}\left(\frac{\varrho t+k}{\mu_{j}}\right),\textrm{ }\forall t\in T
\]

Now to construct an operator out of this. As we do so, the reader
should note the role played by the requirement that $T$ be $H$-branch
invariant.

\vphantom{}\textbf{Proposition} \textbf{10}: Let $H:\mathbb{N}_{0}\rightarrow\mathbb{N}_{0}$
be a $\varrho$-hydra map, and let $T\subseteq\mathbb{R}/\mathbb{Z}$
be $H$-branch invariant. Then, for every $f\in L^{2}\left(T/\mathbb{Z}\right)$
the function:
\[
x\mapsto\frac{1}{\varrho}\sum_{j=0}^{\varrho-1}e^{-2\pi i\frac{b_{j}}{a_{j}}x}\sum_{k=0}^{\mu_{j}-1}\xi_{\mu_{j}}^{-kH\left(j\right)}f\left(\frac{\varrho x+k}{\mu_{j}}\right)
\]
is an element of $L^{2}\left(T/\mathbb{Z}\right)$.

Proof: Let $f\in L^{2}\left(T/\mathbb{Z}\right)$. Then, $f$ can
be viewed as a $1$-periodic element of $L_{\textrm{loc}}^{2}\left(T\right)$
(read {[}$L^{2}$ loak of $T${]}), the space of ``locally square-integrable''
functions $g:T\rightarrow\mathbb{C}$; that is, functions $g:T\rightarrow\mathbb{C}$
so so that:
\[
\sum_{x\in I}\left|g\left(x\right)\right|^{2}<\infty
\]
for every bounded subset $I\subseteq T$. Since $T$ is given to be
invariant under the maps $x\mapsto\frac{\varrho x+k}{\mu_{j}}$, and
since these maps are bijections (seeing as all affine linear maps
on $\mathbb{R}$ are necessarily bijective) for all $j$ and $k$,
$f\in L_{\textrm{loc}}^{2}\left(T\right)$ implies that the functions
$f\left(\frac{\varrho x+k}{\mu_{j}}\right)$ are also in $L_{\textrm{loc}}^{2}\left(T\right)$
for all $j$ and $k$. Since:
\[
\left|e^{-2\pi i\frac{b_{j}}{a_{j}}x}f\left(\frac{\varrho x+k}{\mu_{j}}\right)\right|=\left|f\left(\frac{\varrho x+k}{\mu_{j}}\right)\right|
\]
the function $e^{-2\pi i\frac{b_{j}}{a_{j}}x}f\left(\frac{\varrho x+k}{\mu_{j}}\right)$
is then in $L_{\textrm{loc}}^{2}\left(T\right)$ for all $j,k$. As
such, the function: 
\[
\frac{1}{\varrho}\sum_{j=0}^{\varrho-1}e^{-2\pi i\frac{b_{j}}{a_{j}}x}\sum_{k=0}^{\mu_{j}-1}\xi_{\mu_{j}}^{-kH\left(j\right)}f\left(\frac{\varrho x+k}{\mu_{j}}\right)
\]
is a $\mathbb{C}$-linear combination of functions in $L_{\textrm{loc}}^{2}\left(T\right)$,
and is therefore itself an element of $L_{\textrm{loc}}^{2}\left(T\right)$.

Now, to show that this function is actually in $L^{2}\left(T/\mathbb{Z}\right)$,
all we need to do is prove that it is $1$-periodic, because, then---since
$\left[0,1\right)_{T}=\left[0,1\right)\cap T$ is a bounded subset
of $T$---we will have show that the sum of its square-magnitudes
over $x\in\left[0,1\right)_{T}$ will be finite. To see the periodicity,
we simply replace $x$ with $x+1$. We can make the computations cleaner
by using the temporary notation:
\[
f_{j}\left(x\right)\overset{\textrm{def}}{=}\sum_{k=0}^{\mu_{j}-1}\xi_{\mu_{j}}^{-kH\left(j\right)}f\left(\frac{\varrho x+k}{\mu_{j}}\right)
\]
Then:
\begin{eqnarray*}
f_{j}\left(x+1\right) & = & \sum_{k=0}^{\mu_{j}-1}\xi_{\mu_{j}}^{-kH\left(j\right)}f\left(\frac{\varrho x+\varrho+k}{\mu_{j}}\right)\\
\left(\kappa=\varrho+k\right); & = & \sum_{k=0}^{\mu_{j}-1}\xi_{\mu_{j}}^{-\left(\kappa-\varrho\right)H\left(j\right)}f\left(\frac{\varrho x+\kappa}{\mu_{j}}\right)\\
 & = & \xi_{\mu_{j}}^{-\varrho H\left(j\right)}\underbrace{\sum_{k=0}^{\mu_{j}-1}\xi_{\mu_{j}}^{-\kappa H\left(j\right)}f\left(\frac{\varrho x+\kappa}{\mu_{j}}\right)}_{f_{j}\left(x\right)}\\
 & = & \xi_{\mu_{j}}^{-\varrho H\left(j\right)}f_{j}\left(x\right)
\end{eqnarray*}
As such:
\begin{eqnarray*}
\frac{1}{\varrho}\sum_{j=0}^{\varrho-1}e^{-2\pi i\frac{b_{j}}{a_{j}}\left(x+1\right)}\sum_{k=0}^{\mu_{j}-1}\xi_{\mu_{j}}^{-kH\left(j\right)}f\left(\frac{\varrho\left(x+1\right)+k}{\mu_{j}}\right) & = & \frac{1}{\varrho}\sum_{j=0}^{\varrho-1}e^{-2\pi i\frac{b_{j}}{a_{j}}\left(x+1\right)}f_{j}\left(x+1\right)\\
 & = & \frac{1}{\varrho}\sum_{j=0}^{\varrho-1}e^{-2\pi i\frac{b_{j}}{a_{j}}\left(x+1\right)}\xi_{\mu_{j}}^{-\varrho H\left(j\right)}f_{j}\left(x\right)
\end{eqnarray*}
Cleaning up the exponentials:
\begin{eqnarray*}
e^{-2\pi i\frac{b_{j}}{a_{j}}\left(x+1\right)}\xi_{\mu_{j}}^{-\varrho H\left(j\right)} & = & e^{-2\pi i\frac{b_{j}}{a_{j}}x}e^{-2\pi i\frac{b_{j}}{a_{j}}}e^{-\frac{2\pi i}{\mu_{j}}\varrho H\left(j\right)}\\
 & = & e^{-2\pi i\frac{b_{j}}{a_{j}}x}e^{-2\pi i\frac{b_{j}}{a_{j}}}e^{-\frac{2\pi i}{\frac{\varrho a_{j}}{d_{j}}}\varrho\frac{ja_{j}+b_{j}}{d_{j}}}\\
 & = & e^{-2\pi i\frac{b_{j}}{a_{j}}x}e^{-2\pi i\frac{b_{j}}{a_{j}}}e^{-2\pi i\frac{ja_{j}+b_{j}}{a_{j}}}\\
 & = & e^{-2\pi i\frac{b_{j}}{a_{j}}x}e^{-2\pi i\frac{ja_{j}}{a_{j}}}\\
 & = & e^{-2\pi i\frac{b_{j}}{a_{j}}x}\underbrace{e^{-2\pi ij}}_{1}\\
\left(j\in\mathbb{Z}\right); & = & e^{-2\pi i\frac{b_{j}}{a_{j}}x}
\end{eqnarray*}
we have that: 
\begin{eqnarray*}
\frac{1}{\varrho}\sum_{j=0}^{\varrho-1}e^{-2\pi i\frac{b_{j}}{a_{j}}\left(x+1\right)}f_{j}\left(x+1\right) & = & \frac{1}{\varrho}\sum_{j=0}^{\varrho-1}e^{-2\pi i\frac{b_{j}}{a_{j}}\left(x+1\right)}\xi_{\mu_{j}}^{-\varrho H\left(j\right)}f_{j}\left(x\right)\\
 & = & \frac{1}{\varrho}\sum_{j=0}^{\varrho-1}e^{-2\pi i\frac{b_{j}}{a_{j}}x}f_{j}\left(x\right)
\end{eqnarray*}
and so, our function is indeed $1$-periodic, which shows that it
is thus an element of $L^{2}\left(T/\mathbb{Z}\right)$.

Q.E.D.

\vphantom{}This then motivates the definition of the virtual analogue
of $\mathscr{Q}_{H}$, with which we can then directly formulate the
Singularity Conservation Law.

\vphantom{}\textbf{Definition} \textbf{17}: Let $H:\mathbb{N}_{0}\rightarrow\mathbb{N}_{0}$
be a $\varrho$-hydra map, and let $T\subseteq\mathbb{R}/\mathbb{Z}$
be $H$-branch invariant. Then, we define $\mathscr{\mathfrak{Q}}_{H:T}:L^{2}\left(T/\mathbb{Z}\right)\rightarrow L^{2}\left(T/\mathbb{Z}\right)$,
the \textbf{Dreamcatcher Operator induced by $H$ over $T$} as: 
\begin{equation}
\mathfrak{Q}_{H:T}\left\{ f\right\} \left(x\right)\overset{\textrm{def}}{=}\frac{1}{\varrho}\sum_{j=0}^{\varrho-1}\sum_{k=0}^{\mu_{j}-1}\xi_{\mu_{j}}^{-kH\left(j\right)}e^{-2\pi i\frac{b_{j}}{a_{j}}x}f\left(\frac{\varrho x+k}{\mu_{j}}\right),\textrm{ }\forall f\in L^{2}\left(T/\mathbb{Z}\right)\label{eq:Def of Q_H:T}
\end{equation}

\emph{Remark}s:

i. Note that \textbf{Proposition 10} guarantees that $\mathfrak{Q}_{H:T}$
is well-defined as a linear operator on $L^{2}\left(T/\mathbb{Z}\right)$.

ii. We write $\mathfrak{Q}_{H}$ to denote the case where $T=\mathbb{Q}$.

\vphantom{}Now we are properly positioned to discuss the second reason
for invoking $H$-branch invariance. To begin with, if---for the
moment---we completely forget everything we have done prior to \ref{eq:Def of Q_H:T},
and simply consider the formula: 
\[
\mathfrak{Q}_{H}\left\{ f\right\} \left(x\right)=\frac{1}{\varrho}\sum_{j=0}^{\varrho-1}\sum_{k=0}^{\mu_{j}-1}\xi_{\mu_{j}}^{-kH\left(j\right)}e^{-2\pi i\frac{b_{j}}{a_{j}}x}f\left(\frac{\varrho x+k}{\mu_{j}}\right)
\]
where $f\in L^{2}\left(\mathbb{Q}/\mathbb{Z}\right)$, everything
works out perfectly fine. Since $\frac{\varrho\mathbb{Q}+k}{\mu_{j}}=\mathbb{Q}$,
the operator $\mathfrak{Q}_{H}$ sends every $f\in L^{2}\left(\mathbb{Q}/\mathbb{Z}\right)$
to another function in $L^{2}\left(\mathbb{Q}/\mathbb{Z}\right)$.
On the other hand, in the proof of \textbf{Proposition 10}, we saw
that the pre-composition by the maps $\frac{\varrho x+k}{\mu_{j}}$
by an application of $\mathfrak{Q}_{H:T}$ to a function $f\in L^{2}\left(T/\mathbb{Z}\right)$
might cause $\mathfrak{Q}_{H}\left\{ f\right\} \notin L^{2}\left(T/\mathbb{Z}\right)$
to occur if the $H$-branch invariance of $T$ was not guaranteed.
In other words, viewing $\mathfrak{Q}_{H}$ as an operator on $L^{2}\left(\mathbb{Q}/\mathbb{Z}\right)$,
\emph{we use the $H$-branch invariance of $T$ to guarantee that}
$L^{2}\left(T/\mathbb{Z}\right)$\emph{ is} \emph{an }\textbf{invariant
subspace}\emph{ of $\mathfrak{Q}_{H}$}.

\vphantom{}\textbf{Proposition 11 (The Extension Proposition}): Let
$H:\mathbb{N}_{0}\rightarrow\mathbb{N}_{0}$ be a $\varrho$-hydra
map, let $S,T\subseteq\mathbb{R}/\mathbb{Z}$ be $H$-branch invariant
sets such that $T\subseteq S$. Then: 
\[
\mathfrak{Q}_{H:S}\mid_{L^{2}\left(T/\mathbb{Z}\right)}=\mathfrak{Q}_{H:T}
\]
that is to say, the restriction of $\mathfrak{Q}_{H:S}$ to $L^{2}\left(T/\mathbb{Z}\right)$
is the same as $\mathfrak{Q}_{H:T}$.

Proof: Since $T\subseteq S$, $L^{2}\left(T/\mathbb{Z}\right)$ is
a Hilbert subspace of $L^{2}\left(S/\mathbb{Z}\right)$. Since the
defining formulae for $\mathfrak{Q}_{H:S}$ and $\mathfrak{Q}_{H:T}$
are formally identical---they differ only in their respective domains
of definition---the fact that $T$ is $H$-branch invariant shows
that $L^{2}\left(T/\mathbb{Z}\right)$ is an invariant subspace of
$\mathfrak{Q}_{H:S}$ (i.e., $\mathfrak{Q}_{H:S}\left(L^{2}\left(T/\mathbb{Z}\right)\right)\subseteq L^{2}\left(T/\mathbb{Z}\right)$),
and thus, that $\mathfrak{Q}_{H:S}\mid_{L^{2}\left(T/\mathbb{Z}\right)}=\mathfrak{Q}_{H:T}$
holds.

Q.E.D.

\vphantom{}Now, the \textbf{Singularity Conservation Law}:

\vphantom{}\textbf{Theorem 4: The Singularity Conservation Law for
$\varrho$-Hydra Maps }(\textbf{SCL}): Let $H:\mathbb{N}_{0}\rightarrow\mathbb{N}_{0}$
be a $\varrho$-hydra map, let $T\subseteq\mathbb{R}/\mathbb{Z}$
be $H$-branch invariant, let $\psi\in\mathcal{D}_{2}\left(T\right)$,
and write $R$ to denote the $T$-dreamcatcher of $\psi$. Then:
\[
\psi\in\textrm{Ker}\left(1-\mathscr{Q}_{H}\right)\Rightarrow R\in\textrm{Ker}\left(1-\mathfrak{Q}_{H:T}\right)
\]
That is to say, if $\psi$ is fixed by $\mathscr{Q}_{H}$, then $\psi$'s
dreamcatcher is fixed by $\mathfrak{Q}_{H:T}$.

Proof: Let everything be as given. Then, by the \textbf{VRF} \ref{eq:VRF},
it follows that:
\begin{eqnarray*}
\textrm{VRes}_{t}\left[\mathscr{Q}_{H}\left\{ \psi\right\} \right] & = & \frac{1}{\varrho}\sum_{j=0}^{\varrho-1}e^{-2\pi i\frac{b_{j}}{a_{j}}t}\sum_{k=0}^{\mu_{j}-1}\xi_{\mu_{j}}^{-kH\left(j\right)}\textrm{VRes}_{\frac{\varrho t+k}{\mu_{j}}}\left[\psi\right],\textrm{ }\forall t\in T\\
\left(\psi\textrm{ is fixed by }\mathscr{Q}_{H}\right); & \Updownarrow\\
\underbrace{\textrm{VRes}_{t}\left[\psi\right]}_{R\left(t\right)} & = & \frac{1}{\varrho}\sum_{j=0}^{\varrho-1}e^{-2\pi i\frac{b_{j}}{a_{j}}t}\sum_{k=0}^{\mu_{j}-1}\xi_{\mu_{j}}^{-kH\left(j\right)}\underbrace{\textrm{VRes}_{\frac{\varrho t+k}{\mu_{j}}}\left[\psi\right]}_{R\left(\frac{\varrho t+k}{\mu_{j}}\right)},\textrm{ }\forall t\in T\\
 & \Updownarrow\\
R\left(t\right) & = & \underbrace{\frac{1}{\varrho}\sum_{j=0}^{\varrho-1}e^{-2\pi i\frac{b_{j}}{a_{j}}t}\sum_{k=0}^{\mu_{j}-1}\xi_{\mu_{j}}^{-kH\left(j\right)}R\left(\frac{\varrho t+k}{\mu_{j}}\right)}_{\mathfrak{Q}_{H:T}\left\{ R\right\} \left(t\right)},\textrm{ }\forall t\in T\\
 & \Updownarrow\\
R\left(t\right) & = & \mathfrak{Q}_{H:T}\left\{ R\right\} \left(t\right),\textrm{ }\forall t\in T\\
 & \Updownarrow\\
R & = & \mathfrak{Q}_{H:T}\left\{ R\right\} \textrm{ in }L^{2}\left(T/\mathbb{Z}\right)
\end{eqnarray*}

Q.E.D.

\vphantom{}Consequently, anything we can prove about the fixed points
of $\mathfrak{Q}_{H:T}$ will then establish necessary conditions
on the dreamcatchers of $\mathscr{Q}_{H}$'s fixed points, such as
the set-series for an orbit class of $H$. As was seen in §1.1.2,
the limits used to define the dreamcatcher of a set-series correspond
to the densities and growth asymptotics of the underlying set and
its various cross-sections. So, by chracterizing $\textrm{Ker}\left(1-\mathfrak{Q}_{H:T}\right)$,
we can obtain meaningful information about the orbit classes of $H$.
That being said, the study of $\textrm{Ker}\left(1-\mathfrak{Q}_{H:T}\right)$
becomes terribly complicated whenever we do not have $L^{2}\left(T/\mathbb{Z}\right)$'s
guarantee that it will be an invariant subspace of $\mathfrak{Q}_{H}$.

Our next result---the last of this section, and of the chapter as
a whole---puts a capstone on our current discussion by detailing
how to extend the examination of $\mathfrak{Q}_{H}$'s fixed points
over $L^{2}\left(T/\mathbb{Z}\right)$ to a larger subspace of $L^{2}\left(\mathbb{R}/\mathbb{Z}\right)$.

\vphantom{}\textbf{Corollary 3: (Extended SCL} (\textbf{ESCL})):
Let $H:\mathbb{N}_{0}\rightarrow\mathbb{N}_{0}$ be a $\varrho$-hydra
map, let $T,S\subseteq\mathbb{R}/\mathbb{Z}$ be $H$-branch invariant
sets for which $T\subseteq S$ and for which $S\backslash T$ is also
$H$-branch invariant, and let $\psi\in\mathcal{D}_{1}\left(T\right)$.
Writing $R$ to denote:
\[
R\left(x\right)\overset{\textrm{def}}{=}\begin{cases}
\textrm{VRes}_{x}\left[\psi\right] & \textrm{if }x\in T\\
0 & \textrm{if }x\in S\backslash T
\end{cases}
\]
an extension of $\psi$'s dreamcatcher to $L^{2}\left(S/\mathbb{Z}\right)$,
it follows that:
\[
\psi\in\textrm{Ker}\left(1-\mathscr{Q}_{H}\right)\Rightarrow R\in\textrm{Ker}\left(1-\mathfrak{Q}_{H:S}\right)
\]
That is to say, if $\psi$ is fixed by $\mathscr{Q}_{H}$, then the
dreamcatcher of $\psi$ over $T$ has an extension to $L^{2}\left(S/\mathbb{Z}\right)$
which fixed by $\mathfrak{Q}_{H:S}$, for every $S\supseteq T$ satisfying
the above properties.

Proof: Let $H,T,S,\psi$, and $R$ be as given. If $\psi$ is fixed
by $\mathscr{Q}_{H}$, then $R\mid_{T}$---the restriction of $R$
to $T$---is the dreamcatcher of $\psi$ over $T$, and so---by
the \textbf{SCL}---$R\mid_{T}$ is then fixed by $\mathfrak{Q}_{H:T}$.
Since $L^{2}\left(S/\mathbb{Z}\right)\cong L^{2}\left(T/\mathbb{Z}\right)\oplus L^{2}\left(\left(S\backslash T\right)/\mathbb{Z}\right)$,
and since $\mathfrak{Q}_{H:S}$ is linear, it follows by the \textbf{Extension
Proposition }and the $H$-branch invariance of $S$, $S\backslash T$,
and $T$ that:
\begin{eqnarray*}
\mathfrak{Q}_{H:S}\left\{ R\right\}  & = & \underbrace{\mathfrak{Q}_{H:T}\left\{ R\mid_{T}\right\} }_{=R\mid_{T},\textrm{ by \textbf{SCL-}\ensuremath{\mathbb{Z}}}}\oplus\mathfrak{Q}_{H:S\backslash T}\left\{ R\mid_{S\backslash T}\right\} \\
\left(R\mid_{S\backslash T}\textrm{ is identically }0\right); & = & R\mid_{T}\oplus\mathfrak{Q}_{H:S\backslash T}\left\{ 0\right\} \\
 & = & R\mid_{T}\oplus0\\
\left(R\mid_{S\backslash T}\textrm{ is identically }0\right); & = & R\mid_{T}\oplus R\mid_{S\backslash T}\\
 & = & R
\end{eqnarray*}
and thus, that $R\in\textrm{Ker}\left(1-\mathfrak{Q}_{H:S}\right)$,
as desired.

Q.E.D.

\vphantom{}Finally, we have some terminology for those functions
for which our analyses are free of complications.

\vphantom{}\textbf{Definition 18}: A function $\psi\in\mathcal{A}\left(\mathbb{H}_{+i}\right)$
is said to be \textbf{tame }if $\lim_{y\downarrow0}y\psi\left(x+iy\right)$
exists and is finite for all $x\in\mathbb{Q}$.

\vphantom{}\emph{Remark}s:

i. Every rational function of $e^{2\pi iz}$ whose poles are all simple
(i.e. of degree $1$) is tame.

ii. The sum of $\psi+\varphi$ a tame function $\psi$ with a function
$\varphi\in\mathcal{A}\left(\mathbb{H}_{+i}\right)$ admitting a continuous
extension to the closure of $\mathbb{H}_{+i}$ is always tame.

iii. If $\psi$ is tame, then $\textrm{P}\left(\psi\right)$ contains
$\mathbb{Q}$. Since $\mathbb{Q}$ is always $H$-branch invariant,
and since the complement of $\mathbb{Q}$ in $\mathbb{Q}$---the
empty set---is also $H$-branch invariant, observe that $\psi$ being
tame guarantees that $\psi$ admits a dreamcatcher over $\mathbb{Q}$.

iv. If $\psi$ is a tame fixed point of $\mathscr{Q}_{H}$, the \textbf{Singularity
Conservation Law }immediately applies, and tells us that $\psi$'s
dreamcatcher is a fixed point of $\mathfrak{Q}_{H}$ in $L^{2}\left(\mathbb{Q}/\mathbb{Z}\right)$.

\pagebreak{}

\section{The Finite Dreamcatcher Theorem}

Having formulated the notion of a \textbf{Dreamcatcher}, we have shown
that if a function $\psi\in\mathcal{A}\left(\mathbb{H}_{+i}\right)$
admits a dreamcatcher $R$ over some set $T\subseteq\mathbb{R}/\mathbb{Z}$,
then that dreamcatcher must be fixed by the operator $\mathfrak{Q}_{H:T}$.
That is to say, $R$ of the functional equation:
\begin{equation}
R\left(x\right)=\frac{1}{\varrho}\sum_{j=0}^{\varrho-1}e^{-2\pi i\frac{b_{j}}{a_{j}}x}\sum_{k=0}^{\mu_{j}-1}\xi_{\mu_{j}}^{-kH\left(j\right)}R\left(\frac{\varrho x+k}{\mu_{j}}\right)\label{eq:H-SCE}
\end{equation}
over $L^{2}\left(T/\mathbb{Z}\right)$. The author calls this functional
equation the \textbf{Singularity Conservation Equation induced by
$H$} (\textbf{$H$-SCE}).

The technical issues discussed in Chapter 2 evidence a tension latent
in our methods. On one side, we have the original objects of study---set-series
and permutation operators---on the other, we have a partial reformulation
of the question of investigating $\textrm{Ker}\left(1-\mathscr{Q}_{H}\right)$
as a functional equation over an $L^{2}$ space. The tension between
these two sides comes from the obstacle posed by the technical requirements.
In order to use this approach to draw useful conclusions about a particular
$H$-invariant set $V$, we need $\textrm{P}\left(\psi_{V}\right)$
to contain an $H$-branch invariant subset $T$; moreover, to make
the techniques of the present chapter applicable, we will also need
$T$ to be a subgroup of $\mathbb{Q}/\mathbb{Z}$. This carves a ravinehasm
into the ground, separating the site of our methods (Hilbert space)
from the problem to which we hope to apply them (the dynamics of hydra
maps). And although there is a bridge from one side to the other we
cannot currently state with confidence that it is suitable to bear
the weight of full generality.

As such, rather than continue troubling ourselves about this divide,
let us set aside these concerns and focus on the Dreamcatcher operator
$\mathfrak{Q}_{H}$ and its Singularity Conservation Equation purely
as a problem formulated over $L^{2}\left(\mathbb{Q}/\mathbb{Z}\right)$,
given to us \emph{in medias res},\emph{ }without any reference to
the set-series from whence it was derived. The problem we have been
tasked to solve is straight-forward: \emph{characterize the fixed
points} \emph{of $\mathfrak{Q}_{H}$ in $L^{2}\left(\mathbb{Q}/\mathbb{Z}\right)$}.
In this chapter, we completely characterize the fixed points of $\mathfrak{Q}_{H}$
whose supports are finite---that is, those $R\in\textrm{Ker}\left(1-\mathfrak{Q}_{H}\right)$
for which $R\left(x\right)\neq0$ for at most finitely many distinct
values of $x\in\mathbb{Q}$ modulo $1$. This result is the\textbf{
Finite Dreamcatcher Theorem }(\textbf{FDT}).

\subsection{Fourier Analysis on $\mathbb{Q}/\mathbb{Z}$, $\mathbb{Z}_{p}$,
and Related Number Spaces}

\subsubsection{Pontryagin Duals, Characters, \& Duality Brackets}

This section mostly consists of the formulae, definitions, and notational
conventions to be used throughout the proof of the \textbf{FDT}.

\vphantom{}\textbf{Definition} \textbf{19}:

I. Let $t\in\mathbb{Q}$ and $a\in\mathbb{N}_{1}$. We will say $t$
\textbf{lies on }/\textbf{ is on }$a$ if the denominator of $t$
is a non-zero integer multiple of $a$. Invesely, we say $t$ \textbf{lies
off }/ \textbf{is off} $a$ whenever the denominator of $t$ is co-prime
to $a$.

Note that $t$ is on (resp. off) $a$ if and only if $\left|t\right|_{a}>1$
(resp. $\left|t\right|_{a}\leq1$), where $\left|t\right|_{a}$ is
the $a$-adic absolute value of $t$.

II. Given an interval $I\subseteq\mathbb{R}$, we write $I_{a}$ (resp.,
$I_{\bcancel{a}}$) (read {[}I on a{]} (resp. {[}I off a{]})) to denote
the set of all rational numbers in $I$ that are on (resp., off) $a$.
Additionally, we will write $\mathbb{Q}_{\bcancel{a}}$ (read {[}Q
off a{]}) to denote the set of all rational numbers which are off
$a$:
\[
\mathbb{Q}_{\bcancel{a}}\overset{\textrm{def}}{=}\left\{ x\in\mathbb{Q}:\left|x\right|_{a}\leq1\right\} 
\]

\vphantom{}As discussed above, the fixed points $R\left(x\right)$
of $\mathfrak{Q}_{H}$ are elements of $L^{2}\left(\mathbb{Q}/\mathbb{Z}\right)$.
As such, for the reader's benefit, let us go through a brisk review
of details involved in constructing the Pontryagin dual of $\mathbb{Q}/\mathbb{Z}$,
and realizing the duality bracket and Fourier transforms therein.
Performing harmonic analysis on the profinite integers requires some
rather intense notational prestidigitation in order to make it even
remotely presentable-looking in a paper such as this, and so, the
author has strived to develop notation that balances clarity with
utility. As such, we ask that the reader either commit these notations
to memory or start taking notes so as to avoid having to repeatedly
refer back to them later on.

\vphantom{}\textbf{Fact 1}:\textbf{ }For any integer $\nu\in\mathbb{N}_{2}$
(or, more generally, for any non-zero, non-unit algebraic integer
$\nu$), the locally compact abelian groups $\mathbb{Z}\left[\frac{1}{\nu}\right]$
and $\mathbb{Z}_{\nu}$ and are Pontryagin dual to one another. For
each $t\in\mathbb{Z}\left[\frac{1}{\nu}\right]/\mathbb{Z}$ and $\mathfrak{y}\in\mathbb{Z}_{\nu}$,
the duality bracket is given by:
\[
\left(t,\mathfrak{y}\right)\in\left(\mathbb{Z}\left[\frac{1}{\nu}\right]/\mathbb{Z}\right)\times\mathbb{Z}_{\beta}\mapsto\left\{ t\mathfrak{y}\right\} _{\nu}\in\mathbb{Q}
\]
where $\left\{ \cdot\right\} _{\nu}:\mathbb{Q}_{\nu}\rightarrow\mathbb{Q}$
is the $\nu$-adic fractional part, and where the additive group $\mathbb{Z}\left[\frac{1}{\nu}\right]/\mathbb{Z}$
is identified with its unique representatives in $\left[0,1\right)\cap\mathbb{Z}\left[\frac{1}{\nu}\right]$,
equipped with the law of addition modulo $1$. As such, every character
$\chi:\mathbb{Z}\left[\frac{1}{\nu}\right]/\mathbb{Z}\rightarrow\partial\mathbb{D}$
is of the form:
\[
\chi\left(t\right)=e^{2\pi i\left\{ t\mathfrak{y}\right\} _{p}}
\]
for some unique $\mathfrak{y}\in\mathbb{Z}_{\nu}$, and every character
$\check{\chi}:\mathbb{Z}_{\nu}\rightarrow\partial\mathbb{D}$ is of
the form:
\[
\check{\chi}\left(\mathfrak{y}\right)=e^{2\pi i\left\{ t\mathfrak{y}\right\} _{p}}
\]
for some unique $t\in\mathbb{Z}\left[\frac{1}{\nu}\right]/\mathbb{Z}$,
where---again---$\mathbb{Z}\left[\frac{1}{\nu}\right]/\mathbb{Z}$
is identified with $\left[0,1\right)\cap\mathbb{Z}\left[\frac{1}{\nu}\right]$.

\vphantom{}\textbf{Fact 2}: By the \textbf{Structure Theorem for
Abelian Groups}, the group $\left(\mathbb{Q}/\mathbb{Z},+\right)$
satisfies the group isomorphism:
\begin{equation}
\mathbb{Q}/\mathbb{Z}\cong\bigoplus_{p\in\mathbb{P}}\left(\mathbb{Z}\left[\frac{1}{p}\right]/\mathbb{Z}\right)\label{eq:p-power torsion decomposition of Q/Z}
\end{equation}

\vphantom{}\textbf{Fact 3}: 
\[
\widehat{\mathbb{Q}/\mathbb{Z}}\cong\prod_{p\in\mathbb{P}}\mathbb{Z}_{p}\cong\widetilde{\mathbb{Z}}
\]
is a group isomorphism, where---recall---$\widetilde{\mathbb{Z}}$
denotes the ring of \textbf{profinite integers}. In our viewopint,
a profinite integer $\mathfrak{z}\in\widetilde{\mathbb{Z}}$ is an
infinite sequence of $p$-adic integers, $\mathfrak{z}_{p}$, indexed
by the prime numbers $p$:
\[
\mathfrak{z}=\left\{ \mathfrak{z}_{p}\right\} _{p\in\mathbb{P}}=\left(\mathfrak{z}_{2},\mathfrak{z}_{3},\mathfrak{z}_{5},\ldots\right)
\]
The duality bracket between $\mathbb{Q}/\mathbb{Z}$ and $\widetilde{\mathbb{Z}}$
is then defined by:
\begin{equation}
\left\langle t,\mathfrak{z}\right\rangle \overset{\textrm{def}}{=}\sum_{p\in\mathbb{P}}\left\{ t\mathfrak{z}_{p}\right\} _{p},\textrm{ }\forall t\in\mathbb{Q}/\mathbb{Z},\forall\mathfrak{z}\in\widetilde{\mathbb{Z}}\label{eq:Def of Profinite integer duality bracket}
\end{equation}
where, as usual, $\mathbb{Q}/\mathbb{Z}$ is identified with $\left[0,1\right)_{\mathbb{Q}}$,
and where $\mathfrak{z}_{p}$ denotes the $p$-adic component of $\mathfrak{z}$.
Since every such $t$ has only finitely many prime factors in its
denominator, all but finitely many terms of the sum are equal to $1$
for any given $t$, and thus, the sum defines a function from $\mathbb{Q}/\mathbb{Z}$
to $\mathbb{Q}$.

It can be shown that for every character $\chi$ of $\mathbb{Q}/\mathbb{Z}$,
there is a unique $\mathfrak{z}\in\widetilde{\mathbb{Z}}$ so that
$\chi$ can has the formula:
\[
t\mapsto\chi\left(t\right)=e^{2\pi i\left\langle t,\mathfrak{z}\right\rangle }\overset{\textrm{def}}{=}\prod_{p\in\mathbb{P}}e^{2\pi i\left\{ t\mathfrak{z}_{p}\right\} _{p}},\textrm{ }\forall t\in\mathbb{Q}/\mathbb{Z}
\]
Likewise, for every character $\check{\chi}$ of $\widetilde{\mathbb{Z}}$,
there is a unique $t\in\mathbb{Q}/\mathbb{Z}=\left[0,1\right)_{\mathbb{Q}}$
so that $\check{\chi}$is given by:
\[
\mathfrak{z}\mapsto\check{\chi}\left(\mathfrak{z}\right)=e^{2\pi i\left\langle t,\mathfrak{z}\right\rangle }\overset{\textrm{def}}{=}\prod_{p\in\mathbb{P}}e^{2\pi i\left\{ t\mathfrak{z}_{p}\right\} _{p}},\textrm{ }\forall\mathfrak{z}\in\widetilde{\mathbb{Z}}
\]

\vphantom{}\textbf{Fact 4}: Let $\varrho\in\mathbb{N}_{2}$. Then,
we have the group isomorphism:
\[
\mathbb{Q}_{\bcancel{\varrho}}/\mathbb{Z}=\left(\mathbb{Q}/\mathbb{Z}\right)_{\bcancel{\varrho}}\cong\bigoplus_{p\in\mathbb{P}:p\nmid\varrho}\left(\mathbb{Z}\left[\frac{1}{p}\right]/\mathbb{Z}\right)=\bigoplus_{p\in\mathbb{P}:p\nmid\varrho}\hat{\mathbb{Z}}_{p}
\]
where $\hat{\mathbb{Z}}_{p}\cong\mathbb{Z}\left[\frac{1}{p}\right]/\mathbb{Z}$
is the Pontryagin dual of the additive group of $p$-adic integers,
also known as the \textbf{Prüfer} $p$\textbf{-group}---although
that name is generally reserved for the group-isomorphic copy of $\mathbb{Z}\left[\frac{1}{p}\right]/\mathbb{Z}$
in $\partial\mathbb{D}$ defined by $\left(\left\{ e^{2\pi ix}:x\in\mathbb{Z}\left[\frac{1}{p}\right]/\mathbb{Z}\right\} ,\times\right)$.
Since taking the pontryagin dual of a direct sum yields the direct
product of the pontryagin duals (and vice-versa), the corresponding
group isomorphism of the dual is
\[
\widehat{\mathbb{Q}_{\bcancel{\varrho}}/\mathbb{Z}}=\widehat{\left(\mathbb{Q}/\mathbb{Z}\right)_{\bcancel{\varrho}}}\cong\prod_{\begin{array}{c}
p\in\mathbb{P}\\
p\nmid\varrho
\end{array}}\mathbb{Z}_{p}
\]

\vphantom{}\textbf{Definition} \textbf{20}: Let $\varrho\in\mathbb{N}_{2}$.
Then, the ring of \textbf{profinite integers off $\varrho$}, denoted
$\widetilde{\mathbb{Z}}_{\bcancel{\varrho}}$, is defined by: 
\[
\widetilde{\mathbb{Z}}_{\bcancel{\varrho}}\overset{\textrm{def}}{=}\prod_{\begin{array}{c}
p\in\mathbb{P}\\
p\nmid\varrho
\end{array}}\mathbb{Z}_{p}
\]
Note that $\widetilde{\mathbb{Z}}_{\bcancel{\varrho}}$ is a subring
of $\widetilde{\mathbb{Z}}$.

Now, let $N$ be the number of prime divisors of $\varrho$, and enumerate
the elements of $\mathbb{P}$ as $p_{1},p_{2},\ldots$, where $p_{1},\ldots,p_{N}$
are the prime divisors of $\varrho$. Then, every element of $\widetilde{\mathbb{Z}}_{\bcancel{\varrho}}$
can be written as $\tilde{\mathfrak{z}}=\left(\mathfrak{z}_{p_{N+1}},\mathfrak{z}_{p_{N+2}},\ldots\right)$,
where $\mathfrak{z}_{p_{n}}\in\mathbb{Z}_{p_{n}}$ for every $n\in\mathbb{N}_{N+1}$.
Then, the map:
\begin{equation}
\tilde{\mathfrak{z}}=\left(\mathfrak{z}_{p_{N+1}},\mathfrak{z}_{p_{N+2}},\ldots\right)\mapsto\mathfrak{z}=\left(0,0,\ldots,0,\mathfrak{z}_{p_{N+1}},\mathfrak{z}_{p_{N+2}},\ldots\right)\label{eq:Embedding formula}
\end{equation}
which concatenates an $N$-tuple of $0$s to the left of $\tilde{\mathfrak{z}}$
then defines an embedding of $\widetilde{\mathbb{Z}}_{\bcancel{\varrho}}$
in $\widetilde{\mathbb{Z}}$.  

\vphantom{}\textbf{Fact 5}: For all $x\in\mathbb{Q}/\mathbb{Z}$
and all $\mathfrak{z}\in\widetilde{\mathbb{Z}}$:

I. $\left\langle x,a\mathfrak{z}\right\rangle =\left\langle ax,\mathfrak{z}\right\rangle $,
$\forall a\in\mathbb{Q}$.

II. Let $p\in\mathbb{P}$ and $n\in\mathbb{N}_{1}$. Then, $\left\langle \frac{1}{p^{n}},\mathfrak{z}\right\rangle =\frac{1}{p^{n}}\left[\mathfrak{z}_{p}\right]_{p^{n}}$,
where, recall, $\left[\cdot\right]_{p^{n}}:\mathbb{Z}_{p}\rightarrow\left\{ 0,\ldots,p^{n}-1\right\} $
is ``mod $p^{n}$'' projection operation.

IV. $\left\langle a,b\right\rangle =ab$ for all $a\in\mathbb{Q}$
and all $b\in\mathbb{Z}\subseteq\widetilde{\mathbb{Z}}$.

V. On occasion, we will write $\left\langle x,\mathfrak{z}\right\rangle _{\bcancel{\varrho}}$
to denote the duality bracket on $\widetilde{\mathbb{Z}}_{\bcancel{\varrho}}$.
By the above embedding of $\widetilde{\mathbb{Z}}_{\bcancel{\varrho}}$
in $\widetilde{\mathbb{Z}}$, it then follows that:
\[
\left\langle x,\mathfrak{z}\right\rangle _{\bcancel{\varrho}}\overset{\textrm{def}}{=}\sum_{\begin{array}{c}
p\in\mathbb{P}\\
p\nmid\varrho
\end{array}}\left\{ x\mathfrak{z}\right\} _{p}
\]
Consequently, whenever the domain of $\mathfrak{z}$ is made clear
(such as when working with an integral) we will write $\left\langle x,\mathfrak{z}\right\rangle $
to denote the duality bracket over the indicated subring of $\widetilde{\mathbb{Z}}$.

VI.
\[
\sum_{p\in\mathbb{P}}\left\{ x\right\} _{p}\overset{1}{\equiv}x,\textrm{ }\forall x\in\mathbb{Q}
\]
(Note that this identity proves (IV).)

\subsubsection{The Fourier Transform over $\widetilde{\mathbb{Z}}$;\emph{ }Integration
over the Profinite and $p$-Adic Integers}

This subsection is primarily a repository of definitions, theorems,
and useful formulae. It should be noted that, for all of the following
results, every instance $\mathbb{Q}/\mathbb{Z}$ (resp., $\widetilde{\mathbb{Z}}$)
can be replaced with $G$ (resp. $\check{G}$), where $G$ is any
subgroup of $\mathbb{Q}/\mathbb{Z}$ (resp. $\check{G}$ is the pontryagin
dual of $G$), and letting $d\mathfrak{z}$ be the Haar probability
measure on $\check{G}$.

\vphantom{}\textbf{Definition} \textbf{21}:

I. Given $t\in\mathbb{R}$, the function $\mathbf{1}_{t}:\mathbb{R}\rightarrow\left\{ 0,1\right\} $
is defined by:
\[
\mathbf{1}_{t}\left(x\right)\overset{\textrm{def}}{=}\begin{cases}
1 & \textrm{if }x\overset{1}{\equiv}t\\
0 & \textrm{else}
\end{cases}=\left[x\overset{t}{\equiv}1\right]
\]
That is, $\mathbf{1}_{t}$ is the \textbf{indicator function} for
the set $t+\mathbb{Z}$. Note that $\left\{ \mathbf{1}_{t}\right\} _{t\in\left[0,1\right)_{\mathbb{Q}}}$
is then an orthonormal basis for $L^{2}\left(\mathbb{Q}/\mathbb{Z}\right)$.

II. A \textbf{Fourier series }on $\widetilde{\mathbb{Z}}$ is a series
of the form:
\[
\sum_{t\in\mathbb{Q}/\mathbb{Z}}c_{t}e^{2\pi i\left\langle t,\mathfrak{z}\right\rangle }
\]
where the $c_{t}$s are complex numbers, called the\textbf{ Fourier
coefficients}, and where $\mathfrak{z}$ is a variable taking values
in $\widetilde{\mathbb{Z}}$.

\vphantom{}\emph{Remark}: One might say that the more ``natural''
perspective would be to work with the space of complex-valued functions
on $\widetilde{\mathbb{Z}}$ which are representable by Fourier series.
However, for our purposes, we will spend most of our time working
with the space of Fourier coefficients---complex valued functions
on $\mathbb{Q}/\mathbb{Z}$.

III. The \textbf{Fourier Transform $\mathscr{F}_{\mathbb{Q}/\mathbb{Z}}:L^{1}\left(\mathbb{Q}/\mathbb{Z}\right)\rightarrow L^{\infty}\left(\widetilde{\mathbb{Z}}\right)$
}is defined by:
\[
\mathscr{F}_{\mathbb{Q}/\mathbb{Z}}\left\{ f\right\} \left(\mathfrak{z}\right)\overset{\textrm{def}}{=}\sum_{t\in\left[0,1\right)_{\mathbb{Q}}}f\left(t\right)e^{2\pi i\left\langle t,\mathfrak{z}\right\rangle },\textrm{ }\forall\mathfrak{z}\in\widetilde{\mathbb{Z}}
\]
for all $f\left(x\right)\in L^{1}\left(\mathbb{Q}/\mathbb{Z}\right)$.
We will generally write $\check{f}$ (pronounced {[}$f$-check{]})
to denote the Fourier transform of $f$: 
\[
\check{f}\left(\mathfrak{z}\right)\overset{\textrm{def}}{=}\mathscr{F}_{\mathbb{Q}/\mathbb{Z}}\left\{ f\right\} \left(\mathfrak{z}\right)
\]

\vphantom{}\emph{Remark}: Thus, $\mathscr{F}_{\mathbb{Q}/\mathbb{Z}}$
sends Fourier coefficients (functions on $\mathbb{Q}/\mathbb{Z}$,
a discrete abelian group) to their Fourier series (functions on $\widetilde{\mathbb{Z}}$,
a compact abelian group); this is behind the author's choice to write
$\check{f}$ rather than the more common $\hat{f}$.

IV. The \textbf{Inverse Fourier Transform }$\mathscr{F}_{\mathbb{Q}/\mathbb{Z}}^{-1}:L^{1}\left(\widetilde{\mathbb{Z}}\right)\rightarrow L^{\infty}\left(\mathbb{Q}/\mathbb{Z}\right)$
is defined by:
\[
\mathscr{F}_{\mathbb{Q}/\mathbb{Z}}^{-1}\left\{ F\right\} \left(x\right)\overset{\textrm{def}}{=}\int_{\widetilde{\mathbb{Z}}}F\left(x\right)e^{-2\pi i\left\langle x,\mathfrak{z}\right\rangle }d\mathfrak{z},\textrm{ }\forall x\in\mathbb{Q}/\mathbb{Z}
\]
Here, $d\mathfrak{z}=\prod_{p\in\mathbb{P}}d\mathfrak{z}_{p}$ is
the product Haar measure on $\widetilde{\mathbb{Z}}$ with unit normalization
($\int_{\widetilde{\mathbb{Z}}}d\mathfrak{z}\overset{\textrm{def}}{=}1$),
where $d\mathfrak{z}_{p}$ is the Haar measure on $\mathbb{Z}_{p}$
with unit normalization ($\int_{\mathbb{Z}_{p}}d\mathfrak{z}_{p}\overset{\textrm{def}}{=}1$).

Given an $F:\widetilde{\mathbb{Z}}\rightarrow\mathbb{C}$, we use
the short-hand $\hat{F}$ (pronounced {[}$F$-hat{]}) to denote the
inverse Fourier transform of $F$:
\[
\hat{F}\left(x\right)\overset{\textrm{def}}{=}\mathscr{F}_{\mathbb{Q}/\mathbb{Z}}^{-1}\left\{ F\right\} \left(x\right)
\]

\vphantom{}Since $\mathscr{F}_{\mathbb{Q}/\mathbb{Z}}$ is defined
on $L^{1}\left(\mathbb{Q}/\mathbb{Z}\right)$, and satisfies the \textbf{Parseval-Plancherel
Identity}:
\[
\sum_{t\in\mathbb{Q}/\mathbb{Z}}f\left(t\right)\overline{g\left(t\right)}=\int_{\widetilde{\mathbb{Z}}}\check{f}\left(\mathfrak{z}\right)\overline{\check{g}\left(\mathfrak{z}\right)}d\mathfrak{z},\textrm{ }\forall f,g\in L^{2}\left(\mathbb{Q}/\mathbb{Z}\right)
\]
an application of the standard ``extend by continuity'' argument
shows that $\mathscr{F}_{\mathbb{Q}/\mathbb{Z}}$ can be extended
to an operator on $L^{1}\left(\mathbb{Q}/\mathbb{Z}\right)\cap L^{2}\left(\mathbb{Q}/\mathbb{Z}\right)$,
and from there to an operator on $L^{r}\left(\mathbb{Q}/\mathbb{Z}\right)$,
for any $r\in\mathbb{R}\cap\left[1,2\right]$. From this, one obtains
the usual theory of the Fourier Transform:

\vphantom{}\textbf{Theorem 5: The Fourier Inversion Theorem}: Let
$r\in\mathbb{R}\cap\left[1,2\right]$. Then:

I. For every $f\left(x\right)\in L^{r}\left(\mathbb{Q}/\mathbb{Z}\right)$,
$\mathscr{F}_{\mathbb{Q}/\mathbb{Z}}^{-1}\left\{ \check{f}\right\} \left(x\right)$
is equal to $f\left(x\right)$ for every\footnote{In general, these statements only hold almost everywhere; i.e., ``for
almost every $x\in\mathbb{Q}/\mathbb{Z}$''. However, since the measure
on $\mathbb{Q}/\mathbb{Z}$ is the counting measure, the only set
of measure zero in $\mathbb{Q}/\mathbb{Z}$ is the empty set, and
thus, ``equality for almost every $x\in\mathbb{Q}/\mathbb{Z}$''
is synonymous with ``equality \emph{for every }$x\in\mathbb{Q}/\mathbb{Z}$''.} $x\in\mathbb{Q}/\mathbb{Z}$; that is to say:
\[
f\left(x\right)\in L^{r}\left(\mathbb{Q}/\mathbb{Z}\right)\Rightarrow f\left(x\right)=\int_{\widetilde{\mathbb{Z}}}\check{f}\left(\mathfrak{z}\right)e^{-2\pi i\left\langle x,\mathfrak{z}\right\rangle },\textrm{ }\forall x\in\mathbb{Q}/\mathbb{Z}
\]

II. For every $F\left(\mathfrak{z}\right)\in L^{r}\left(\widetilde{\mathbb{Z}}\right)$,
$\mathscr{F}_{\mathbb{Q}/\mathbb{Z}}\left\{ \hat{F}\right\} \left(\mathfrak{z}\right)$
is equal to $F\left(\mathfrak{z}\right)$ for almost every $\mathfrak{z}\in\widetilde{\mathbb{Z}}$;
that is to say:
\[
F\left(\mathfrak{z}\right)\in L^{r}\left(\widetilde{\mathbb{Z}}\right)\Rightarrow F\left(\mathfrak{z}\right)=\sum_{t\in\left[0,1\right)_{\mathbb{Q}}}\hat{F}\left(t\right)e^{2\pi i\left\langle t,\mathfrak{z}\right\rangle },\textrm{ }\tilde{\forall}\mathfrak{z}\in\widetilde{\mathbb{Z}}
\]

\vphantom{}\textbf{Theorem 6: Parseval's Theorem}: $\mathscr{F}_{\mathbb{Q}/\mathbb{Z}}:L^{2}\left(\mathbb{Q}/\mathbb{Z}\right)\rightarrow L^{2}\left(\widetilde{\mathbb{Z}}\right)$
is an isometric isomorphism of Hilbert spaces. Moreover, for every
$f\left(x\right)\in L^{2}\left(\mathbb{Q}/\mathbb{Z}\right)$ (resp.,
$F\left(\mathfrak{z}\right)\in L^{2}\left(\widetilde{\mathbb{Z}}\right)$),
$\check{f}\left(\mathfrak{z}\right)$ (resp, $\hat{F}\left(x\right)$)
is defined for almost every $\mathfrak{z}\in\widetilde{\mathbb{Z}}$
(resp. for every $x\in\mathbb{Q}/\mathbb{Z}$).

\vphantom{}\textbf{Notations and Basic Formulae}:

I. Given $t\in\left[0,1\right)_{\mathbb{Q}}$:
\[
\mathscr{F}_{\mathbb{Q}/\mathbb{Z}}\left\{ \mathbf{1}_{t}\right\} \left(\mathfrak{z}\right)=e^{2\pi i\left\langle t,\mathfrak{z}\right\rangle },\textrm{ }\forall\mathfrak{z}\in\widetilde{\mathbb{Z}}
\]
As such, $\check{\mathbf{1}}_{t}\left(\mathfrak{z}\right)$ (pronounced
{[}one-$t$-check of zed{]}) will frequently be used to denote $e^{2\pi i\left\langle t,\mathfrak{z}\right\rangle }$:
\[
\check{\mathbf{1}}_{t}\left(\mathfrak{z}\right)\overset{\textrm{def}}{=}e^{2\pi i\left\langle t,\mathfrak{z}\right\rangle }
\]
Moreover, by \textbf{Parseval's Theorem}, it then follows that $\left\{ \check{\mathbf{1}}_{t}\right\} _{t\in\left[0,1\right)_{\mathbb{Q}}}$is
an orthonormal basis for $L^{2}\left(\widetilde{\mathbb{Z}}\right)$.

II. (\textbf{$p$-adic integration techniques}) Let $p\in\mathbb{P}$,
let $f:\mathbb{Z}_{p}\rightarrow\mathbb{C}$, and let $a,b\in\mathbb{Q}$,
with $a\neq0$, let $\mathfrak{c}\in\mathbb{Z}_{p}$, and let $n\in\mathbb{N}_{0}$.
Then:
\[
\int_{\mathbb{Z}_{p}}f\left(a\mathfrak{y}+b\right)d\mathfrak{y}=\frac{1}{\left|a\right|_{p}}\int_{a\mathbb{Z}_{p}+b}f\left(\mathfrak{x}\right)d\mathfrak{x}
\]
\[
\int_{a\mathbb{Z}_{p}+b}f\left(\mathfrak{x}\right)d\mathfrak{x}=\int_{a\mathbb{Z}_{p}}f\left(\mathfrak{x}+b\right)d\mathfrak{x}=\left|a\right|_{p}\int_{\mathbb{Z}_{p}}f\left(a\mathfrak{y}+b\right)d\mathfrak{y}
\]

\[
\int_{\mathbb{Z}_{p}}f\left(\mathfrak{y}+\mathfrak{c}\right)d\mathfrak{y}=\int_{\mathbb{Z}_{p}+\mathfrak{c}}f\left(\mathfrak{y}\right)d\mathfrak{y}=\int_{\mathbb{Z}_{p}}f\left(\mathfrak{y}\right)d\mathfrak{y}
\]
\[
\int_{\mathbb{Z}_{p}}f\left(\mathfrak{y}\right)d\mathfrak{y}=\sum_{k=0}^{p^{n}-1}\int_{k+p^{n}\mathbb{Z}_{p}}f\left(\mathfrak{y}\right)d\mathfrak{y}
\]

III. Let $p\in\mathbb{P}$. Then:
\[
\int_{\mathbb{Z}_{p}}e^{2\pi i\left\{ x\mathfrak{y}\right\} _{p}}d\mathfrak{y}=\int_{\mathbb{Z}_{p}}e^{-2\pi i\left\{ x\mathfrak{y}\right\} _{p}}d\mathfrak{y}=\mathbf{1}_{\mathbb{Z}_{p}}\left(x\right),\textrm{ }\forall x\in\mathbb{Q}
\]

IV. Let $Q\subseteq\mathbb{P}$ be set containing all but finitely
many primes (i.e., $\left|\mathbb{P}\backslash Q\right|<\infty$).
Then:
\[
\prod_{q\in Q}\mathbf{1}_{\mathbb{Z}_{q}}\left(x\right)=\mathbf{1}_{\mathbb{Z}\left[\prod_{p\in\mathbb{P}\backslash Q}\frac{1}{p}\right]}\left(x\right)
\]

V.

\[
\int_{\widetilde{\mathbb{Z}}}e^{2\pi i\left\langle x,\mathfrak{z}\right\rangle }d\mathfrak{z}=\int_{\widetilde{\mathbb{Z}}}e^{-2\pi i\left\langle x,\mathfrak{z}\right\rangle }d\mathfrak{z}=\mathbf{1}_{\mathbb{Z}}\left(x\right),\textrm{ }\forall x\in\mathbb{Q}
\]

VI. Let $A$ be a measurable set in $\mathbb{Z}_{p}$ of measure $\textrm{meas}_{p}\left(A\right)$.
Then:
\[
\int_{\mathbb{Z}_{p}}\left[\mathfrak{y}\in A\right]d\mathfrak{y}=\textrm{meas}_{p}\left(A\right)
\]
Additionally, let $B$ be a measurable set in $\widetilde{\mathbb{Z}}$
of measure $\textrm{meas}\left(B\right)$. Then:
\[
\int_{\widetilde{\mathbb{Z}}}\left[\mathfrak{z}\in B\right]d\mathfrak{z}=\textrm{meas}\left(B\right)
\]

VII. Recalling the ``prime factorization'' of $\mathbb{Q}/\mathbb{Z}$
into a direct sum of the Prüfer $p$-groups $\hat{\mathbb{Z}}_{p}$
of prime order $p$:
\[
\mathbb{Q}/\mathbb{Z}\cong\bigoplus_{p\in\mathbb{P}}\hat{\mathbb{Z}}_{p}
\]
note that, given any integer $a\in\mathbb{N}_{2}$, we can collect
the terms of this direct sum into two groups: $\left(\mathbb{Q}/\mathbb{Z}\right)_{a}$,
the part of $\mathbb{Q}/\mathbb{Z}$ that lies ``on $a$'':

\[
\hat{\mathbb{Z}}_{a}\cong\left(\mathbb{Q}/\mathbb{Z}\right)_{a}=\bigoplus_{\begin{array}{c}
p\in\mathbb{P}\\
p\mid a
\end{array}}\hat{\mathbb{Z}}_{p}
\]
and $\left(\mathbb{Q}/\mathbb{Z}\right)_{\bcancel{a}}$, the part
of $\mathbb{Q}/\mathbb{Z}$ that lies ``off $a$'':
\[
\left(\mathbb{Q}/\mathbb{Z}\right)_{\bcancel{a}}=\bigoplus_{\begin{array}{c}
p\in\mathbb{P}\\
p\nmid a
\end{array}}\hat{\mathbb{Z}}_{p}
\]
Consequently, we can write:
\[
\mathbb{Q}/\mathbb{Z}\cong\hat{\mathbb{Z}}_{a}\oplus\left(\mathbb{Q}/\mathbb{Z}\right)_{\bcancel{a}}\cong\hat{\mathbb{Z}}_{a}\times\left(\mathbb{Q}/\mathbb{Z}\right)_{\bcancel{a}}
\]
where we used the fact the the direct sum of finitely many things
is the same as the direct product of said finitely many things. Taking
Pontryagin duals gives:
\[
\widetilde{\mathbb{Z}}\cong\mathbb{Z}_{a}\times\widetilde{\mathbb{Z}}_{\bcancel{a}}
\]
We will \emph{constantly }use this type of decomposition to evaluate
integrals over $\widetilde{\mathbb{Z}}$. Indeed, for our purposes,
the integrals:
\[
\int_{\widetilde{\mathbb{Z}}}f\left(\mathfrak{z}\right)d\mathfrak{z}
\]
we shall consider will be such that there is an integer $a$ so that,
upon identifying $\mathfrak{z}\in\widetilde{\mathbb{Z}}$ with $\left(\mathfrak{z}_{a},\mathfrak{z}_{\bcancel{a}}\right)\in\mathbb{Z}_{a}\times\widetilde{\mathbb{Z}}_{\bcancel{a}}$
and re-writing the integral as:
\[
\int_{\widetilde{\mathbb{Z}}}f\left(\mathfrak{z}\right)d\mathfrak{z}=\int_{\mathbb{Z}_{a}\times\widetilde{\mathbb{Z}}_{\bcancel{a}}}f\left(\mathfrak{z}\right)d\mathfrak{z}=\int_{\widetilde{\mathbb{Z}}_{\bcancel{a}}}\left(\int_{\mathbb{Z}_{a}}f\left(\mathfrak{z}_{a},\mathfrak{z}_{\bcancel{a}}\right)d\mathfrak{z}_{a}\right)d\mathfrak{z}_{\bcancel{a}}
\]
the inner integral over $\mathbb{Z}_{a}$ will be amenable to a change-of-variables
(as in (II)). Doing so will allow us to simplify the integrand at
the cost of altering the domain of integration. Combining the integral
back into a single integral over $\widetilde{\mathbb{Z}}$ then entails
re-writing the altered domain of integration in terms of an indicator
function/iverson bracket. The general rule is that a change of variables
of the form $\mathfrak{y}=a\mathfrak{z}+b$ (where $a$ and $b$ are
integers) will have no effect on the $p$-adic components of $\widetilde{\mathbb{Z}}$
for the primes $p$ which are co-prime to $a$, in which case, $\mathfrak{z}_{p}\in\mathbb{Z}_{p}\mapsto a\mathfrak{z}_{p}+b\in\mathbb{Z}_{p}$
is then a (continuous) measure-preserving bijection of $\mathbb{Z}_{p}$.

For example, consider:
\[
\int_{\widetilde{\mathbb{Z}}}f\left(3\mathfrak{z}+1\right)d\mathfrak{z}=\int_{\widetilde{\mathbb{Z}}_{\bcancel{3}}}\left(\int_{\mathbb{Z}_{3}}f\left(3\mathfrak{z}_{3}+1,3\mathfrak{z}_{\bcancel{3}}+1\right)d\mathfrak{z}_{3}\right)d\mathfrak{z}_{\bcancel{3}}
\]
Substituting $\mathfrak{y}=3\mathfrak{z}+1$ in the integral on the
left is the same as substituting $\mathfrak{y}_{3}=3\mathfrak{z}_{3}+1$
and $3\mathfrak{y}_{\bcancel{3}}+1=3\mathfrak{z}_{\bcancel{3}}+1$
in the integral on the right. The map $\mathfrak{z}_{\bcancel{3}}\in\widetilde{\mathbb{Z}}_{\bcancel{3}}\mapsto3\mathfrak{z}_{\bcancel{3}}+1\in\widetilde{\mathbb{Z}}_{\bcancel{3}}$
is a measure-preserving bijection, and so, that component of the integral
will be unchanged if we undo it:
\begin{eqnarray*}
\int_{\widetilde{\mathbb{Z}}}f\left(3\mathfrak{z}+1\right)d\mathfrak{z} & = & \int_{3\widetilde{\mathbb{Z}}_{\bcancel{3}}+1}\frac{1}{\left|3\right|_{3}}\int_{3\mathbb{Z}_{3}+1}f\left(\mathfrak{y}_{3},\mathfrak{y}_{\bcancel{3}}\right)d\mathfrak{y}_{3}d\mathfrak{y}_{\bcancel{3}}\\
\left(3\widetilde{\mathbb{Z}}_{\bcancel{3}}+1=\widetilde{\mathbb{Z}}_{\bcancel{3}}\right); & = & 3\int_{\widetilde{\mathbb{Z}}_{\bcancel{3}}}\int_{3\mathbb{Z}_{3}+1}f\left(\mathfrak{y}_{3},\mathfrak{y}_{\bcancel{3}}\right)d\mathfrak{y}_{3}d\mathfrak{y}_{\bcancel{3}}\\
\left(\mathfrak{y}_{3}\in3\mathbb{Z}_{3}+1\Leftrightarrow\mathfrak{y}_{3}\overset{3}{\equiv}1\right); & = & 3\int_{\widetilde{\mathbb{Z}}_{\bcancel{3}}}\int_{\mathbb{Z}_{3}}\left[\mathfrak{y}_{3}\overset{3}{\equiv}1\right]f\left(\mathfrak{y}_{3},\mathfrak{y}_{\bcancel{3}}\right)d\mathfrak{y}_{3}d\mathfrak{y}_{\bcancel{3}}\\
 & = & 3\int_{\widetilde{\mathbb{Z}}}\left[\mathfrak{y}_{3}\overset{3}{\equiv}1\right]f\left(\mathfrak{y}\right)d\mathfrak{y}
\end{eqnarray*}

\pagebreak{}

\subsection{The Finite Dreamcatcher Theorem - Heuristics and Preparatory Work}

First, a definition.

\vphantom{}\textbf{Definition 22}: Given $f\in L^{2}\left(\mathbb{Q}/\mathbb{Z}\right)$,
the \textbf{support }of $f$, denoted $\textrm{supp}\left(f\right)$,
is defined by:
\[
\textrm{supp}\left(f\right)\overset{\textrm{def}}{=}\left\{ x\in\mathbb{Q}:f\left(x\right)\neq0\right\} 
\]
Furthermore, $f$ is said to be \textbf{finitely supported }(on $\mathbb{Q}/\mathbb{Z}$)
when $\textrm{supp}\left(f\right)$ contains only finitely many distinct
equivalence classes of rational numbers modulo $1$.

The goal of this chapter is to prove the following theorem:

\vphantom{}\textbf{Theorem 7: The Finite Dreamcatcher Theorem}: Let
$H:\mathbb{N}_{0}\rightarrow\mathbb{N}_{0}$ be a prime\footnote{In all likelihood, the primality hypothesis on $H$ is unecessary.
We invoke it only primarily because it simplifies the $p$-adic computations
(pun intended).}, regulated\footnote{Recall that this means there is a $j\in\left\{ 0,\ldots,\varrho-1\right\} $
so that $\mu_{j}=1$.} $\varrho$-hydra map, and let $R\in\textrm{Ker}\left(1-\mathfrak{Q}_{H}\right)\in L^{2}\left(\mathbb{Q}/\mathbb{Z}\right)$.
If $R$ is finitely supported, then $R$ is of the form:
\[
R\left(x\right)=R\left(0\right)\mathbf{1}_{\mathbb{Z}}\left(x\right),\textrm{ }\forall x\in\mathbb{Q}
\]

\subsubsection{Outline of the Proof of the FDT}

Here is the outline of what remains to be done:

\vphantom{}Step 0: Compute the \textbf{Basis Formula }for $\mathscr{\mathfrak{Q}}_{H}$,
an expression of $\mathfrak{Q}_{H}$ in terms of $\left\{ \mathbf{1}_{t}\right\} _{t\in\left[0,1\right)_{\mathbb{Q}}}$---the
``standard'' orthonormal basis of $L^{2}\left(\mathbb{Q}/\mathbb{Z}\right)$.
Conjugating everything by the Fourier transform $\mathscr{F}_{\mathbb{Q}/\mathbb{Z}}:L^{2}\left(\mathbb{Q}/\mathbb{Z}\right)\rightarrow L^{2}\left(\widetilde{\mathbb{Z}}\right)$,
an analogous formula is derived for the operator $\check{\mathfrak{Q}}_{H}:L^{2}\left(\widetilde{\mathbb{Z}}\right)\rightarrow L^{2}\left(\widetilde{\mathbb{Z}}\right)$
defined by $\check{\mathfrak{Q}}_{H}\overset{\textrm{def}}{=}\mathscr{F}_{\mathbb{Q}/\mathbb{Z}}\circ\mathfrak{Q}_{H}\circ\mathscr{F}_{\mathbb{Q}/\mathbb{Z}}^{-1}$.

As an example to whet the reader's appetite, in the case of the $H_{p}$
maps, where $p$ is an odd prime, the \textbf{Basis Formulae} for
$\mathfrak{Q}_{p}$ and $\check{\mathfrak{Q}}_{p}$ (the author's
notations for $\mathfrak{Q}_{H_{p}}$ and $\check{\mathfrak{Q}}_{H_{p}}$,
respectively) are:
\[
\mathfrak{Q}_{p}\left\{ \mathbf{1}_{t}\right\} \left(x\right)=\frac{1}{2}\left(\mathbf{1}_{\frac{t}{2}}\left(x\right)+\mathbf{1}_{\frac{t+1}{2}}\left(x\right)\right)+\frac{e^{-\pi it}}{2}\left(\mathbf{1}_{\frac{pt}{2}}\left(x\right)-\mathbf{1}_{\frac{pt+1}{2}}\left(x\right)\right),\textrm{ }\forall x\in\mathbb{Q}
\]

\[
\check{\mathfrak{Q}}_{p}\left\{ \check{\mathbf{1}}_{t}\right\} \left(\mathfrak{z}\right)=\left[\mathfrak{z}_{2}\overset{2}{\equiv}0\right]\check{\mathbf{1}}_{t}\left(\frac{\mathfrak{z}}{2}\right)+\left[\mathfrak{z}_{2}\overset{2}{\equiv}1\right]\check{\mathbf{1}}_{t}\left(\frac{p\mathfrak{z}-1}{2}\right),\textrm{ }\forall\mathfrak{z}\in\widetilde{\mathbb{Z}}
\]
where $t\in\mathbb{Q}$ is arbitrary.

\emph{Remark}: As befits a ``Step 0'', the \textbf{Basis Formulae
}will be stated in the \textbf{formula repository }of the next subsection---§3.2.3.

\vphantom{}Step 1: Using the Fourier transform and its inverse, the
assumption of $H$'s reulgarity allows us to prove that that if $R\in\textrm{Ker}\left(1-\mathfrak{Q}_{H}\right)$
has finite support, then $R$ must vanish for all $t\in\mathbb{Q}$
which are on $\varrho$ (that is, for all $t$ with $\left|t\right|_{\varrho}>1$).
Thus, fixed points of $\mathfrak{Q}_{H}$ with finite support are
supported ``off $\varrho$'' (the \textbf{Off-$\varrho$ Theorem}).

\vphantom{}Step 2: Using $\mathfrak{Q}_{H}$'s \textbf{Singularity
Conservation Equation }\ref{eq:H-SCE}, we derive an \textbf{Auxiliary
Functional Equation} \ref{eq:AFE} satisfied by every $R\in\textrm{Ker}\left(1-\mathfrak{Q}_{H}\right)\subseteq L^{2}\left(\mathbb{Q}/\mathbb{Z}\right)$
\[
R\left(\varrho x\right)=e^{2\pi i\frac{b_{\iota}}{a_{\iota}}x}\sum_{n=0}^{\varrho-1}\xi_{\varrho}^{-n\iota}R\left(x+\frac{n}{\varrho}\right),\textrm{ }\forall x\in\mathbb{Q}
\]
where $\iota$ is any regulated index of $H$. An elementary argument
like that in the \textbf{Off-$\varrho$ Theorem} then establishes
yet another functional equation (the \textbf{Off-$\varrho$ Functional
Equation} (\textbf{ORFE}))\textbf{ }satisfied by every $R\in\textrm{Ker}\left(1-\mathfrak{Q}_{H}\right)\subseteq L^{2}\left(\mathbb{Q}/\mathbb{Z}\right)$
supported off $\varrho$:
\[
R\left(\varrho x\right)=e^{2\pi i\frac{b_{\iota}}{a_{\iota}}x}R\left(x\right),\textrm{ }\forall x\in\mathbb{Q}_{\bcancel{\varrho}}
\]
that is, the above equation holds for all $x\in\mathbb{Q}$ with $\left|x\right|_{\varrho}\leq1$.

\vphantom{}Step 3: To prove the \textbf{FDT} we show that being ``supported
off $\varrho$'' forces $R$'s support to either be trivial or infinite.
The argument is a ``Rob Peter to pay Saul'' scenario. Most of the
main features of the idea can be seen by considering the case of the
Collatz map itself. First, however, a definition:

\vphantom{}\textbf{Definition} \textbf{23}: Fix a $\varrho$-hydra
map $H$, and let $t\in\mathbb{Q}$. We write $S_{H}\left(t\right)$
to denote:
\begin{equation}
S_{H}\left(t\right)\overset{\textrm{def}}{=}\textrm{supp}\left(\mathfrak{Q}_{H}\left\{ \mathbf{1}_{t}\right\} \right)\label{eq:Definition of S_H (t)}
\end{equation}
In a minor abuse of terminology, we call $S_{H}\left(t\right)$ the
\textbf{image of $t$ under $\mathfrak{Q}_{H}$}; elements of $S_{H}\left(t\right)$
are said to be \textbf{in the image of $t$ (under $\mathfrak{Q}_{H}$)}.

\vphantom{}Now, let $R\in\textrm{Ker}\left(1-\mathfrak{Q}_{3}\right)$
be supported off $\varrho=2$. Then, $R\left(t\right)=0$ for all
rational numbers $t$ that, when written in simplest form, have an
even denominator. Moreover, letting $T$ denote $t\in\left[0,1\right)\cap\textrm{supp}\left(R\right)$
we can write $R\left(x\right)$ as the linear combination:
\[
R\left(x\right)=\sum_{t\in T}R\left(t\right)\mathbf{1}_{t}\left(x\right),\textrm{ }\forall x\in\mathbb{Q}
\]

So, consider what happens for the case of a particular value of $t\in\textrm{supp}\left(R\right)\subseteq\mathbb{Q}_{\bcancel{2}}$.
Suppose, for instance, that $t=\frac{1}{5}$ is in the support of
$R$. Then, since $R$ is fixed by $\mathfrak{Q}_{3}$, the $\mathbf{1}_{\frac{1}{5}}\left(x\right)$
present in $R$ will get sent to $\mathfrak{Q}_{3}\left\{ \mathbf{1}_{\frac{1}{5}}\right\} \left(x\right)$
by $\mathfrak{Q}_{3}$:
\[
R\left(x\right)=\mathfrak{Q}_{3}\left\{ R\right\} \left(x\right)=\sum_{t\in T}R\left(t\right)\mathfrak{Q}_{3}\left\{ \mathbf{1}_{t}\right\} \left(x\right)=R\left(\frac{1}{5}\right)\mathfrak{Q}_{3}\left\{ \mathbf{1}_{\frac{1}{5}}\right\} \left(x\right)+\sum_{t\in T\backslash\left\{ \frac{1}{5}\right\} }R\left(t\right)\mathfrak{Q}_{3}\left\{ \mathbf{1}_{t}\right\} \left(x\right)
\]
 Using the \textbf{Basis Formula} for $\mathfrak{Q}_{3}$:

\[
\mathfrak{Q}_{3}\left\{ \mathbf{1}_{t}\right\} \left(x\right)=\frac{1}{2}\left(\mathbf{1}_{\frac{t}{2}}\left(x\right)+\mathbf{1}_{\frac{t+1}{2}}\left(x\right)\right)+\frac{e^{-\pi it}}{2}\left(\mathbf{1}_{\frac{3t}{2}}\left(x\right)-\mathbf{1}_{\frac{3t+1}{2}}\left(x\right)\right)
\]
it is easy to show that, for any $t\in\mathbb{Q}$, every element
of $S_{H_{3}}\left(t\right)$ is congruent mod $1$ to one of $\frac{t}{2}$,
$\frac{t+1}{2}$, $\frac{3t}{2}$, or $\frac{3t+1}{2}$. As such,
the images of $t=\frac{1}{5}$ under $\mathfrak{Q}_{3}$ are:
\[
S_{H_{3}}\left(\frac{1}{5}\right)\overset{1}{\equiv}\left\{ \frac{1}{10},\frac{3}{5},\frac{3}{10},\frac{4}{5}\right\} 
\]
Since $\mathfrak{Q}_{3}\left\{ \mathbf{1}_{\frac{1}{5}}\right\} \left(x\right)$
is present in $R$, we see that $\frac{1}{10},\frac{3}{5},\frac{3}{10},\frac{4}{5}$
are all potentially elements of $\textrm{supp}\left(R\right)$. However,
$\frac{1}{10}$ and $\frac{3}{10}$ are both on $2$ (they have even
denominators). Since we know that $R$ is supported off $2$, it must
be that $R\left(\frac{1}{10}\right)=R\left(\frac{3}{10}\right)=0$.
Thus, the fact that $R\left(\frac{1}{10}\right)=0$ yet $\mathfrak{Q}_{3}\left\{ \mathbf{1}_{\frac{1}{5}}\right\} \left(\frac{1}{10}\right)\neq0$
tells us that there must be some $\tau\in\textrm{supp}\left(R\right)$
(with $\tau\overset{1}{\cancel{\equiv}}\frac{1}{5}$) so that the
image of $\tau$ under $\mathfrak{Q}_{3}$ contains $\frac{1}{10}$:
$\frac{1}{10}\in S_{H}\left(\tau\right),\textrm{ }\forall\tau\in T^{\prime}$.
The idea is that there will be a set $T^{\prime}\subseteq T$ of such
$\tau$s (all of which are distinct from one another mod $1$, with
$\frac{1}{5}\in T^{\prime}$) so that, for all $t\in T$:
\[
\mathfrak{Q}_{3}\left\{ \mathbf{1}_{t}\right\} \left(\frac{1}{10}\right)\neq0\Leftrightarrow t\in T^{\prime}
\]
and that:
\[
\sum_{\tau\in T^{\prime}}R\left(\tau\right)\mathfrak{Q}_{3}\left\{ \mathbf{1}_{\tau}\right\} \left(\frac{1}{10}\right)=0
\]
Consequently:
\[
R\left(\frac{1}{10}\right)=\sum_{t\in T\backslash T^{\prime}}R\left(t\right)\underbrace{\mathfrak{Q}_{3}\left\{ \mathbf{1}_{t}\right\} \left(\frac{1}{10}\right)}_{0}+\overbrace{\sum_{\tau\in T^{\prime}}\underbrace{R\left(\tau\right)\mathfrak{Q}_{3}\left\{ \mathbf{1}_{\tau}\right\} \left(\frac{1}{10}\right)}_{\neq0}}^{0}=0
\]
and order is restored. A set $T^{\prime\prime}\subseteq T$ must also
exist whose relation to $\frac{3}{10}$ is the same as that of $T^{\prime}$
to $\frac{1}{10}$.

Since the \textbf{Basis Formula }tells us \emph{exactly }what $S_{H}\left(\tau\right)$
looks like, a condition such as $\frac{1}{10}\in S_{H}\left(\tau\right)$
or $\frac{3}{10}\in S_{H}\left(\tau\right)$ allows us to explicitly
solve for what $\tau$ might be. For $\frac{1}{10}$, for instance,
the equations: 
\begin{eqnarray*}
\frac{\tau}{2} & = & \frac{1}{10}\\
\frac{\tau+1}{2} & = & \frac{1}{10}\\
\frac{3\tau}{2} & = & \frac{1}{10}\\
\frac{3\tau+1}{2} & = & \frac{1}{10}
\end{eqnarray*}
have as their solutions mod $1$: 
\[
\tau\in\left\{ \frac{1}{5},\frac{1}{15},\frac{11}{15}\right\} 
\]
Since $\tau$ had to be distinct from $\frac{1}{5}$ mod $1$, this
forces $\tau$ to be either $\frac{1}{15},\frac{11}{15}$. The key
observation to make is that both of these rational numbers have $15=3\times5$
as their denominators. That we multiplied them by $3$ is no accident.
This is the ``Rob Peter to pay Saul'' aspect of the argument. Here,
``Peter'' is $3$ (from $3x+1$; in particular, $3=\mu_{1}=\frac{\varrho a_{1}}{d_{1}}=\frac{2\times3}{2}$)
and ``Saul'' is $2$ (that is, $\varrho$). Starting with a $t$
($\frac{1}{5}$) which was ''off $2$'' (that is, its denominator
was odd), we saw that image of $\frac{1}{5}$ under $\mathfrak{Q}_{3}$
contained elements with even denominators ($\frac{1}{10}$, $\frac{3}{10}$).
In order to keep $R$ from violating its ``off-$2$-ness'', we had
to find at least one $\tau\in\textrm{supp}\left(R\right)$ whose image
under $\mathfrak{Q}_{3}$ also contained $\frac{1}{10}$, so that
the $\frac{1}{10}$ might be cancelled out (this is ``paying off
Saul''). In doing so, however, we found this forced $\tau$ to accrue
an extra factor of $3$ (stolen from Peter) which was not already
present in the original $t$ ($\frac{1}{5}$).

This argument can be repeated inductively: find an image of $t=\frac{1}{15}\in T$
with an even denominator; solve the congruences to show that there
must be a $\tau\in T$ whose denominator is a multiple of $9=3^{2}$.
At the $n$th step of this process, we find there must be an element
of $T$ whose denominator is a multiple of $3^{n}$. Consequently,
$T$ contains elements of arbitrarily large $3$-adic magnitude, which
then forces $T$ to be infinite. For a general $\varrho$-hydra map,
to keep $R$'s support off $\varrho$, $T$ must contain elements
whose denominators contain arbitrarily high powers of $\mu_{j}$ for
at least one non-regulated index $j\in\left\{ 0,\ldots,\varrho-1\right\} $.
That is, to pay off our debt of $\varrho$s to Saul, we must rob $\mu$s
from Peter's vault, and stash them in the denominators of elements
of $T$. The trick to making this work has to do with the \textbf{ORFE},
which can be used to show that if $T$ contains a non-zero element
(i.e., if there is a non-integer $t\in\mathbb{Q}\backslash\mathbb{Z}$
so that $R\left(t\right)\neq0$), then $T$ must contain a unit of
$\mathbb{Z}_{\varrho}$. As such, any finitely-supported $R$ must
be of the form $R\left(x\right)=R\left(0\right)\mathbf{1}_{\mathbb{Z}}\left(x\right)$,
exactly as asserted by the \textbf{FDT}, thus completing the proof.

\subsubsection{Formula Repository}

For ease of readability, the author has chosem to compile the bulk
of the proof's extraneous computations (and the resultant formulae)
in this subsection. For ease of access, we will merely \emph{state}
the formulae\emph{ }here, and leave their proofs at the end of the
of the Chapter (§3.4).

\vphantom{}\textbf{Iverson Bracket Formulae }(\textbf{IBF}): Let
$p$ be a prime number, let $q\in\left\{ 1\right\} \cup\mathbb{P}$,
and let $\nu$ be an integer. Then: 
\begin{equation}
\mathbf{1}_{pt}\left(qx\right)=\left[qx\overset{1}{\equiv}pt\right]=\sum_{k=0}^{q-1}\mathbf{1}_{\frac{pt+k}{q}}\left(x\right),\textrm{ }\forall x,t\in\mathbb{Q}\label{eq:Iverson bracket decomposition formula}
\end{equation}

\begin{equation}
\sum_{k=0}^{p-1}\xi_{p}^{-k\nu}\mathbf{1}_{t}\left(\frac{qx+k}{p}\right)=\xi_{p}^{\nu\left(qx-pt\right)}\left[qx\overset{1}{\equiv}pt\right],\textrm{ }\forall x,t\in\mathbb{Q}\label{eq:Adic integral computation formula}
\end{equation}

\begin{equation}
\left[\mathfrak{z}_{p}\overset{p}{\equiv}\nu\right]=\frac{1}{p}\sum_{k=0}^{p-1}\xi_{p}^{-\nu k}\check{\mathbf{1}}_{\frac{k}{p}}\left(\mathfrak{z}\right),\textrm{ }\forall\mathfrak{z}\in\widetilde{\mathbb{Z}}\label{eq: (rho,j) decomposition operator over Z-bar}
\end{equation}

\emph{Remark}: \ref{eq: (rho,j) decomposition operator over Z-bar}
is the analogue of the operator $\varpi_{\varrho,j}$ over $L^{2}\left(\widetilde{\mathbb{Z}}\right)$.

\textbf{\vphantom{}Basis Formula for $\mathfrak{Q}_{H}$} (\textbf{BF}-$\mathbb{Q}/\mathbb{Z}$):
Let $H:\mathbb{N}_{0}\rightarrow\mathbb{N}_{0}$ be a prime $\varrho$-hydra
map. Then, for all $t\in\mathbb{Q}$:
\begin{equation}
\mathfrak{Q}_{H}\left\{ \mathbf{1}_{t}\right\} \left(x\right)=\frac{1}{\varrho}\sum_{j=0}^{\varrho-1}e^{-2\pi i\frac{b_{j}}{d_{j}}t}\sum_{k=0}^{\varrho-1}\xi_{\varrho}^{jk}\mathbf{1}_{\frac{\mu_{j}t+k}{\varrho}}\left(x\right),\textrm{ }\forall x\in\mathbb{Q}\label{eq:Basis Formula}
\end{equation}

\vphantom{}We will also need to work with the analogue of $\mathfrak{Q}_{H}$
over $L^{2}\left(\widetilde{\mathbb{Z}}\right)$.

\vphantom{}\textbf{Definition} \textbf{24}: Let $\check{\mathfrak{Q}}_{H}:L^{2}\left(\widetilde{\mathbb{Z}}\right)\rightarrow L^{2}\left(\widetilde{\mathbb{Z}}\right)$
denote the linear operator obtained by conjugating $\mathfrak{Q}_{H}$
by the Fourier transform $\mathscr{F}_{\mathbb{Q}/\mathbb{Z}}:L^{2}\left(\mathbb{Q}/\mathbb{Z}\right)\rightarrow L^{2}\left(\widetilde{\mathbb{Z}}\right)$:
\[
\check{\mathfrak{Q}}_{H}\overset{\textrm{def}}{=}\mathscr{F}_{\mathbb{Q}/\mathbb{Z}}\circ\mathfrak{Q}_{H}\circ\mathscr{F}_{\mathbb{Q}/\mathbb{Z}}^{-1}
\]
For those who enjoy commutative diagrams, the situation is:
\[
\begin{array}{ccc}
L^{2}\left(\mathbb{Q}/\mathbb{Z}\right) & \overset{\mathfrak{Q}_{H}}{\rightarrow} & L^{2}\left(\mathbb{Q}/\mathbb{Z}\right)\\
\updownarrow &  & \updownarrow\\
L^{2}\left(\widetilde{\mathbb{Z}}\right) & \overset{\check{\mathfrak{Q}}_{H}}{\rightarrow} & L^{2}\left(\widetilde{\mathbb{Z}}\right)
\end{array}
\]
where the up arrows are $\mathscr{F}_{\mathbb{Q}/\mathbb{Z}}^{-1}$
and the down arrows are $\mathscr{F}_{\mathbb{Q}/\mathbb{Z}}$.

Essential here is the fact that $\mathscr{F}_{\mathbb{Q}/\mathbb{Z}}$
is an isometric isomorphism of Hilbert Spaces, preserving both magnitude
and orthogonality. These properties simplify the computations, by
allowing for things to be done over $L^{2}\left(\widetilde{\mathbb{Z}}\right)$,
where the arithmetical issues are more transparent.

\vphantom{}\textbf{Basis Formula for $\check{\mathfrak{Q}}_{H}$}
(\textbf{BF-}$\widetilde{\mathbb{Z}}$): Let $H:\mathbb{N}_{0}\rightarrow\mathbb{N}_{0}$
be a prime $\varrho$-hydra map. Then:
\begin{equation}
\check{\mathfrak{Q}}_{H}\left\{ \check{\mathbf{1}}_{t}\right\} \left(\mathfrak{z}\right)=\sum_{j=0}^{\varrho-1}\left[\mathfrak{z}_{\varrho}\overset{\varrho}{\equiv}-j\right]\check{\mathbf{1}}_{t}\left(\frac{a_{j}\mathfrak{z}-b_{j}}{d_{j}}\right),\textrm{ }\forall t\in\mathbb{Q}\label{eq:Basis formula for Q-check; prime rho}
\end{equation}

\vphantom{}\textbf{General Formula for $\check{\mathfrak{Q}}_{H}$}
(\textbf{GF-$\widetilde{\mathbb{Z}}$)}:

\begin{equation}
\check{\mathfrak{Q}}_{H}\left\{ F\right\} \left(\mathfrak{z}\right)=\sum_{j=0}^{\varrho-1}\left[\mathfrak{z}_{\varrho}\overset{\varrho}{\equiv}-j\right]F\left(\frac{a_{j}\mathfrak{z}-b_{j}}{d_{j}}\right),\textrm{ }\forall F\in L^{2}\left(\widetilde{\mathbb{Z}}\right)\label{eq:General formula for Q-check; prime rho}
\end{equation}

\pagebreak{}

\subsection{Proof of the Finite Dreamcatcher Theorem}

We begin by moving off $\varrho$.

\vphantom{}\textbf{Theorem 8: The Off-$\varrho$ Theorem}: Let $H:\mathbb{N}_{0}\rightarrow\mathbb{N}_{0}$
be a prime, regulated $\varrho$-hydra map, and let $R\in\textrm{Ker}\left(1-\mathfrak{Q}_{H}\right)$
have finite support. Then $\left|t\right|_{\varrho}\leq1$ for all
$t\in\textrm{supp}\left(R\right)$; equivalently, $R\left(t\right)=0$
for all $t\in\mathbb{Q}$ with $\left|t\right|_{\varrho}>1$.

Before we can get to the proof, we need to do a short computation.
Here, let everything be as given in \textbf{Theorem 8}:

\vphantom{}\textbf{Lemma} \textbf{4}: Let $T\overset{\textrm{def}}{=}\textrm{supp}\left(R\right)$,
enumerate the elements of $T$ as as $t_{1},\ldots,t_{N}$, and let
$\iota$ be a regulated index of $H$. Then:
\begin{equation}
\underbrace{\sum_{n=1}^{N}R\left(t_{n}\right)\mathbf{1}_{t_{n}}\left(x\right)}_{R\left(x\right)}=\sum_{n=1}^{N}R\left(t_{n}\right)e^{2\pi i\frac{\varrho^{m}-1}{\varrho-1}\frac{\varrho b_{\iota}+\left(\varrho-d_{\iota}\right)\iota}{d_{\iota}}t_{n}}\mathbf{1}_{\varrho^{m}t_{n}}\left(x\right),\textrm{ }\forall x\in\mathbb{Q},\textrm{ }\forall m\in\mathbb{N}_{0}\label{eq:R-functional equation; regulated; state 1}
\end{equation}

Proof: Since $R$ is finitely supported, we can write $R$ as:
\begin{equation}
R\left(x\right)=\sum_{n=1}^{N}R\left(t_{n}\right)\mathbf{1}_{t_{n}}\left(x\right)\label{eq:Basis formula for R}
\end{equation}
Taking the Fourier transform of \ref{eq:Basis formula for R} gives:
\begin{equation}
\check{R}\left(\mathfrak{z}\right)=\sum_{n=1}^{N}R\left(t_{n}\right)e^{2\pi i\left\langle t_{n},\mathfrak{z}\right\rangle }\label{eq:Basis Formula for R-check}
\end{equation}
Where, unlike the case of an $R$ with infinite support, the finite
support of our $R$ then guarantees that $\check{R}$ is defined for
all $\mathfrak{z}\in\widetilde{\mathbb{Z}}$. Since $R\in\textrm{Ker}\left(1-\mathfrak{Q}_{H}\right)$,
we know that $\check{R}\in\textrm{Ker}\left(1-\check{\mathfrak{Q}}_{H}\right)$,
and as such, by \ref{eq:General formula for Q-check; prime rho},
$\check{R}$ satisfies the functional equation:
\begin{equation}
\check{R}\left(\mathfrak{z}\right)=\sum_{j=0}^{\varrho-1}\left[\mathfrak{z}_{\varrho}\overset{\varrho}{\equiv}-j\right]\check{R}\left(\frac{a_{j}\mathfrak{z}-b_{j}}{d_{j}}\right),\textrm{ }\forall\mathfrak{z}\in\widetilde{\mathbb{Z}}\label{eq:Functional Equation for R-check}
\end{equation}

Now, since $H$ is regulated, there is an index $\iota\in\left\{ 0,\ldots,\varrho-1\right\} $
for which $a_{\iota}=1$. Replacing $\mathfrak{z}$ in \ref{eq:Functional Equation for R-check}
with $\varrho\mathfrak{z}-\iota$ gives:
\begin{equation}
\check{R}\left(\varrho\mathfrak{z}-\iota\right)=\sum_{j=0}^{\varrho-1}\left[\varrho\mathfrak{z}_{\varrho}-\iota\overset{\varrho}{\equiv}-j\right]\check{R}\left(\frac{a_{j}\left(\varrho\mathfrak{z}-\iota\right)-b_{j}}{d_{j}}\right)\label{eq:R-check functional equation; regulated index; state 0}
\end{equation}
Since:
\[
\left[\varrho\mathfrak{z}_{\varrho}-\iota\overset{\varrho}{\equiv}-j\right]=\left[\iota\overset{\varrho}{\equiv}j\right]
\]
the only non-zero term in the right-hand side of \ref{eq:R-check functional equation; regulated index; state 0}
is the $j=\iota$ term, which simplifies to:
\[
\check{R}\left(\frac{\overbrace{a_{\iota}}^{1}\left(\varrho\mathfrak{z}-\iota\right)-b_{\iota}}{d_{\iota}}\right)=\check{R}\left(\mathfrak{z}-\frac{\iota+b_{\iota}}{d_{\iota}}\right)
\]
leaving us with:
\begin{equation}
\check{R}\left(\varrho\mathfrak{z}-\iota\right)=\check{R}\left(\mathfrak{z}-\frac{\iota+b_{\iota}}{d_{\iota}}\right),\textrm{ }\forall\mathfrak{z}\in\widetilde{\mathbb{Z}}\label{eq:R-check functional equation; regulated index; state 1}
\end{equation}
Replacing $\mathfrak{z}$ with $\mathfrak{z}+\frac{\iota+b_{\iota}}{d_{\iota}}$
then yields:
\begin{equation}
\check{R}\left(\varrho\mathfrak{z}+\frac{\varrho b_{\iota}+\left(\varrho-d_{\iota}\right)\iota}{d_{\iota}}\right)=\check{R}\left(\mathfrak{z}\right),\textrm{ }\forall\mathfrak{z}\in\widetilde{\mathbb{Z}}\label{eq:R-check functional equation; regulated index; state 3}
\end{equation}
Consequently, letting:
\begin{equation}
\alpha\left(\mathfrak{z}\right)\overset{\textrm{def}}{=}\varrho\mathfrak{z}+\frac{\varrho b_{\iota}+\left(\varrho-d_{\iota}\right)\iota}{d_{\iota}}\label{eq:Definition of alpha for proof of Off-Rho Theorem}
\end{equation}
we have that $\check{R}$ is invariant under pre-compositions of $\alpha$.
Noting the pattern:
\begin{eqnarray*}
\alpha\left(\mathfrak{z}\right) & = & \varrho\mathfrak{z}+\frac{\varrho b_{\iota}+\left(\varrho-d_{\iota}\right)\iota}{d_{\iota}}\\
\alpha^{\circ2}\left(\mathfrak{z}\right) & = & \varrho^{2}\mathfrak{z}+\varrho\frac{\varrho b_{\iota}+\left(\varrho-d_{\iota}\right)\iota}{d_{\iota}}+\frac{\varrho b_{\iota}+\left(\varrho-d_{\iota}\right)\iota}{d_{\iota}}\\
\alpha^{\circ3}\left(\mathfrak{z}\right) & = & \varrho^{3}\mathfrak{z}+\varrho^{2}\frac{\varrho b_{\iota}+\left(\varrho-d_{\iota}\right)\iota}{d_{\iota}}+\varrho\frac{\varrho b_{\iota}+\left(\varrho-d_{\iota}\right)\iota}{d_{\iota}}+\frac{\varrho b_{\iota}+\left(\varrho-d_{\iota}\right)\iota}{d_{\iota}}\\
 & \vdots
\end{eqnarray*}
it follows that the $m$th iterate of $\alpha$ is given by:
\begin{equation}
\alpha^{\circ m}\left(\mathfrak{z}\right)=\varrho^{m}\mathfrak{z}+\frac{\varrho b_{\iota}+\left(\varrho-d_{\iota}\right)\iota}{d_{\iota}}\sum_{n=0}^{m-1}\varrho^{n}=\varrho^{m}\mathfrak{z}+\frac{\varrho^{m}-1}{\varrho-1}\frac{\varrho b_{\iota}+\left(\varrho-d_{\iota}\right)\iota}{d_{\iota}}\label{eq:mth iterate of alpha}
\end{equation}
Using \ref{eq:R-check functional equation; regulated index; state 3},
we can write:
\[
\check{R}\left(\mathfrak{z}\right)=\check{R}\left(\alpha\left(\mathfrak{z}\right)\right)=\check{R}\left(\alpha^{\circ2}\left(\mathfrak{z}\right)\right)=\cdots=\check{R}\left(\alpha^{\circ m}\left(\mathfrak{z}\right)\right)=\cdots
\]
for all $m$ and all $\mathfrak{z}$. With the above formula, we can
write $\check{R}\left(\mathfrak{z}\right)=\check{R}\left(\alpha^{\circ m}\left(\mathfrak{z}\right)\right)$
as:
\begin{equation}
\check{R}\left(\mathfrak{z}\right)=\check{R}\left(\varrho^{m}\mathfrak{z}+\frac{\varrho^{m}-1}{\varrho-1}\frac{\varrho b_{\iota}+\left(\varrho-d_{\iota}\right)\iota}{d_{\iota}}\right),\textrm{ }\forall\mathfrak{z}\in\widetilde{\mathbb{Z}}\label{eq:R-check functional equation; regulated index; state 4}
\end{equation}

Next, we prepare for the to return to $L^{2}\left(\mathbb{Q}/\mathbb{Z}\right)$.
Substituting \ref{eq:Basis Formula for R-check} into \ref{eq:R-check functional equation; regulated index; state 4}
we obtain:
\begin{equation}
\sum_{n=1}^{N}R\left(t_{n}\right)e^{2\pi i\left\langle t_{n},\mathfrak{z}\right\rangle }=\sum_{n=1}^{N}R\left(t_{n}\right)e^{2\pi i\left\langle t_{n},\varrho^{m}\mathfrak{z}+\frac{\varrho b_{\iota}+\left(\varrho-d_{\iota}\right)\iota}{d_{\iota}}\right\rangle }\label{eq:R-check functional equation; regulated index; state 5}
\end{equation}
Since the duality bracket $\left\langle \cdot,\cdot\right\rangle $
is $\mathbb{Z}$-bilinear, we can simplify the right-hand side's duality
bracket as follows:
\begin{eqnarray*}
\left\langle t_{n},\varrho^{m}\mathfrak{z}+\frac{\varrho^{m}-1}{\varrho-1}\frac{\varrho b_{\iota}+\left(\varrho-d_{\iota}\right)\iota}{d_{\iota}}\right\rangle  & = & \left\langle t_{n},\frac{\varrho^{m}-1}{\varrho-1}\frac{\varrho b_{\iota}+\left(\varrho-d_{\iota}\right)\iota}{d_{\iota}}\right\rangle +\left\langle t_{n},\varrho^{m}\mathfrak{z}\right\rangle \\
 & \overset{1}{\equiv} & t_{n}\frac{\varrho^{m}-1}{\varrho-1}\frac{\varrho b_{\iota}+\left(\varrho-d_{\iota}\right)\iota}{d_{\iota}}+\left\langle \varrho^{m}t_{n},\mathfrak{z}\right\rangle 
\end{eqnarray*}
and so, \ref{eq:R-check functional equation; regulated index; state 5}
becomes:
\begin{equation}
\sum_{n=1}^{N}R\left(t_{n}\right)e^{2\pi i\left\langle t_{n},\mathfrak{z}\right\rangle }=\sum_{n=1}^{N}R\left(t_{n}\right)e^{2\pi i\frac{\varrho^{m}-1}{\varrho-1}\frac{\varrho b_{\iota}+\left(\varrho-d_{\iota}\right)\iota}{d_{\iota}}t_{n}}e^{2\pi i\left\langle \varrho^{m}t_{n},\mathfrak{z}\right\rangle }\label{eq:R-check functional equation; regulated index; state 6}
\end{equation}
Applying $\mathscr{F}_{\mathbb{Q}/\mathbb{Z}}^{-1}$, the $e^{2\pi i\left\langle a,\mathfrak{z}\right\rangle }$s
get transformed to $\mathbf{1}_{a}\left(x\right)$s, and we obtain:
\[
\underbrace{\sum_{n=1}^{N}R\left(t_{n}\right)\mathbf{1}_{t_{n}}\left(x\right)}_{R\left(x\right)}=\sum_{n=1}^{N}R\left(t_{n}\right)e^{2\pi i\frac{\varrho^{m}-1}{\varrho-1}\frac{\varrho b_{\iota}+\left(\varrho-d_{\iota}\right)\iota}{d_{\iota}}t_{n}}\mathbf{1}_{\varrho^{m}t_{n}}\left(x\right)
\]
which holds for all $x\in\mathbb{Q}$ and all integers $m\geq0$,
as desired.

Q.E.D.

\vphantom{}\textbf{Proof of Theorem 8}: We proceed by way of contradiction:
suppose that there is a $t\in T$ with $\left|t\right|_{\varrho}>1$.
Without loss of generality, we can assume that $t_{1}$ is the element
of $T$ whose $\varrho$-adic magnitude is $>1$.

Since $t_{1}\in T$ implies $R\left(t_{1}\right)\neq0$, plugging
$x=t_{1}$ into \ref{eq:R-functional equation; regulated; state 1}
yields:
\[
0\neq R\left(t_{1}\right)=\sum_{n=1}^{N}R\left(t_{n}\right)e^{2\pi i\frac{\varrho^{m}-1}{\varrho-1}\frac{\varrho b_{\iota}+\left(\varrho-d_{\iota}\right)\iota}{d_{\iota}}t_{n}}\mathbf{1}_{\varrho^{m}t_{n}}\left(t_{1}\right)
\]
where $m\in\mathbb{N}_{0}$ is arbitrary. Note that if there is no
$n$ for which $\mathbf{1}_{\varrho^{m}t_{n}}\left(t_{1}\right)=1$,
the right-hand side evaluates to $0$, which is impossible. Thus,
for every $m$, there must be at least one $n$ for which $\mathbf{1}_{\varrho^{m}t_{n}}\left(t_{1}\right)=1$.
We denote this $n$ by $n_{m}$. Since the statement ``$\mathbf{1}_{\varrho^{m}t_{n_{m}}}\left(t_{1}\right)=1$''
is the same as `` ``$\varrho^{m}t_{n_{m}}-t_{1}\in\mathbb{Z}$'',
we have that $t_{n_{m}}$satisfies $\varrho^{m}t_{n_{m}}-t_{1}\in\mathbb{Z}$
for all $m\in\mathbb{N}_{0}$. Since every element of $\mathbb{Z}$
is also an element of $\mathbb{Z}_{\varrho}$, it must be that: 
\[
1\geq\left|\varrho^{m}t_{n_{m}}-t_{1}\right|_{\varrho}
\]
By the \textbf{Strong Triangle Inequality} property satisfied by the
$\varrho$-adic absolute value, either $\varrho^{m}t_{n_{m}}$ or
$t_{1}$ have the same $\varrho$-adic absolute value, or, the $\varrho$-adic
absolute value of their difference is equal to the maximum of their
respective $\varrho$-adic absolute values; that is:
\[
\left|\varrho^{m}t_{n_{m}}\right|_{\varrho}\neq\left|t_{1}\right|_{\varrho}\Rightarrow\left|\varrho^{m}t_{n_{m}}-t_{1}\right|_{\varrho}=\max\left\{ \left|\varrho^{m}t_{n_{m}}\right|_{\varrho},\left|t_{1}\right|_{\varrho}\right\} 
\]
Since it was given that $\left|t_{1}\right|_{\varrho}>1$, this implication
reads:
\[
\left|\varrho^{m}t_{n_{m}}\right|_{\varrho}\neq\left|t_{1}\right|_{\varrho}\Rightarrow1\geq\left|\varrho^{m}t_{n_{m}}-t_{1}\right|_{\varrho}=\max\left\{ \left|\varrho^{m}t_{n_{m}}\right|_{\varrho},\left|t_{1}\right|_{\varrho}\right\} >1
\]
which is impossible; $1$ cannot be strictly greater than itself.

Thus, $t_{1}\in T$ with $\left|t_{1}\right|_{\varrho}>1$ implies
that, for each $m\geq0$, there is a $t_{n_{m}}\in T$ so that:
\begin{eqnarray*}
\left|t_{1}\right|_{\varrho} & = & \left|\varrho^{m}t_{n_{m}}\right|_{\varrho}\\
 & = & \left|\varrho^{m}\right|_{\varrho}\left|t_{n_{m}}\right|_{\varrho}\\
 & = & \varrho^{-m}\left|t_{n_{m}}\right|_{\varrho}\\
 & \Updownarrow\\
\varrho^{m}\left|t_{1}\right|_{\varrho} & = & \left|t_{n_{m}}\right|_{\varrho}
\end{eqnarray*}
But then, this means that each of the $t_{n_{m}}$s must be distinct
w.r.t. $m$; indeed, using it, we can write:
\[
\left|t_{n_{m+1}}\right|_{\varrho}=\varrho^{m+1}\left|t_{1}\right|_{\varrho}=\varrho\times\varrho^{m}\left|t_{1}\right|_{\varrho}=\varrho\left|t_{n_{m}}\right|_{\varrho}
\]
Since this must hold for all $m$, this then forces $T$---the support
of $R$---to contain \emph{infinitely} many elements whenever there
is a $t_{1}\in T$ for which $\left|t_{1}\right|_{\varrho}>1$. But,
$T$ was given to be finite. This is a contradiction.

Thus, if $T$ is finite, $\left|t\right|_{\varrho}\leq1$ for all
$t\in T$, which proves the \textbf{Off-$\varrho$ Theorem}.

Q.E.D.

\vphantom{}Having shown that finitely-supported fixed points of $\mathfrak{Q}_{H}$
much vanish on $\varrho$, to complete the proof of the \textbf{FDT},
we need only show that this vanishing is sufficient to force the fixed
point to vanish for all non-integer inputs. To do this, we first need
to show how the \textbf{Off-$\varrho$ Theorem} affects the structure
of the support of a fixed point of $\mathfrak{Q}_{H}$.

\vphantom{}\textbf{Lemma 5: The Auxiliary Functional Equation} (\textbf{AFE}):
Let $H$ be a prime, regulated $\varrho$-hydra map, let $\iota$
be a regulated index of $H$, and let $R\in\textrm{Ker}\left(1-\mathfrak{Q}_{H}\right)$.
Then, $R$ satisfies the \textbf{Auxiliary Functional Equation}:
\begin{equation}
R\left(\varrho x\right)=e^{2\pi i\frac{b_{\iota}}{a_{\iota}}x}\sum_{n=0}^{\varrho-1}\xi_{\varrho}^{-n\iota}R\left(x+\frac{n}{\varrho}\right),\textrm{ }\forall x\in\mathbb{Q}\label{eq:AFE}
\end{equation}

Proof: Let $H$, $\iota$, and $R$ be as given. Now, taking (\ref{eq:H-SCE}):
\[
R\left(x\right)=\frac{1}{\varrho}\sum_{j=0}^{\varrho-1}e^{-2\pi i\frac{b_{j}}{a_{j}}x}\sum_{k=0}^{\mu_{j}-1}\xi_{\mu_{j}}^{-kH\left(j\right)}R\left(\frac{\varrho x+k}{\mu_{j}}\right)
\]
replace $x$ with $x+\frac{n}{\varrho}$, multiply both sides by $\xi_{\varrho}^{-n\iota}$,
and then sum everything from $n=0$ to $n=\varrho-1$. This produces:
\begin{eqnarray*}
\sum_{n=0}^{\varrho-1}\xi_{\varrho}^{-n\iota}R\left(x+\frac{n}{\varrho}\right) & = & \sum_{n=0}^{\varrho-1}\frac{\xi_{\varrho}^{-n\iota}}{\varrho}\sum_{j=0}^{\varrho-1}e^{-2\pi i\frac{b_{j}}{a_{j}}\left(x+\frac{n}{\varrho}\right)}\sum_{k=0}^{\mu_{j}-1}\xi_{\mu_{j}}^{-kH\left(j\right)}R\left(\frac{\varrho\left(x+\frac{n}{\varrho}\right)+k}{\mu_{j}}\right)\\
 & = & \sum_{n=0}^{\varrho-1}\frac{\xi_{\varrho}^{-n\iota}}{\varrho}\sum_{j=0}^{\varrho-1}e^{-2\pi i\frac{b_{j}}{a_{j}}x}e^{-2\pi i\frac{b_{j}}{a_{j}}\frac{n}{\varrho}}\sum_{k=0}^{\mu_{j}-1}\xi_{\mu_{j}}^{-kH\left(j\right)}R\left(\frac{\varrho x+n+k}{\mu_{j}}\right)\\
\left(m=n+k\right); & = & \sum_{n=0}^{\varrho-1}\frac{\xi_{\varrho}^{-n\iota}}{\varrho}\sum_{j=0}^{\varrho-1}e^{-2\pi i\frac{b_{j}}{a_{j}}x}e^{-2\pi i\frac{b_{j}}{d_{j}}\frac{nd_{j}}{\varrho a_{j}}}\sum_{m=0}^{\mu_{j}-1}\xi_{\mu_{j}}^{\left(n-m\right)H\left(j\right)}R\left(\frac{\varrho x+m}{\mu_{j}}\right)\\
\left(\frac{d_{j}}{\varrho a_{j}}=\frac{1}{\mu_{j}}\right); & = & \sum_{n=0}^{\varrho-1}\frac{\xi_{\varrho}^{-n\iota}}{\varrho}\sum_{j=0}^{\varrho-1}\frac{e^{-\frac{2\pi in}{\mu_{j}}\frac{b_{j}}{d_{j}}}\xi_{\mu_{j}}^{nH\left(j\right)}}{e^{\frac{2\pi ib_{j}x}{a_{j}}}}\sum_{m=0}^{\mu_{j}-1}\xi_{\mu_{j}}^{-mH\left(j\right)}R\left(\frac{\varrho x+m}{\mu_{j}}\right)\\
 & = & \sum_{n=0}^{\varrho-1}\frac{\xi_{\varrho}^{-n\iota}}{\varrho}\sum_{j=0}^{\varrho-1}\frac{e^{\frac{2\pi in}{\mu_{j}}\left(H\left(j\right)-\frac{b_{j}}{d_{j}}\right)}}{e^{2\pi i\frac{b_{j}}{a_{j}}x}}\sum_{m=0}^{\mu_{j}-1}\xi_{\mu_{j}}^{-mH\left(j\right)}R\left(\frac{\varrho x+m}{\mu_{j}}\right)\\
\left(H\left(j\right)=\frac{ja_{j}+b_{j}}{d_{j}}\right); & = & \sum_{n=0}^{\varrho-1}\frac{\xi_{\varrho}^{-n\iota}}{\varrho}\sum_{j=0}^{\varrho-1}\frac{e^{\frac{2\pi in}{\mu_{j}}\left(\frac{ja_{j}+b_{j}}{d_{j}}-\frac{b_{j}}{d_{j}}\right)}}{e^{2\pi i\frac{b_{j}}{a_{j}}x}}\sum_{m=0}^{\mu_{j}-1}\xi_{\mu_{j}}^{-mH\left(j\right)}R\left(\frac{\varrho x+m}{\mu_{j}}\right)\\
\left(\frac{1}{\mu_{j}}=\frac{d_{j}}{\varrho a_{j}}\right); & = & \sum_{n=0}^{\varrho-1}\frac{\xi_{\varrho}^{-n\iota}}{\varrho}\sum_{j=0}^{\varrho-1}\frac{e^{\frac{2\pi ijn}{\varrho}}}{e^{2\pi i\frac{b_{j}}{a_{j}}x}}\sum_{m=0}^{\mu_{j}-1}\xi_{\mu_{j}}^{-mH\left(j\right)}R\left(\frac{\varrho x+m}{\mu_{j}}\right)\\
 & = & \sum_{n=0}^{\varrho-1}\frac{\xi_{\varrho}^{-n\iota}}{\varrho}\sum_{j=0}^{\varrho-1}e^{-2\pi i\frac{b_{j}}{a_{j}}x}\xi_{\varrho}^{jn}\sum_{m=0}^{\mu_{j}-1}\xi_{\mu_{j}}^{-mH\left(j\right)}R\left(\frac{\varrho x+m}{\mu_{j}}\right)\\
 & = & \sum_{j=0}^{\varrho-1}e^{-2\pi i\frac{b_{j}}{a_{j}}x}\underbrace{\left(\frac{1}{\varrho}\sum_{n=0}^{\varrho-1}\xi_{\varrho}^{n\left(j-\iota\right)}\right)}_{\left[j\overset{\varrho}{\equiv}\iota\right]}\sum_{m=0}^{\mu_{j}-1}\xi_{\mu_{j}}^{-mH\left(j\right)}R\left(\frac{\varrho x+m}{\mu_{j}}\right)\\
\left(\textrm{only }j=\iota\textrm{ stays}\right); & = & e^{-2\pi i\frac{b_{\iota}}{a_{\iota}}x}\sum_{m=0}^{\mu_{\iota}-1}\xi_{\mu_{\iota}}^{-mH\left(\iota\right)}R\left(\frac{\varrho x+m}{\mu_{\iota}}\right)
\end{eqnarray*}

Since $\mu_{\iota}=1$, we are left with:
\[
\sum_{n=0}^{\varrho-1}\xi_{\varrho}^{-n\iota}R\left(x+\frac{n}{\varrho}\right)=e^{-\frac{2\pi ib_{\iota}}{a_{\iota}}x}R\left(\varrho x\right)
\]
Multiplying by $e^{\frac{2\pi ib_{\iota}}{a_{\iota}}x}$ on both sides
then gives:
\[
R\left(\varrho x\right)=e^{2\pi i\frac{b_{\iota}}{a_{\iota}}x}\sum_{n=0}^{\varrho-1}\xi_{\varrho}^{-n\iota}R\left(x+\frac{n}{\varrho}\right)
\]
which is the desired functional equation.

Q.E.D.

\vphantom{}\textbf{Corollary 4: The Off-$\varrho$ Functional Equation}
(\textbf{ORFE}): Let $H$ be a prime, regulated $\varrho$-hydra map,
let $\iota$ be a regulated index of $H$, and let $R\in\textrm{Ker}\left(1-\mathfrak{Q}_{H}\right)$
be supported off $\varrho$ (that is, $\left|t\right|_{\varrho}\leq1$
for all $t\in\textrm{supp}\left(R\right)$). Then:
\begin{equation}
R\left(\varrho x\right)=e^{2\pi i\frac{b_{\iota}}{a_{\iota}}x}R\left(x\right),\textrm{ }\forall x\in\mathbb{Q}_{\bcancel{\varrho}}\label{eq:off-rho functional equation}
\end{equation}

Proof: The proof begins with the \textbf{AFE }\ref{eq:AFE}:
\[
R\left(\varrho x\right)=e^{2\pi i\frac{b_{\iota}}{a_{\iota}}x}\sum_{n=0}^{\varrho-1}\xi_{\varrho}^{-n\iota}R\left(x+\frac{n}{\varrho}\right)
\]
Now fix $x\in\mathbb{Q}_{\bcancel{\varrho}}$, and write $x$ in irreducible
form as $x=\frac{\beta}{\alpha}$. By the \textbf{Off-$\varrho$ Theorem},
we know that $R\left(t\right)$ necessarily vanishes for all $t\in\mathbb{Q}$
with $\left|t\right|_{\varrho}>1$. Thus:
\[
R\left(\frac{\beta}{\alpha}+\frac{n}{\varrho}\right)=0\textrm{ }\forall n:\left|\frac{\beta}{\alpha}+\frac{n}{\varrho}\right|_{\varrho}>1
\]
Since $x$ is off $\varrho$, $\gcd\left(\alpha,\varrho\right)=1$.
As such, $\left|\alpha\right|_{\varrho}=1$, and so:
\[
\left|\frac{\beta}{\alpha}+\frac{n}{\varrho}\right|_{\varrho}=\left|\frac{\beta\varrho+n\alpha}{\varrho\alpha}\right|_{\varrho}=\frac{\left|\beta\varrho+n\alpha\right|_{\varrho}}{\left|\varrho\alpha\right|_{\varrho}}=\frac{\left|\beta\varrho+n\alpha\right|_{\varrho}}{\frac{1}{\varrho}}=\varrho\left|\beta\varrho+n\alpha\right|_{\varrho}
\]
Note that:
\begin{eqnarray*}
\varrho\left|\beta\varrho+n\alpha\right|_{\varrho} & \leq & 1\\
 & \Updownarrow\\
\left|\beta\varrho+n\alpha\right|_{\varrho} & \leq & \frac{1}{\varrho}\\
 & \Updownarrow\\
\underbrace{\beta\varrho}_{0}+n\alpha & \overset{\varrho}{\equiv} & 0\\
\left(\gcd\left(\alpha,\varrho\right)=1\right); & \Updownarrow\\
n & \overset{\varrho}{\equiv} & 0\\
\left(n\in\left\{ 0,\ldots,\varrho-1\right\} \right); & \Updownarrow\\
n & = & 0
\end{eqnarray*}
Thus, for every $x$ off $\varrho$, the \textbf{AFE }becomes:
\[
R\left(\varrho x\right)=e^{2\pi i\frac{b_{\iota}}{a_{\iota}}x}\sum_{n=0}^{\varrho-1}\xi_{\varrho}^{-n\iota}R\left(x+\frac{n}{\varrho}\right)=e^{2\pi i\frac{b_{\iota}}{a_{\iota}}x}R\left(x\right)
\]
since the terms in the $n$-sum vanish for all $n$ except $n=0$.
This is the desired functional equation.

Q.E.D.

\vphantom{}\emph{Remark}:\emph{ }In particular, this shows that:
\[
x\in\textrm{supp}\left(R\right)\Rightarrow\left[\varrho x\right]_{1}\in\textrm{supp}\left(R\right)
\]
where, recall, $\left[\varrho x\right]_{1}$ is the unique element
of $\left[0,1\right)$ which is congruent to $\varrho x$ mod $1$.
Consequently, the \textbf{ORFE }implies that \emph{the support of
$R$ is invariant under multiplication by $\varrho$ mod $1$}.

\vphantom{}\textbf{Proposition 12}: Let $H$ be a regulated $\varrho$-hydra
map, let $R\in\textrm{Ker}\left(1-\mathfrak{Q}_{H}\right)$ is supported
off $\varrho$, and let $T\overset{\textrm{def}}{=}\textrm{supp}\left(R\right)\cap\left[0,1\right)$.
If there is a non-integer element $x_{0}\in T$, then $T$ contains
an element $\tau$ with $\left|\tau\right|_{\varrho}=1$ (i.e., both
the numerator and denominator of $\tau$ are co-prime to $\varrho$).

Proof: Let $H$, $R$, and $x_{0}$ be as given. Then, we can write
$x_{0}\in\left(0,1\right)$ in simplest form as $\frac{\beta}{\alpha}$.
Since $R$ is supported off $\varrho$, $\alpha$ must be co-prime
to $\varrho$. Moreover, it follows that $R$ satisfies the \textbf{ORFE};
as such, $\varrho^{n}x_{0}$ mod $1$ is in $T$, for all $n\in\mathbb{N}_{0}$.
Thus, if $\beta$ is a multiple of $\varrho$, we can multiply it
by $\varrho$ a finite number of times (namely, by any $n$ so that
$\varrho^{n}\overset{\alpha}{\equiv}\left[\varrho^{\textrm{val}_{\varrho}\left(\beta\right)}\right]_{\alpha}^{-1}$)
and eventually obtain an element of $T$ whose numerator \emph{and
}denominator are both co-prime to $\varrho$. This element is the
desired $\tau$.

Q.E.D.

\vphantom{}\textbf{Proposition 13}:\textbf{ }Fix a $\tau$ as described
in \textbf{Proposition 12}. Then, there must be a $t\in T\backslash\left\{ \tau\right\} $
so that $\frac{\tau}{\varrho}\in S_{H}\left(\tau\right)\cap S_{H}\left(t\right)$.
(That is, the term on $\varrho$ created by $\mathfrak{Q}_{H}$ at
$\tau$ must be cancelled out by at least one other term created by
$\mathfrak{Q}_{H}$ at an element of $T$ \emph{other than }$t$.)

Proof: Using the \textbf{Basis Formula }for $\mathfrak{Q}_{H}$, observe
that:
\[
\mathfrak{Q}_{H}\left\{ \mathbf{1}_{\tau}\right\} \left(x\right)=\frac{1}{\varrho}\sum_{j=0}^{\varrho-1}e^{-2\pi i\frac{b_{j}}{d_{j}}\tau}\sum_{k=0}^{\varrho-1}\xi_{\varrho}^{jk}\mathbf{1}_{\frac{\mu_{j}\tau+k}{\varrho}}\left(x\right)
\]
Letting 
\[
f_{\tau:j}\left(x\right)\overset{\textrm{def}}{=}\frac{1}{\varrho}\sum_{k=0}^{\varrho-1}\xi_{\varrho}^{jk}\mathbf{1}_{\frac{\mu_{j}\tau+k}{\varrho}}\left(x\right)
\]
observe that:
\begin{eqnarray*}
\check{f}_{\tau:j}\left(\mathfrak{z}\right) & = & \frac{1}{\varrho}\sum_{k=0}^{\varrho-1}\xi_{\varrho}^{jk}\check{\mathbf{1}}_{\frac{\mu_{j}\tau+k}{\varrho}}\left(\mathfrak{z}\right)\\
 & = & \left(\frac{1}{\varrho}\sum_{k=0}^{\varrho-1}\xi_{\varrho}^{jk}\check{\mathbf{1}}_{\frac{k}{\varrho}}\left(\mathfrak{z}\right)\right)\check{\mathbf{1}}_{\frac{\mu_{j}\tau}{\varrho}}\left(\mathfrak{z}\right)\\
\left(\textrm{use \ref{eq:Iverson bracket decomposition formula}}\right); & = & \left[\mathfrak{z}_{\varrho}\overset{\varrho}{\equiv}-j\right]\check{\mathbf{1}}_{\frac{\mu_{j}\tau}{\varrho}}\left(\mathfrak{z}\right)
\end{eqnarray*}
and so:
\[
\check{\mathfrak{Q}}_{H}\left\{ \check{\mathbf{1}}_{\tau}\right\} \left(\mathfrak{z}\right)=\sum_{j=0}^{\varrho-1}e^{-2\pi i\frac{b_{j}}{d_{j}}\tau}\underbrace{\left[\mathfrak{z}_{\varrho}\overset{\varrho}{\equiv}-j\right]\check{\mathbf{1}}_{\frac{\mu_{j}\tau}{\varrho}}\left(\mathfrak{z}\right)}_{\check{f}_{\tau:j}\left(\mathfrak{z}\right)}
\]
Since the Iverson brackets are pair-wise orthogonal in $L^{2}\left(\widetilde{\mathbb{Z}}\right)$
with respect to $j$, so too are the $\check{f}_{\tau:j}$s, and hence---since
the Fourier Transform is an isometric isomorphism of Hilbert Spaces---so
too are the $f_{\tau:j}$s over $L^{2}\left(\mathbb{Q}/\mathbb{Z}\right)$.
Thus, the supports of the $f_{\tau:j}$s are pair-wise disjoint with
respect to $j$. As such:
\[
\textrm{supp}\left(\mathfrak{Q}_{H}\left\{ \mathbf{1}_{\tau}\right\} \right)=\bigcup_{j=0}^{\varrho-1}\textrm{supp}\left(f_{\tau:j}\right)
\]

Next, since $H$ is regulated, there is a regulated index $\iota\in\left\{ 0,\ldots,\varrho-1\right\} $
so that $\mu_{\iota}=1$. So, picking $j=\iota$, we find that:
\[
f_{\tau:\iota}\left(x\right)=\frac{1}{\varrho}\sum_{k=0}^{\varrho-1}\xi_{\varrho}^{jk}\mathbf{1}_{\frac{\tau+k}{\varrho}}\left(x\right)
\]
Since the rational numbers $\frac{\tau-k}{\varrho}$ are pair-wise
distinct mod $1$ with respect to $k$, upon picking $k=0$, we find
that $\frac{\tau}{\varrho}\in\textrm{supp}\left(f_{\tau:\iota}\right)\subseteq\textrm{supp}\left(\mathfrak{Q}_{H}\left\{ \mathbf{1}_{\tau}\right\} \right)$.
But this is a problem! Since:
\[
R\left(x\right)=\mathfrak{Q}_{H}\left\{ R\right\} \left(x\right)=\sum_{t\in T}R\left(t\right)\mathfrak{Q}_{H}\left\{ \mathbf{1}_{t}\right\} \left(x\right)
\]
pulling out $\tau$ from $T$ yields:
\[
R\left(x\right)=\sum_{t\in T}R\left(t\right)\mathfrak{Q}_{H}\left\{ \mathbf{1}_{t}\right\} \left(x\right)=R\left(\tau\right)\mathfrak{Q}_{H}\left\{ \mathbf{1}_{\tau}\right\} \left(x\right)+\sum_{t\in T\backslash\left\{ \tau\right\} }R\left(t\right)\mathfrak{Q}_{H}\left\{ \mathbf{1}_{t}\right\} \left(x\right),\textrm{ }\forall x\in\mathbb{Q}
\]
where $R\left(\tau\right)\neq0$, since $\tau\in T$. Setting $x=\frac{\tau}{\varrho}$
gives:
\[
R\left(\frac{\tau}{\varrho}\right)=\sum_{t\in T}R\left(t\right)\mathfrak{Q}_{H}\left\{ \mathbf{1}_{t}\right\} \left(\frac{\tau}{\varrho}\right)=\underbrace{R\left(\tau\right)}_{\neq0}\underbrace{\mathfrak{Q}_{H}\left\{ \mathbf{1}_{\tau}\right\} \left(\frac{\tau}{\varrho}\right)}_{\neq0}+\sum_{t\in T\backslash\left\{ \tau\right\} }R\left(t\right)\mathfrak{Q}_{H}\left\{ \mathbf{1}_{t}\right\} \left(\frac{\tau}{\varrho}\right)
\]
So, if:
\[
\mathfrak{Q}_{H}\left\{ \mathbf{1}_{t}\right\} \left(\frac{\tau}{\varrho}\right)=0,\textrm{ }\forall t\in T\backslash\left\{ \tau\right\} 
\]
we then have that 
\[
R\left(\frac{\tau}{\varrho}\right)=R\left(\tau\right)\mathfrak{Q}_{H}\left\{ \mathbf{1}_{\tau}\right\} \left(\frac{\tau}{\varrho}\right)\neq0
\]
and thus---since $\left|\tau\right|_{\varrho}=1$---that $\frac{\tau}{\varrho}$
is an element in the support of $R$ which is on $\varrho$, which
contradicts the fact that $R$ is supported off $\varrho$.

Thus, in order to avoid a contradiction, it must be that there is
a $t\in T\backslash\left\{ \tau\right\} $ for which $\mathfrak{Q}_{H}\left\{ \mathbf{1}_{t}\right\} \left(\frac{\tau}{\varrho}\right)\neq0$;
that is, there is a $t\in T\backslash\left\{ \tau\right\} $ so that:
\[
\frac{\tau}{\varrho}\in S_{H}\left(\tau\right)\cap S_{H}\left(t\right)
\]
as desired.

Q.E.D.

\vphantom{}\textbf{Theorem 9 (The Off-$\varrho$ Infinity Theorem})\textbf{:}
Let $H$ be a prime, regulated $\varrho$-hydra map. Suppose $R\in\textrm{Ker}\left(1-\mathfrak{Q}_{H}\right)$
is supported off $\varrho$, and that $R$ is not of the form:
\[
R\left(x\right)=R\left(0\right)\mathbf{1}_{\mathbb{Z}}\left(x\right),\textrm{ }\forall x\in\mathbb{Q}
\]
Then, $R$ is \emph{not} finitely supported.

Proof: Let $H$ and $R$ be as given, and suppose that $R$ is not
of the form $R\left(x\right)=R\left(0\right)\mathbf{1}_{\mathbb{Z}}\left(x\right)$.
Since $R\left(x\right)=R\left(0\right)\mathbf{1}_{\mathbb{Z}}\left(x\right)$
if and only if $\textrm{supp}\left(R\right)\subseteq\mathbb{Z}$,
it must be that there is an non-integer $x_{0}\in\mathbb{Q}\backslash\mathbb{Z}$
contained in $\textrm{supp}\left(R\right)$. Since $R$ is supported
off $\varrho$, it must be that $\left|x_{0}\right|_{\varrho}\leq1$.

Now, let $T\overset{\textrm{def}}{=}\textrm{supp}\left(R\right)\cap\left[0,1\right)$.
Then, applying the argument from \textbf{Proposition 12}, after multiplying
$x_{0}$ by a sufficiently large power of $\varrho$ and simplifying
mod $1$, we obtain a $\tau\in T$ with \textbf{$\left|\tau\right|_{\varrho}=1$}.
Next, by \textbf{Proposition 13}, there is a $t\in T\backslash\left\{ \tau\right\} $
so that $\frac{\tau}{\varrho}\in S_{H}\left(\tau\right)\cap S_{H}\left(t\right)$.
Turning to the \textbf{Basis Formula }for $\mathfrak{Q}_{H}$, note
that:

\[
0\neq\mathfrak{Q}_{H}\left\{ \mathbf{1}_{t}\right\} \left(\frac{\tau}{\varrho}\right)=\frac{1}{\varrho}\sum_{j=0}^{\varrho-1}e^{-2\pi i\frac{b_{j}}{d_{j}}t}\sum_{k=0}^{\varrho-1}\xi_{\varrho}^{jk}\mathbf{1}_{\frac{\mu_{j}t+k}{\varrho}}\left(\frac{\tau}{\varrho}\right)
\]
Hence, there must be $j,k\in\left\{ 0,\ldots,\varrho-1\right\} $
for which:
\begin{eqnarray*}
\mathbf{1}_{\frac{\mu_{j}t+k}{\varrho}}\left(\frac{\tau}{\varrho}\right) & = & 1\\
 & \Updownarrow\\
\frac{\mu_{j}t+k}{\varrho} & \overset{1}{\equiv} & \frac{\tau}{\varrho}\\
\left(\exists\nu\in\mathbb{Z}\right); & \Updownarrow\\
\mu_{j}t+\varrho\nu+k & = & \tau
\end{eqnarray*}
else, $\mathfrak{Q}_{H}\left\{ \mathbf{1}_{t}\right\} \left(\frac{\tau}{\varrho}\right)$
would vanish, which is not allowed. So: 
\[
\mu_{j}t+\underbrace{\varrho\nu+k}_{\in\mathbb{Z}}=\tau\in T\subseteq\textrm{supp}\left(R\right)
\]
Now, note that $j$ \emph{cannot }be a regulated index; that is, $\mu_{j}\neq1$.
Indeed, were $\mu_{j}=1$, then we would have $t=\tau-\left(\varrho\nu+k\right)$,
which implies that $t$ is congruent to $\tau$ mod $1$. However,
since \emph{both} $t$ and $\tau$ are in $T\subseteq\left[0,1\right)$,
$t\overset{1}{\equiv}\tau$ forces $t=\tau$, which is impossible
($t\in T\backslash\left\{ \tau\right\} $).

So, with $\mu_{j}\neq1$, since $R$ is $1$-periodic, and since $\varrho\nu+k\in\mathbb{Z}$,
it follows that $\mu_{j}t+\varrho\nu+k\in\textrm{supp}\left(R\right)$
occurs if and only if: 
\[
\mu_{j}t=\mu_{j}t+\varrho\nu+k-\left(\varrho\nu+k\right)\in\textrm{supp}\left(R\right)
\]
Thus:
\[
\tau^{\prime}\overset{\textrm{def}}{=}\tau-\left(\varrho\nu+k\right)
\]
is an element of $\textrm{supp}\left(R\right)$ so that $\mu_{j}t=\tau^{\prime}$.
Since $t$ was guaranteed to be in $T$ (and thus, in $\textrm{supp}\left(R\right)$),
this shows that: 
\[
t=\frac{\tau^{\prime}}{\mu_{j}}
\]
Since $\tau$ was obtained from a non-integer $x_{0}\in\textrm{supp}\left(R\right)$
with $\left|x_{0}\right|_{\varrho}\leq1$ by multiplication by power
of $\varrho$ and reduction mod $1$, observe that $\tau$ is necessarily
a non-integer as well, and---since $\tau$ is congruent to $\tau^{\prime}$
mod $1$---so is $\tau^{\prime}$. In particular, $\tau^{\prime}\neq0$.
Since $\mu_{j}\neq1$, it then follows that the $\mu_{j}$-adic absolute
value of $t$ is non-trivial, and is given by: 
\[
\left|t\right|_{\mu_{j}}=\left|\frac{\tau^{\prime}}{\mu_{j}}\right|_{\mu_{j}}=\mu_{j}\left|\tau^{\prime}\right|_{\mu_{j}}
\]
Thus, just as was foretold, in order to pay off our $\varrho$-debt
to Saul, there must be $t,\tau^{\prime}\in\textrm{supp}\left(R\right)$
and a non-regulated index $j\in\left\{ 0,\ldots,\varrho-1\right\} $
so that $\left|t\right|_{\mu_{j}}=\mu_{j}\left|\tau^{\prime}\right|_{\mu_{j}}$.

In summary: if $R$ (as given in the statement of the theorem) vanishes
on a non-zero equivalence class in $\mathbb{Q}/\mathbb{Z}$ ($x_{0}$),
we can obtain $t,\tau^{\prime}\in\textrm{supp}\left(R\right)\backslash\mathbb{Z}$,
with $\left|t^{\prime}\right|_{\mu_{j}}=\mu_{j}\left|\tau^{\prime}\right|_{\mu_{j}}$.
We can then proceed inductively; let $x_{1}=t$; then, we obtain $t^{\prime},\tau^{\prime\prime}\in\textrm{supp}\left(R\right)\backslash\mathbb{Z}$
so that $\left|t^{\prime}\right|_{\mu_{j^{\prime}}}=\mu_{j^{\prime}}\left|\tau^{\prime\prime}\right|_{\mu_{j^{\prime}}}$
for some non-regulated index $j^{\prime}$. Moreover, in this construction,
we have that $\tau^{\prime\prime}\overset{1}{\equiv}\tau^{\prime}$,
and thus, that $\left|\tau^{\prime\prime}\right|_{p}=\left|\tau^{\prime}\right|_{p}$
for all primes $p$. Proceeding in this manner, we obtain an infinite
sequence $\left\{ x_{n}\right\} _{n\in\mathbb{N}_{0}}\subseteq\textrm{supp}\left(R\right)\backslash\mathbb{Z}$
with the property that, for each $n$, there is a non-regulated index
$j_{n}\in\left\{ 0,\ldots,\varrho-1\right\} $ so that: 
\[
\left|x_{n+1}\right|_{\mu_{j_{n}}}=\mu_{j_{n}}\left|x_{n}\right|_{\mu_{j_{n}}}
\]
Since there are only finitely many possible values for the $j_{n}$s
to take, the \textbf{Pigeonhole Principle }guarantees that as $n\rightarrow\infty$,
there must be at least one non-regulated index $j\in\left\{ 0,\ldots,\varrho-1\right\} $
so that $j_{n}=j$ for infinitely many $n$. Consequently, the $x_{n}$s
contain a subsequence which is strictly increasing in $\mu_{j}$-adic
magnitude; this forces all the elements of that subsequence to be
distinct. Thus, $\textrm{supp}\left(R\right)\backslash\mathbb{Z}\subseteq\textrm{supp}\left(R\right)$
contains infinitely many elements; that is, $\textrm{supp}\left(R\right)$
is \emph{not }finitely supported---as desired.

Q.E.D.

\vphantom{}\textbf{Proof of the Finite Dreamcatcher Theorem}: Let
$H$ be a prime, regulated $\varrho$-hydra map, let $\iota$ be a
regulated index of $H$, let $R\in\textrm{Ker}\left(1-\mathfrak{Q}_{H}\right)$,
and suppose that $R$ has finite support. By the \textbf{Off-$\varrho$
Theorem}, the finiteness of $R$'s support forces $R$ to be supported
off $\varrho$. Consequently, the \textbf{Off-$\varrho$ Infinity
Theorem} then forces $R$ to be of the form:
\[
R\left(x\right)=R\left(0\right)\mathbf{1}_{\mathbb{Z}}\left(x\right),\textrm{ }\forall x\in\mathbb{Q}
\]
else, $R$ would not be finitely supported---a contradiction. Thus,
$R$ has exactly the form asserted by the\textbf{ FDT}. This proves
the \textbf{FDT}.

Q.E.D.

\pagebreak{}

\subsection{Proofs of the Formulae from §3.2.2}

\subsubsection*{Proof of IBF}

(I) Proofs of \ref{eq:Iverson bracket decomposition formula} and
\ref{eq:Adic integral computation formula}:
\begin{eqnarray*}
\sum_{k=0}^{p-1}\xi_{p}^{-k\nu}\mathbf{1}_{t}\left(\frac{qx+k}{p}\right) & = & \int_{\widetilde{\mathbb{Z}}}\sum_{k=0}^{p-1}\xi_{p}^{-k\nu}e^{-2\pi i\left\langle \frac{qx+k}{p}-t,\mathfrak{z}\right\rangle }d\mathfrak{z}\\
 & = & \int_{\widetilde{\mathbb{Z}}}\sum_{k=0}^{p-1}\xi_{p}^{-k\nu}e^{-2\pi i\left\langle \frac{k}{p},\mathfrak{z}\right\rangle }e^{-2\pi i\left\langle \frac{qx}{p}-t,\mathfrak{z}\right\rangle }d\mathfrak{z}\\
 & = & \int_{\widetilde{\mathbb{Z}}}\sum_{k=0}^{p-1}\xi_{p}^{-k\nu}\xi_{p}^{-k\left[\mathfrak{z}_{p}\right]_{p}}e^{-2\pi i\left\langle \frac{qx}{p}-t,\mathfrak{z}\right\rangle }d\mathfrak{z}\\
 & = & \int_{\widetilde{\mathbb{Z}}}p\left[\mathfrak{z}_{p}\overset{p}{\equiv}-\nu\right]e^{-2\pi i\left\langle \frac{qx}{p}-t,\mathfrak{z}\right\rangle }d\mathfrak{z}\\
 & = & \mathbf{1}_{\mathbb{Z}\left[\frac{1}{p}\right]}\left(\frac{qx}{p}-t\right)\int_{\mathbb{Z}_{p}}p\left[\mathfrak{y}\overset{p}{\equiv}-\nu\right]e^{-2\pi i\left\{ \left(\frac{qx}{p}-t\right)\mathfrak{y}\right\} _{p}}d\mathfrak{y}\\
 & = & \mathbf{1}_{\frac{1}{p}\mathbb{Z}\left[\frac{1}{p}\right]}\left(qx-pt\right)\int_{p\mathbb{Z}_{p}-\nu}pe^{-2\pi i\left\{ \left(\frac{qx}{p}-t\right)\mathfrak{y}\right\} _{p}}d\mathfrak{y}\\
 & = & p\mathbf{1}_{\mathbb{Z}\left[\frac{1}{p}\right]}\left(qx-pt\right)\int_{p\mathbb{Z}_{p}}e^{-2\pi i\left\{ \left(\frac{qx}{p}-t\right)\left(\mathfrak{y}-\nu\right)\right\} _{p}}d\mathfrak{y}\\
 & = & pe^{2\nu\pi i\left\{ \frac{qx-pt}{p}\right\} _{p}}\mathbf{1}_{\mathbb{Z}\left[\frac{1}{p}\right]}\left(qx-pt\right)\int_{p\mathbb{Z}_{p}}e^{-2\pi i\left\{ \frac{qx-pt}{p}\mathfrak{y}\right\} _{p}}d\mathfrak{y}\\
\left(\mathfrak{x}=\frac{\mathfrak{y}}{p}\right); & = & pe^{2\nu\pi i\left\{ \frac{qx-pt}{p}\right\} _{p}}\mathbf{1}_{\mathbb{Z}\left[\frac{1}{p}\right]}\left(qx-pt\right)\frac{1}{\left|1/p\right|_{p}}\int_{\mathbb{Z}_{p}}e^{-2\pi i\left\{ \left(qx-pt\right)\mathfrak{x}\right\} _{p}}d\mathfrak{y}\\
 & = & pe^{2\nu\pi i\left\{ \frac{qx-pt}{p}\right\} _{p}}\mathbf{1}_{\mathbb{Z}\left[\frac{1}{p}\right]}\left(qx-pt\right)\frac{1}{p}\mathbf{1}_{\mathbb{Z}_{p}}\left(qx-pt\right)\\
 & = & e^{2\nu\pi i\left\{ \frac{qx-pt}{p}\right\} _{p}}\mathbf{1}_{\mathbb{Z}}\left(qx-pt\right)\\
\left(\textrm{eval. property}\right); & = & e^{2\nu\pi i\frac{qx-pt}{p}}\mathbf{1}_{\mathbb{Z}}\left(qx-pt\right)\\
 & = & \xi_{p}^{\nu\left(qx-pt\right)}\left[qx\overset{1}{\equiv}pt\right]
\end{eqnarray*}
Now:
\begin{eqnarray*}
qx & \overset{1}{\equiv} & pt\\
\left(\exists k\in\mathbb{Z}\right); & \Updownarrow\\
qx & = & pt+k\\
 & \Updownarrow\\
x & = & \frac{pt+k}{q}
\end{eqnarray*}
Thus:
\[
\left[qx\overset{1}{\equiv}pt\right]=\sum_{k=0}^{q-1}\left[x\overset{1}{\equiv}\frac{pt-k}{q}\right]=\sum_{k=0}^{q-1}\mathbf{1}_{\frac{pt+k}{q}}\left(x\right)
\]
which proves \ref{eq:Iverson bracket decomposition formula} and so:
\begin{eqnarray*}
\xi_{p}^{\nu\left(qx-pt\right)}\left[qx\overset{1}{\equiv}pt\right] & = & \sum_{k=0}^{q-1}\xi_{p}^{\nu\left(qx-pt\right)}\mathbf{1}_{\frac{pt-k}{q}}\left(x\right)\\
 & = & \sum_{k=0}^{q-1}\xi_{p}^{\nu\left(q\frac{pt-k}{q}-pt\right)}\mathbf{1}_{\frac{pt-k}{q}}\left(x\right)\\
 & = & \sum_{k=0}^{q-1}\xi_{p}^{-\nu k}\mathbf{1}_{\frac{pt-k}{q}}\left(x\right)
\end{eqnarray*}
which proves \ref{eq:Adic integral computation formula}.

Q.E.D.

\vphantom{}(II) Proof of \ref{eq: (rho,j) decomposition operator over Z-bar}:

\begin{eqnarray*}
\left[\mathfrak{z}_{\varrho}\overset{\varrho}{\equiv}j\right] & = & \sum_{\tau\in\left[0,1\right)_{\mathbb{Q}}}\left\langle \left[\mathfrak{z}_{\varrho}\overset{\varrho}{\equiv}j\right]\mid\check{\mathbf{1}}_{\tau}\right\rangle _{\widetilde{\mathbb{Z}}}\check{\mathbf{1}}_{\tau}\left(\mathfrak{z}\right)\\
 & = & \sum_{\tau\in\left[0,1\right)_{\mathbb{Q}}}\left(\int_{\widetilde{\mathbb{Z}}}\left[\mathfrak{z}_{\varrho}\overset{\varrho}{\equiv}j\right]e^{-2\pi i\left\langle \tau,\mathfrak{z}\right\rangle }d\mathfrak{z}\right)\check{\mathbf{1}}_{\tau}\left(\mathfrak{z}\right)\\
\left(\varrho\in\mathbb{P}\right) & = & \sum_{\tau\in\left[0,1\right)_{\mathbb{Q}}}\left(\mathbf{1}_{\mathbb{Z}\left[\frac{1}{\varrho}\right]}\left(\tau\right)\int_{\varrho\mathbb{Z}_{\varrho}}e^{-2\pi i\left\{ \tau\left(\mathfrak{y}+j\right)\right\} _{\varrho}}d\mathfrak{y}\right)\check{\mathbf{1}}_{\tau}\left(\mathfrak{z}\right)\\
\left(\mathfrak{x}=\frac{\mathfrak{y}}{\varrho}\right); & = & \sum_{\tau\in\left[0,1\right)_{\mathbb{Q}}}\left(\frac{1}{\varrho}\mathbf{1}_{\mathbb{Z}\left[\frac{1}{\varrho}\right]}\left(\tau\right)\int_{\mathbb{Z}_{\varrho}}e^{-2\pi i\left\{ \tau\left(\varrho\mathfrak{x}+j\right)\right\} _{\varrho}}d\mathfrak{x}\right)\check{\mathbf{1}}_{\tau}\left(\mathfrak{z}\right)\\
 & = & \sum_{\tau\in\left[0,1\right)_{\mathbb{Q}}}\left(\frac{1}{\varrho}\mathbf{1}_{\mathbb{Z}\left[\frac{1}{\varrho}\right]}\left(\tau\right)e^{-2\pi i\left\{ \tau j\right\} _{\varrho}}\mathbf{1}_{\mathbb{Z}_{\varrho}}\left(\varrho\tau\right)\right)\check{\mathbf{1}}_{\tau}\left(\mathfrak{z}\right)\\
 & = & \frac{1}{\varrho}\sum_{\tau\in\left[0,1\right)_{\mathbb{Q}}}e^{-2\pi i\left\{ \tau j\right\} _{\varrho}}\mathbf{1}_{\varrho\mathbb{Z}\left[\frac{1}{\varrho}\right]}\left(\varrho\tau\right)\mathbf{1}_{\mathbb{Z}_{\varrho}}\left(\varrho\tau\right)\check{\mathbf{1}}_{\tau}\left(\mathfrak{z}\right)\\
 & = & \frac{1}{\varrho}\sum_{\tau\in\left[0,1\right)_{\mathbb{Q}}}e^{-2\pi i\left\{ \tau j\right\} _{\varrho}}\mathbf{1}_{\mathbb{Z}\left[\frac{1}{\varrho}\right]}\left(\varrho\tau\right)\mathbf{1}_{\mathbb{Z}_{\varrho}}\left(\varrho\tau\right)\check{\mathbf{1}}_{\tau}\left(\mathfrak{z}\right)\\
\left(\mathbf{1}_{\mathbb{Z}\left[\frac{1}{\varrho}\right]}\times\mathbf{1}_{\mathbb{Z}_{\varrho}}=\mathbf{1}_{\mathbb{Z}}\right); & = & \frac{1}{\varrho}\sum_{\tau\in\left[0,1\right)_{\mathbb{Q}}}e^{-2\pi i\left\{ \tau j\right\} _{\varrho}}\mathbf{1}_{\mathbb{Z}}\left(\varrho\tau\right)\check{\mathbf{1}}_{\tau}\left(\mathfrak{z}\right)\\
\left(\mathbf{1}_{\mathbb{Z}}\left(\varrho\tau\right)=\mathbf{1}_{\frac{\mathbb{Z}}{\varrho}}\left(\tau\right)\right); & = & \frac{1}{\varrho}\sum_{\tau\in\left[0,1\right)_{\mathbb{Q}}}e^{-2\pi i\left\{ \tau j\right\} _{\varrho}}\mathbf{1}_{\frac{\mathbb{Z}}{\varrho}}\left(\tau\right)\check{\mathbf{1}}_{\tau}\left(\mathfrak{z}\right)\\
\left(\tau\in\frac{\mathbb{Z}}{\varrho}\Rightarrow e^{2\pi i\left\{ \tau j\right\} _{\varrho}}=e^{2\pi i\tau j}\right); & = & \frac{1}{\varrho}\sum_{\tau\in\left[0,1\right)_{\mathbb{Q}}}e^{-2\pi i\tau j}\mathbf{1}_{\frac{\mathbb{Z}}{\varrho}}\left(\tau\right)\check{\mathbf{1}}_{\tau}\left(\mathfrak{z}\right)\\
\left(\frac{\mathbb{Z}}{\varrho}\cap\left[0,1\right)=\left\{ 0,\frac{1}{\varrho},\ldots,\frac{\varrho-1}{\varrho}\right\} \right); & = & \frac{1}{\varrho}\sum_{k=0}^{\varrho-1}e^{-\frac{2\pi ikj}{\varrho}}\check{\mathbf{1}}_{\frac{k}{\varrho}}\left(\mathfrak{z}\right)\\
 & = & \frac{1}{\varrho}\sum_{k=0}^{\varrho-1}\xi_{\varrho}^{-jk}\check{\mathbf{1}}_{\frac{k}{\varrho}}\left(\mathfrak{z}\right)
\end{eqnarray*}

Q.E.D.

\subsubsection*{Proof of BF-$\mathbb{Q}/\mathbb{Z}$}

Let $t\in\left[0,1\right)_{\mathbb{Q}}$ be arbitrary. Then:
\begin{eqnarray*}
\mathscr{\mathfrak{Q}}_{H}\left\{ \mathbf{1}_{t}\right\} \left(x\right) & = & \frac{1}{\varrho}\sum_{j=0}^{\varrho-1}\sum_{k=0}^{\mu_{j}-1}\xi_{\mu_{j}}^{-kH\left(j\right)}e^{-2\pi i\frac{b_{j}}{a_{j}}x}\mathbf{1}_{t}\left(\frac{\varrho x+k}{\mu_{j}}\right)\\
 & = & \frac{1}{\varrho}\sum_{j=0}^{\varrho-1}\sum_{k=0}^{\mu_{j}-1}\xi_{\mu_{j}}^{-kH\left(j\right)}e^{-2\pi i\frac{b_{j}}{a_{j}}x}\int_{\widetilde{\mathbb{Z}}}\check{\mathbf{1}}_{t}\left(\mathfrak{z}\right)e^{-2\pi i\left\langle \frac{\varrho x+k}{\mu_{j}},\mathfrak{z}\right\rangle }d\mathfrak{z}\\
 & = & \frac{1}{\varrho}\sum_{j=0}^{\varrho-1}\sum_{k=0}^{\mu_{j}-1}\xi_{\mu_{j}}^{-kH\left(j\right)}e^{-2\pi i\left\langle x,\frac{b_{j}}{a_{j}}\right\rangle }\int_{\widetilde{\mathbb{Z}}}\check{\mathbf{1}}_{t}\left(\mathfrak{z}\right)e^{-2\pi i\left(\left\langle \frac{k}{\mu_{j}},\mathfrak{z}\right\rangle +\left\langle x,\frac{\varrho\mathfrak{z}}{\mu_{j}}\right\rangle \right)}d\mathfrak{z}\\
 & = & \frac{1}{\varrho}\sum_{j=0}^{\varrho-1}\int_{\widetilde{\mathbb{Z}}}\check{\mathbf{1}}_{t}\left(\mathfrak{z}\right)\sum_{k=0}^{\mu_{j}-1}\xi_{\mu_{j}}^{-kH\left(j\right)}\xi_{\mu_{j}}^{-k\left[\mathfrak{z}_{\mu_{j}}\right]_{\mu_{j}}}e^{-2\pi i\left\langle x,\frac{\varrho\mathfrak{z}}{\mu_{j}}+\frac{b_{j}}{a_{j}}\right\rangle }d\mathfrak{z}\\
 & = & \sum_{j=0}^{\varrho-1}\frac{\mu_{j}}{\varrho}\int_{\widetilde{\mathbb{Z}}}\left[\mathfrak{z}_{\mu_{j}}\overset{\mu_{j}}{\equiv}-H\left(j\right)\right]\check{\mathbf{1}}_{t}\left(\mathfrak{z}\right)e^{-2\pi i\left\langle x,\frac{d_{j}\mathfrak{z}+b_{j}}{a_{j}}\right\rangle }d\mathfrak{z}\\
 & = & \sum_{j=0}^{\varrho-1}\frac{\mu_{j}}{\varrho}\int_{\mu_{j}\widetilde{\mathbb{Z}}-H\left(j\right)}\check{\mathbf{1}}_{t}\left(\mathfrak{z}\right)e^{-2\pi i\left\langle x,\frac{d_{j}\mathfrak{z}+b_{j}}{a_{j}}\right\rangle }d\mathfrak{z}\\
\left(\mathfrak{y}=\frac{\mathfrak{z}+H\left(j\right)}{\mu_{j}}\right); & = & \sum_{j=0}^{\varrho-1}\frac{1}{\varrho}\int_{\widetilde{\mathbb{Z}}}\check{\mathbf{1}}_{t}\left(\mu_{j}\mathfrak{y}-H\left(j\right)\right)e^{-2\pi i\left\langle x,\frac{d_{j}\left(\mu_{j}\mathfrak{y}-H\left(j\right)\right)+b_{j}}{a_{j}}\right\rangle }d\mathfrak{y}\\
 & = & \frac{1}{\varrho}\sum_{j=0}^{\varrho-1}\int_{\widetilde{\mathbb{Z}}}\check{\mathbf{1}}_{t}\left(\mu_{j}\mathfrak{y}-H\left(j\right)\right)e^{-2\pi i\left\langle x,\frac{\varrho a_{j}\mathfrak{y}-ja_{j}-b_{j}+b_{j}}{a_{j}}\right\rangle }d\mathfrak{y}\\
 & = & \frac{1}{\varrho}\sum_{j=0}^{\varrho-1}\int_{\widetilde{\mathbb{Z}}}\check{\mathbf{1}}_{t}\left(\mu_{j}\mathfrak{y}-H\left(j\right)\right)e^{-2\pi i\left\langle x,\varrho\mathfrak{y}-j\right\rangle }d\mathfrak{y}\\
 & = & \frac{1}{\varrho}\sum_{j=0}^{\varrho-1}e^{2\pi ijx}e^{-2\pi itH\left(j\right)}\int_{\widetilde{\mathbb{Z}}}\check{\mathbf{1}}_{\mu_{j}t}\left(\mathfrak{y}\right)e^{-2\pi i\left\langle \varrho x,\mathfrak{y}\right\rangle }d\mathfrak{y}\\
 & = & \frac{1}{\varrho}\sum_{j=0}^{\varrho-1}e^{2\pi ijx}e^{-2\pi itH\left(j\right)}\mathscr{F}_{\mathbb{Q}/\mathbb{Z}}^{-1}\left\{ \check{\mathbf{1}}_{\mu_{j}t}\right\} \left(\varrho x\right)\\
 & = & \frac{1}{\varrho}\sum_{j=0}^{\varrho-1}e^{-2\pi itH\left(j\right)}e^{2\pi ijx}\mathbf{1}_{\mu_{j}t}\left(\varrho x\right)\\
 & = & \frac{1}{\varrho}\sum_{j=0}^{\varrho-1}e^{-2\pi itH\left(j\right)}e^{2\pi ijx}\sum_{k=0}^{\varrho-1}\mathbf{1}_{\frac{\mu_{j}t+k}{\varrho}}\left(x\right)\\
\left(\textrm{eval. prop.}\right); & = & \frac{1}{\varrho}\sum_{j=0}^{\varrho-1}e^{-2\pi it\overbrace{H\left(j\right)}^{\frac{ja_{j}+b_{j}}{d_{j}}}}\sum_{k=0}^{\varrho-1}e^{2\pi ij\frac{\mu_{j}t+k}{\varrho}}\mathbf{1}_{\frac{\mu_{j}t+k}{\varrho}}\left(x\right)\\
 & = & \frac{1}{\varrho}\sum_{j=0}^{\varrho-1}\sum_{k=0}^{\varrho-1}e^{2\pi i\left(-\frac{ja_{j}+b_{j}}{d_{j}}t+\frac{j\mu_{j}t}{\varrho}+\frac{jk}{\varrho}\right)}\mathbf{1}_{\frac{\mu_{j}t+k}{\varrho}}\left(x\right)\\
 & = & \frac{1}{\varrho}\sum_{j=0}^{\varrho-1}\sum_{k=0}^{\varrho-1}e^{2\pi i\left(\frac{jk}{\varrho}-\frac{b_{j}t}{d_{j}}\right)}\mathbf{1}_{\frac{\mu_{j}t+k}{\varrho}}\left(x\right)\\
 & = & \frac{1}{\varrho}\sum_{j=0}^{\varrho-1}e^{-2\pi i\frac{b_{j}}{d_{j}}t}\sum_{k=0}^{\varrho-1}\xi_{\varrho}^{jk}\mathbf{1}_{\frac{\mu_{j}t+k}{\varrho}}\left(x\right)
\end{eqnarray*}
which is the desired formula.

Q.E.D.

\subsubsection*{Proof of BF-$\widetilde{\mathbb{Z}}$}

Let $t\in\mathbb{Q}$. Beginning with:

\[
\mathfrak{Q}_{H}\left\{ \mathbf{1}_{t}\right\} \left(x\right)=\frac{1}{\varrho}\sum_{j=0}^{\varrho-1}e^{-2\pi i\frac{b_{j}}{d_{j}}t}\sum_{k=0}^{\varrho-1}\xi_{\varrho}^{jk}\mathbf{1}_{\frac{\mu_{j}t+k}{\varrho}}\left(x\right)
\]
and, recall, since: 
\[
\check{\mathbf{1}}_{t}\left(\mathfrak{z}\right)\overset{\textrm{def}}{=}e^{2\pi i\left\langle t,\mathfrak{z}\right\rangle }=\mathscr{F}_{\mathbb{Q}/\mathbb{Z}}\left\{ \mathbf{1}_{t}\right\} \left(\mathfrak{z}\right)
\]
the computation is \emph{almost} trivial: 
\begin{eqnarray*}
\check{\mathfrak{Q}}_{H}\left\{ \check{\mathbf{1}}_{t}\right\} \left(\mathfrak{z}\right) & = & \frac{1}{\varrho}\sum_{j=0}^{\varrho-1}e^{-\frac{2\pi ib_{j}t}{d_{j}}}\sum_{k=0}^{\varrho-1}\xi_{\varrho}^{jk}\check{\mathbf{1}}_{\frac{\mu_{j}t+k}{\varrho}}\left(\mathfrak{z}\right)\\
 & = & \frac{1}{\varrho}\sum_{j=0}^{\varrho-1}e^{-\frac{2\pi ib_{j}t}{d_{j}}}\sum_{k=0}^{\varrho-1}\xi_{\varrho}^{jk}\check{\mathbf{1}}_{\frac{k}{\varrho}}\left(\mathfrak{z}\right)\check{\mathbf{1}}_{\frac{\mu_{j}t}{\varrho}}\left(\mathfrak{z}\right)\\
 & = & \sum_{j=0}^{\varrho-1}e^{2\pi i\left\langle t,-\frac{b_{j}}{d_{j}}\right\rangle }\check{\mathbf{1}}_{t}\left(\frac{\mu_{j}\mathfrak{z}}{\varrho}\right)\frac{1}{\varrho}\sum_{k=0}^{\varrho-1}\xi_{\varrho}^{k\left(j+\left[\mathfrak{z}_{\varrho}\right]_{\varrho}\right)}\\
 & = & \sum_{j=0}^{\varrho-1}\check{\mathbf{1}}_{t}\left(-\frac{b_{j}}{d_{j}}\right)\check{\mathbf{1}}_{t}\left(\frac{\mu_{j}\mathfrak{z}}{\varrho}\right)\left[\mathfrak{z}_{\varrho}\overset{\varrho}{\equiv}-j\right]\\
 & = & \sum_{j=0}^{\varrho-1}\left[\mathfrak{z}_{\varrho}\overset{\varrho}{\equiv}-j\right]\check{\mathbf{1}}_{t}\left(\frac{\mu_{j}\mathfrak{z}}{\varrho}-\frac{b_{j}}{d_{j}}\right)\\
 & = & \sum_{j=0}^{\varrho-1}\left[\mathfrak{z}_{\varrho}\overset{\varrho}{\equiv}-j\right]\check{\mathbf{1}}_{t}\left(\frac{a_{j}\mathfrak{z}-b_{j}}{d_{j}}\right)
\end{eqnarray*}
which proves \ref{eq:Basis formula for Q-check; prime rho}.

Q.E.D.

\subsubsection*{Proof of GF-$\widetilde{\mathbb{Z}}$}

Since $\left\{ \mathbf{1}_{t}\left(x\right):t\in\left[0,1\right)_{\mathbb{Q}}\right\} $
is an orthonormal basis of $L^{2}\left(\mathbb{Q}/\mathbb{Z}\right)$,
the isometric isomorphism property of $\mathscr{F}_{\mathbb{Q}/\mathbb{Z}}$
shows that $\left\{ \check{\mathbf{1}}_{t}\left(\mathfrak{z}\right):t\in\left[0,1\right)_{\mathbb{Q}}\right\} $
is then an orthonormal basis for $L^{2}\left(\widetilde{\mathbb{Z}}\right)$.
As such, letting $F\in L^{2}\left(\widetilde{\mathbb{Z}}\right)$
be arbitrary, there are unique complex constants $\left\{ c_{t}\right\} _{t\in\left[0,1\right)_{\mathbb{Q}}}$
(the Fourier coefficients of $F$) so that:
\[
F\left(\mathfrak{z}\right)=\sum_{t\in\left[0,1\right)_{\mathbb{Q}}}c_{t}\check{\mathbf{1}}_{t}\left(\mathfrak{z}\right)
\]
Thus, using \ref{eq:Basis formula for Q-check; prime rho}, the linearity
of $\check{\mathfrak{Q}}_{H}$ gives:
\begin{eqnarray*}
\check{\mathfrak{Q}}_{H}\left\{ F\right\} \left(\mathfrak{z}\right) & = & \sum_{t\in\left[0,1\right)_{\mathbb{Q}}}c_{t}\check{\mathfrak{Q}}_{H}\left\{ \check{\mathbf{1}}_{t}\right\} \left(\mathfrak{z}\right)\\
\left(\textrm{use }\ref{eq:Basis formula for Q-check; prime rho}\right); & = & \sum_{t\in\left[0,1\right)_{\mathbb{Q}}}c_{t}\sum_{j=0}^{\varrho-1}\left[\mathfrak{z}_{\varrho}\overset{\varrho}{\equiv}-j\right]\check{\mathbf{1}}_{t}\left(\frac{a_{j}\mathfrak{z}-b_{j}}{d_{j}}\right)\\
 & = & \sum_{j=0}^{\varrho-1}\left[\mathfrak{z}_{\varrho}\overset{\varrho}{\equiv}-j\right]\sum_{t\in\left[0,1\right)_{\mathbb{Q}}}c_{t}\check{\mathbf{1}}_{t}\left(\frac{a_{j}\mathfrak{z}-b_{j}}{d_{j}}\right)\\
 & = & \sum_{j=0}^{\varrho-1}\left[\mathfrak{z}_{\varrho}\overset{\varrho}{\equiv}-j\right]F\left(\frac{a_{j}\mathfrak{z}-b_{j}}{d_{j}}\right)
\end{eqnarray*}
which is the desired formula.

Q.E.D.

\pagebreak{}

\section{Rationality Theorems, \&c.}

While we do not have, at present, any \textbf{Strong Rationality Theorem
}(\textbf{SRT}) results, the \textbf{Finite Dreamcatcher Theorem }(\textbf{FDT})
does immediately imply a \textbf{Weak Rationality Theorem }(\textbf{WRT})
for regulated prime\footnote{This result almost certainly extends to all regulated $\varrho$-hydra
maps, but we have refrained from pursuing that generalization for
the sake of keeping the already-intricate computations from growing
out of hand.} $\varrho$-hydra maps. The purpose of this brief chapter is to record
this and related applications of the \textbf{FDT }to set-series and
their represented sets.

\subsection{Equivalence modulo zero-density sets}

One extremely important feature of dreamcatcher analysis we have yet
to discuss is their inability to detect changes in a set-series modulo
a set of density zero. Specifically, given any $V\subseteq\mathbb{N}_{0}$,
$\psi_{V}$'s virtual poles remain unchanged when $V$ is altered
by a set of zero density.

\vphantom{}\textbf{Definition} \textbf{25}:

I. Let $\mathcal{M}_{0}$ denote the set of all subsets of $\mathbb{N}_{0}$
with a well-defined natural density.

II. Let $\Delta$ denote the symmetric difference operator:
\[
V\Delta W\overset{\textrm{def}}{=}\left(V\backslash W\right)\cup\left(W\backslash V\right)
\]
Then, define the relation $\approx$ between sets $V,W\in\mathcal{M}_{0}$
by:
\[
V\approx W\Leftrightarrow d\left(V\Delta W\right)=0
\]

\vphantom{}\emph{Remark}:\emph{ }The properties of natural density
guarantee that $\mathcal{M}_{0}$ is a sigma-algebra. Thus, $d$ defines
a measure on $\mathbb{N}_{0}$, where $\mathcal{M}_{0}$ is the sigma-algebra
of $d$-measurable sets. Since $\left(\mathbb{N}_{0},\mathcal{M}_{0},d\right)$
is a measure space, it then follows that the relation $\approx$ is
an equivalence relation on $\left(\mathbb{N}_{0},\mathcal{M}_{0},d\right)$.
In fact, $V\approx W$ if and only if $V$ and $W$ are equal ``almost
everywhere'' (i.e., with respect to $d$) as measurable sets in $\left(\mathbb{N}_{0},\mathcal{M}_{0},d\right)$.

\vphantom{}\textbf{Lemma 6 }(\textbf{Equivalence mod null sets)}:
Let $V,W\in\mathcal{M}_{0}$ have positive natural density. Then,
$V\approx W$ implies $\textrm{P}\left(\psi_{V}\right)=\textrm{P}\left(\psi_{W}\right)$
($V$ and $W$ have the same edgepoints), and $R_{V}\left(t\right)=R_{W}\left(t\right)$
for all $t\in\textrm{P}\left(\psi_{V}\right)\cup\textrm{P}\left(\psi_{W}\right)$
($V$ and $W$ have the same dreamcatchers).

Proof: If $V=W$, then there is nothing to prove. So, suppose $V\neq W$.
Like our work in Chapter 1, it suffices to prove the Lemma using the
ordinary set-series of $V$ and $W$.

I. Suppose one of $V$ or $W$ contains the other as a proper subset;
without loss of generality, suppose $W\subset V$. Then:
\[
\varsigma_{V}\left(z\right)-\varsigma_{V}\left(z\right)=\varsigma_{V\backslash W}\left(z\right)
\]
Since $W\subset V$, note that $V\Delta W=\left(V\backslash W\right)\cup\left(W\backslash V\right)=V\backslash W$.
As such, the condition that $V\approx W$ shows that $0=d\left(V\Delta W\right)=d\left(V\backslash W\right)$.
Then, for any $t\in\mathbb{R}$, note that, by digitality:
\begin{eqnarray*}
\left|\limsup_{r\uparrow1}\left(1-r\right)\varsigma_{V\backslash W}\left(re^{2\pi it}\right)\right| & \leq & \limsup_{r\uparrow1}\left(1-r\right)\varsigma_{V\backslash W}\left(\left|re^{2\pi it}\right|\right)\\
 & = & \limsup_{r\uparrow1}\left(1-r\right)\varsigma_{V\backslash W}\left(r\right)\\
 & = & \overline{d}\left(V\backslash W\right)\\
\left(d\left(V\backslash W\right)\textrm{ exists and is }0\right); & = & 0
\end{eqnarray*}
Thus:
\[
\lim_{r\uparrow1}\left(1-r\right)\left(\varsigma_{V}\left(re^{2\pi it}\right)-\varsigma_{W}\left(re^{2\pi it}\right)\right)=\lim_{r\uparrow1}\left(1-r\right)\varsigma_{V\backslash W}\left(re^{2\pi it}\right)=0,\textrm{ }\forall t\in\mathbb{R}
\]
As such, for any $t\in\mathbb{R}$, $\textrm{VRes}_{e^{2\pi it}}\left[\varsigma_{V}\right]=\lim_{r\uparrow1}\left(1-r\right)\varsigma_{V}\left(re^{2\pi it}\right)$
exists if and only if $\textrm{VRes}_{e^{2\pi it}}\left[\varsigma_{W}\right]=\lim_{r\uparrow1}\left(1-r\right)\varsigma_{W}\left(re^{2\pi it}\right)$
exists, and if they exist, they must be equal. Consequently, $\textrm{VRes}_{t}\left[\psi_{V}\right]$
exists if and only if $\textrm{VRes}_{t}\left[\psi_{W}\right]$ exists
and that $\textrm{VRes}_{t}\left[\psi_{V}\right]=\textrm{VRes}_{t}\left[\psi_{W}\right]$
for all $t$ for which either exists ($\textrm{P}\left(\psi_{V}\right)=\textrm{P}\left(\psi_{W}\right)$
and $R_{V}\left(t\right)=R_{W}\left(t\right)$ for all $t\in\textrm{P}\left(\psi_{V}\right)\cup\textrm{P}\left(\psi_{W}\right)$).
Thus, \textbf{Lemma 6 }holds whenever $W\subset V$ or $V\subset W$.
$\checkmark$

II. Suppose that neither $V$ nor $W$ fully contains the other.

If $V$ and $W$ are disjoint, then $V\Delta W=V\cup W$; but then,
since $V\approx W$, the sub-additivity of the measure $d$ allows
us to write:
\[
0\overset{V\approx W}{=}d\left(V\Delta W\right)=d\left(V\cup W\right)\overset{\textrm{sub-add}}{=}d\left(V\right)+d\left(W\right)\overset{V\approx W}{=}2d\left(V\right)
\]
which forces $d\left(V\right)=0$. However, this contradicts the fact
that $V$ was given to have a positive natural density. Thus, it must
be that $V$and $W$ are not disjoint: $V\cap W$.

Since $V\cap W$ is non-empty, and since neither $V$ nor $W$ fully
contain the other, $V\cap W$ is necessarily a proper subset of both
$V$ and $W$. Moreover:
\[
V\Delta\left(V\cap W\right)=\left(V\backslash\left(V\cap W\right)\right)\cup\underbrace{\left(\left(V\cap W\right)\backslash V\right)}_{\varnothing}=V\backslash W\subseteq\left(V\backslash W\right)\cup\left(W\backslash V\right)=V\Delta W
\]
and, by symmetry (swap the positions of $V$ and $W$ in the above
equalities) and the commutativity of $\Delta$:
\[
W\Delta\left(V\cap W\right)=W\Delta V=V\Delta W
\]
Thus $0=d\left(V\Delta W\right)=d\left(V\Delta\left(V\cap W\right)\right)$,
which implies $V\approx V\cap W$; likewise, $0=d\left(V\Delta W\right)=d\left(W\Delta\left(V\cap W\right)\right)$,
which implies $W\approx V\cap W$. Since $V\cap W\subset V$, Case
(I) then shows that the boundary data of $\varsigma_{V}\left(z\right)$
and $\varsigma_{V\cap W}\left(z\right)$ are the same; likewise, $V\cap W\subset W$
implies that the boundary data of $\varsigma_{W}\left(z\right)$ and
$\varsigma_{V\cap W}\left(z\right)$ are the same. Thus, as Euclid
(or Lincoln) would say, it is then self-evident that the boundary
data of $\varsigma_{V}\left(z\right)$ and $\varsigma_{W}\left(z\right)$
must be the same, as desired. $\checkmark$

Q.E.D.

\vphantom{}\textbf{Definition} \textbf{26}: Let $\widetilde{\mathcal{M}}_{0}$
denote the set of equivalence classes of $\mathcal{M}_{0}$ under
the relation $\approx$. It is known (Bell. \emph{The Symmetric Difference
Metric}) that\textbf{ }$d$ defines a complete metric on $\widetilde{\mathcal{M}}_{0}$.

\subsection{Consequences of the Dreamcatcher Theorem}

Throughout this section, let $H:\mathbb{N}_{0}\rightarrow\mathbb{N}_{0}$
be a regulated prime $\varrho$-hydra map. There are three results
to be given; the first two are \textbf{WRT}s.

\vphantom{}\textbf{Theorem 10} (\textbf{WRT 1}): Let $V\in\mathcal{M}_{0}$
be $H$-invariant. If $V$ is a rational set, then one of $V$ and
$\mathbb{N}_{0}/V$ is finite (or empty).

Proof: If $V$ is either empty or finite, then we have nothing to
prove. So, suppose $V$ is infinite. Since $V$ is a rational set,
its Fourier set-series is a rational function of $e^{2\pi iz}$. Consequently,
the virtual poles of $\psi_{V}$ and residues thereof correspond with
its classical poles and residues. As a rational function of $e^{2\pi iz}$,
the set of poles of $\psi_{V}$ is a discrete subset of $\mathbb{R}$,
and in particular, there are only finitely many poles in $\left[0,1\right)$.
Thus, the dreamcatcher $R_{V}$ of $\psi_{V}$ is finitely supported
in $\mathbb{Q}/\mathbb{Z}$, and so, by the \textbf{FDT}, $R_{V}$
has the form:
\[
R_{V}\left(x\right)=R_{V}\left(0\right)\mathbf{1}_{0}\left(x\right),\textrm{ }\forall x\in\mathbb{Q}
\]
Since $V$ is an infinite rational set, it necessarily has positive
density. Consequently, by the \textbf{HLTT}, $0<d\left(V\right)=R_{V}\left(0\right)$.

Now, writing $\psi_{V}$ out explicitly as a rational function of
$e^{2\pi iz}$, we have:
\[
\psi_{V}\left(z\right)=Q\left(e^{2\pi iz}\right)+\frac{e^{2\eta\pi iz}P\left(e^{2\pi iz}\right)}{1-e^{2\alpha\pi iz}}
\]
for digital polynomials $Q,P\in\mathbb{C}\left[z\right]$, and integers
$\eta\in\mathbb{N}_{0}$ and $\alpha\in\mathbb{N}_{1}$ (where $\alpha$
is the period of $V$), where $P\left(z\right)$ and $1-z^{\alpha}$
(and hence, $P\left(e^{2\pi iz}\right)$ and $1-e^{2\alpha\pi iz}$)
are co-prime as polynomials, meaning that they share no common zeroes.
Because of this co-primality, it follows that all the roots of unity
$z\in\left\{ 1,\xi_{\alpha},\ldots,\xi_{\alpha}^{\alpha-1}\right\} $
are poles of $\frac{P\left(z\right)}{1-z^{\alpha}}$, and thus, that
$\psi_{V}\left(z\right)$ has poles at every $z\in\frac{1}{\alpha}\mathbb{Z}$.

Now, \emph{suppose $\alpha\geq2$}. Since every pole $z=x_{0}\in\mathbb{Q}$
of $\psi_{V}\left(z\right)$ is necessarily also a \emph{virtual }pole
of $\psi_{V}\left(z\right)$, and since every virtual pole $z=x_{0}\in\mathbb{Q}$
of $\psi_{V}\left(z\right)$ lies in the support of $\psi_{V}$'s
dreamcatcher, $R_{V}$, it must be that $R_{V}\left(x\right)\neq0$
has poles for all $x\in\frac{1}{\alpha}\mathbb{Z}$. Since $\alpha\geq2$,
$\frac{1}{\alpha}$ is then a non-integer rational number. But then,
for the virtual pole $x=\frac{1}{\alpha}$:
\[
0\neq R_{V}\left(\frac{1}{\alpha}\right)=\underbrace{R_{V}\left(0\right)}_{>0}\overbrace{\mathbf{1}_{0}\left(\frac{1}{\alpha}\right)}^{0}=0\textrm{ }
\]
which is impossible! Thus, the assumption that $\alpha\geq2$ leads
to a contradiction, and so, $\alpha$ \emph{must be equal to $1$}.

Since $\alpha=1$ is the period of $V$, this shows that $V$ contains
all sufficiently large integers. Hence, $\mathbb{N}_{0}\backslash V$
is finite (or empty). Thus, if $V$ is neither finite nor empty, $\mathbb{N}_{0}/V$
must be either finite or empty. This proves the theorem.

Q.E.D

\vphantom{}\textbf{Theorem 11} (\textbf{WRT 2}):\textbf{ }Let $V\in\mathcal{M}_{0}$
be $H$-invariant. If there is a rational set $W$ so that $V\approx W$,
then $d\left(V\right)\in\left\{ 0,1\right\} $.

Proof: Since $V\approx W$, $\psi_{V}$ and $\psi_{W}$ have identical
dreamcatchers, both of which are meromorphic. Applying the argument
from \textbf{Theorem 10 }then shows that one of $W$ and $\mathbb{N}_{0}\backslash W$
is finite or empty. Since $V\approx W$, observe that if $W$ is finite
or empty, then $d\left(W\right)=0$ and so, $d\left(V\right)=0$.
On the other hand, if $\mathbb{N}_{0}\backslash W$ is finite or empty,
then $1-d\left(W\right)=d\left(\mathbb{N}_{0}\backslash W\right)=0$,
and so, $d\left(W\right)=1$, and so, $V\approx W$ implies $d\left(V\right)=d\left(W\right)=1$.
This proves the theorem.

Q.E.D.

\vphantom{}Finally, here is the most general consequence of the \textbf{FDT}:

\vphantom{}\textbf{Theorem 12}: Let $V\in\mathcal{M}_{0}$ be $H$-invariant
set. For any given integers $\alpha,\beta$ with $0\leq\beta\leq\alpha-1$,
let: 
\[
V_{\alpha,\beta}\overset{\textrm{def}}{=}\left(\alpha\mathbb{Z}+\beta\right)\cap V
\]
denote the set of all elements of $V$ congruent to $\beta$ mod $\alpha$.
Additionally, let: 
\[
V_{\alpha,\beta}\left(N\right)\overset{\textrm{def}}{=}\left\{ v\in V_{\alpha,\beta}:v\leq N\right\} 
\]
and let $\left|V_{\alpha,\beta}\left(N\right)\right|$ denote the
number of elements of $V_{\alpha,\beta}\left(N\right)$.

Now, suppose both $T_{V}\overset{\textrm{def}}{=}\mathbb{Q}\cap\textrm{P}\left(\psi_{V}\right)$
and $T_{V}^{c}\overset{\textrm{def}}{=}\mathbb{Q}\backslash\textrm{P}\left(\psi_{V}\right)$
are $H$-branch invariant and discrete. Letting $T_{V}^{\prime}\overset{\textrm{def}}{=}T_{V}/\mathbb{Z}=T_{V}\cap\left[0,1\right)_{\mathbb{Q}}$,
write elements $t\in T_{V}^{\prime}$ as $t=\frac{\beta}{\alpha}$,
where $0\leq\alpha-1\leq\beta$ and $\gcd\left(\alpha,\beta\right)=1$
whenever $\beta\neq0$. Then:
\begin{equation}
\lim_{N\rightarrow\infty}\frac{1}{N}\sum_{k=0}^{\alpha-1}\left|V_{\alpha,k}\left(N\right)\right|\xi_{\alpha}^{\beta k}=0,\textrm{ }\forall\frac{\beta}{\alpha}\in T_{V}^{\prime}\backslash\left\{ 0\right\} \label{eq:Dreamcatcher Application - Random Walk Limit}
\end{equation}

Proof: Let $V$ and $T_{V}^{\prime}$ be as given. By the\textbf{
Extended SCL}-$\mathbb{Z}$, the $T_{V}$-dreamcatcher $R_{V}$ of
$\psi_{V}$ extends to a function in $L^{2}\left(\mathbb{Q}/\mathbb{Z}\right)$
which is fixed by $\mathfrak{Q}_{H}$. Since $T_{V}$ is discrete,
$R_{V}$ is finitely supported on $T_{V}$, and hence, on $\mathbb{Q}/\mathbb{Z}$.
By the \textbf{Finite Dreamcatcher Theorem}, this extension must then
be of the form $R_{V}\left(0\right)\mathbf{1}_{\mathbb{Z}}\left(x\right)$
for all $x\in\mathbb{Q}$. Consequently, restricting back to $x\in T_{V}$:
\[
R_{V}\left(x\right)=R_{V}\left(0\right)\mathbf{1}_{\mathbb{Z}}\left(x\right)=\frac{d\left(V\right)}{2\pi}\left[x\in\mathbb{Z}\right],\textrm{ }\forall x\in T_{V}
\]
where, recall, the existence of $d\left(V\right)$ is guaranteed by
the fact that $V\in\mathcal{M}_{0}$.

Now, letting $t=\frac{\beta}{\alpha}\in T_{V}^{\prime}\backslash\left\{ 0\right\} $
be arbitrary:
\[
0=\frac{d\left(V\right)}{2\pi}\left[t\in\mathbb{Z}\right]=R_{V}\left(t\right)\overset{\textrm{def}}{=}\lim_{y\downarrow0}y\psi_{V}\left(t+iy\right)=\lim_{y\downarrow0}y\sum_{v\in V}e^{2\pi iv\left(t+iy\right)}
\]
Hence:
\[
0=\lim_{y\downarrow0}y\sum_{v\in V}e^{2\pi iv\left(t+iy\right)}=\frac{1}{2\pi}\lim_{N\rightarrow\infty}\frac{1}{N}\sum_{k=0}^{\alpha-1}\xi_{\alpha}^{k\beta}\left|V_{\alpha,k}\left(N\right)\right|
\]
and so:
\[
\lim_{N\rightarrow\infty}\frac{1}{N}\sum_{k=0}^{\alpha-1}\xi_{\alpha}^{k\beta}\left|V_{\alpha,k}\left(N\right)\right|=0,\textrm{ }\forall\frac{\beta}{\alpha}\in T_{V}^{\prime}\backslash\left\{ 0\right\} 
\]
as desired.

Q.E.D.

\pagebreak{}

\section{Dénouement}

\subsection{Technical Issues}

A particularly appealing way of formulating the \textbf{ESCL} (\textbf{Corollary
3})\textbf{ }is as follows:

\vphantom{}\textbf{Corollary 4: Extended SCL} (\textbf{ESCL}) - \textbf{Alternative
Formulation}: Let $H:\mathbb{N}_{0}\rightarrow\mathbb{N}_{0}$ be
a $\varrho$-gydra map, and let $V\subseteq\mathbb{N}_{0}$ be $H$-invariant.
Then, the implication:
\[
\psi_{V}\in\textrm{Ker}\left(1-\mathscr{Q}_{H}\right)\Rightarrow R_{V}\in\textrm{Ker}\left(1-\mathfrak{Q}_{H:\mathbb{Q}}\right)
\]
holds true whenever both $\mathbb{Q}\cap\textrm{P\ensuremath{\left(V\right)}}$
and $\mathbb{Q}\backslash\textrm{P\ensuremath{\left(V\right)}}$ are
$H$-branch invariant.

Proof: Let $T=\textrm{P}\left(V\right)\cap\mathbb{Q}$ and $S=\mathbb{Q}$.
If $\mathbb{Q}\cap\textrm{P}\left(V\right)$ and $\mathbb{Q}\backslash\textrm{P\ensuremath{\left(V\right)}}$
are both $H$-branch invariant, then the above implication follows
by the \textbf{ESCL}.

Q.E.D.

\vphantom{}The importance of this is that it shows that dreamcatcher
methods are applicable to a given $\psi_{V}$ only if the set of edgepoints
and non-edgepoints of $\psi_{V}$ in $\mathbb{Q}$ are invariant under
the inverses of the branches of $H$. These technical obstacles to
the applicability of dreamcatcher techniques seem to occur quite often,
and are the technique's Achilles'-heel. At present, the author has
been unable to wash away these difficulties. Nevertheless, the author
hopes that the lines of inquiry detailed in next few definitions and
results might prove useful in shedding light on the matter.

\vphantom{}\textbf{Definition} \textbf{27}: Let $H:\mathbb{N}_{0}\rightarrow\mathbb{N}_{0}$
be a $\varrho$-hydra map. Let: 
\[
\mathcal{T}_{H}\overset{\textrm{def}}{=}\left\{ T\subseteq\mathbb{Q}/\mathbb{Z}:T\textrm{ \& }\mathbb{Q}/T\textrm{ are both }H\textrm{-branch invariant}\right\} 
\]
and then let:
\[
\mathcal{D}_{H}\overset{\textrm{def}}{=}\bigcup_{T\in\mathcal{T}_{H}}\mathcal{D}_{2}\left(T\right)
\]
denote the set of all $\psi\in\mathcal{A}\left(\mathbb{H}_{+i}\right)$
which admit an $T$-dreamcatcher for some $T\in\mathcal{T}_{H}$.

\vphantom{}\textbf{Fact}: $\mathcal{D}_{H}$ contains every function
in $\mathbb{C}\left(e^{2\pi iz}\right)$---that is, every rational
function of $e^{2\pi iz}$ with complex coefficients.

\vphantom{}\textbf{Lemma} \textbf{7}: Let $\left\{ \psi_{m}\right\} _{m\in\mathbb{N}_{1}}$
be a sequence in $\mathcal{D}_{H}$. Let $T_{m}\overset{\textrm{def}}{=}\mathbb{Q}\cap\textrm{P}\left(\psi_{m}\right)$.
For each $m$, let $R_{m}\left(x\right)$ be the extension of $\textrm{VRes}_{x}\left[\psi_{m}\right]\mid_{x\in T_{m}}$
to $\mathbb{Q}$ by way of the formula: 
\[
R_{m}\left(x\right)\overset{\textrm{def}}{=}\begin{cases}
\textrm{VRes}_{x}\left[\psi_{m}\right] & \textrm{if }x\in T_{m}\\
0 & \textrm{if }x\in\mathbb{Q}\backslash T_{m}
\end{cases}
\]
 Then, writing $\psi_{m}\left(z\right)=\sum_{k=0}^{\infty}c_{m,k}e^{2k\pi iz}$,
it follows that:
\begin{equation}
\left\Vert R_{m}-R_{n}\right\Vert _{\mathbb{Q}/\mathbb{Z}}^{2}\leq\frac{1}{2\pi}\lim_{y\downarrow0}y\sum_{k=0}^{\infty}\left|c_{m,k}-c_{n,k}\right|^{2}e^{-2k\pi y}\label{eq:Hardy-space/L^2 norm inequality}
\end{equation}

Proof: For brevity, write:
\begin{equation}
\Delta_{m,n}\left(x\right)=\lim_{y\downarrow0}y\left(\psi_{m}\left(x+iy\right)-\psi_{n}\left(x+iy\right)\right)=R_{m}\left(x\right)-R_{n}\left(x\right)\label{eq:Def of =002206_m,n}
\end{equation}
Then, just like in the proof of the \textbf{Square-Sum Lemma}:
\begin{eqnarray*}
\left\Vert \Delta_{m,n}\right\Vert _{L^{2}}^{4} & = & \left(\sum_{x\in\left[0,1\right)_{\mathbb{Q}}}\Delta_{m,n}\left(x\right)\overline{\Delta_{m,n}\left(x\right)}\right)^{2}\\
 & = & \left(\sum_{x\in\left[0,1\right)_{\mathbb{Q}}}\lim_{y\downarrow0}\left(\sum_{k=0}^{\infty}y\left(c_{m,k}-c_{n,k}\right)e^{2k\pi i\left(x+iy\right)}\right)\overline{\Delta_{m,n}\left(x\right)}\right)^{2}\\
 & = & \left(\lim_{y\downarrow0}\sum_{k=0}^{\infty}\sqrt{y}\left(c_{m,k}-c_{n,k}\right)e^{-k\pi y}\sum_{x\in\left[0,1\right)_{\mathbb{Q}}}\sqrt{y}\overline{\Delta_{m,n}\left(x\right)}e^{2k\pi ix}e^{-k\pi y}\right)^{2}\\
\left(\textrm{CSI}\right); & \leq & \lim_{y\downarrow0}\sum_{k=0}^{\infty}\left|c_{m,k}-c_{n,k}\right|^{2}e^{-2k\pi y}\\
 &  & \times\lim_{y^{\prime}\downarrow0}\sum_{j=0}^{\infty}y^{\prime}\left|\sum_{x\in\left[0,1\right)_{\mathbb{Q}}}\overline{\Delta_{m,n}\left(x\right)}e^{2j\pi ix}\right|^{2}e^{-2j\pi y^{\prime}}
\end{eqnarray*}
As such, like with the proof of the square summability lemma:
\[
\lim_{y^{\prime}\downarrow0}\sum_{j=0}^{\infty}y^{\prime}\left|\sum_{x\in\left[0,1\right)_{\mathbb{Q}}}\overline{\Delta_{m,n}\left(x\right)}e^{2j\pi ix}\right|^{2}e^{-2j\pi y^{\prime}}=\frac{\left\Vert \Delta_{m,n}\right\Vert _{L^{2}}^{2}}{2\pi}
\]
and so:
\[
\left\Vert \Delta_{m,n}\right\Vert _{L^{2}}^{4}\leq\frac{\left\Vert \Delta_{m,n}\right\Vert _{L^{2}}^{2}}{2\pi}\lim_{y\downarrow0}y\sum_{k=0}^{\infty}\left|c_{m,k}-c_{n,k}\right|^{2}e^{-2k\pi y}
\]
and so, either $\left\Vert R_{m}-R_{n}\right\Vert ^{2}=\left\Vert \Delta_{m,n}\right\Vert _{L^{2}}^{2}=0$
or:
\[
\left\Vert R_{m}-R_{n}\right\Vert ^{2}\leq\frac{1}{2\pi}\lim_{y\downarrow0}y\sum_{k=0}^{\infty}\left|c_{m,k}-c_{n,k}\right|^{2}e^{-2k\pi y}
\]
If $\left\Vert R_{m}-R_{n}\right\Vert ^{2}=0$, then the above inequality
is automatically true, since the upper bound (even if divergent to
$\infty$) is necessarily $\geq0$.

Q.E.D.

\vphantom{}\textbf{Definition} \textbf{28}: For every $\psi\in\mathcal{A}\left(\mathbb{H}_{+i}\right)$,
define:
\begin{equation}
\left\Vert \psi\right\Vert _{2,\frac{1}{y}}\overset{\textrm{def}}{=}\sqrt{\lim_{y\downarrow0}y\int_{0}^{1}\left|\psi\left(x+iy\right)\right|^{2}dx}\label{eq:(1/y)-asymptotic 2-hardy space norm}
\end{equation}

\emph{Remarks}:

I. To see where this formula comes from, recall that for any $\psi\left(z\right)=\sum_{n=0}^{\infty}c_{n}e^{2n\pi iz}$,
the \textbf{Parseval--Gutzmer formula }holds:
\[
\int_{0}^{1}\left|\psi\left(x+iy\right)\right|^{2}dx=\sum_{n=0}^{\infty}\left|c_{n}\right|^{2}e^{-2n\pi y}
\]
Hence:
\[
\lim_{y\downarrow0}y\int_{0}^{1}\left|\sum_{n=0}^{\infty}a_{n}e^{2n\pi i\left(x+iy\right)}-\sum_{n=0}^{\infty}b_{n}e^{2n\pi i\left(x+iy\right)}\right|^{2}dx=\lim_{y\downarrow0}y\sum_{n=0}^{\infty}\left|a_{n}-b_{n}\right|^{2}e^{-2n\pi y}
\]

II. \textbf{Lemma 7 }can then be stated as:

\begin{equation}
\left\Vert R_{m}-R_{n}\right\Vert _{\mathbb{Q}/\mathbb{Z}}\leq\frac{1}{\sqrt{2\pi}}\left\Vert \psi_{m}-\psi_{n}\right\Vert _{2,\frac{1}{y}}\label{eq:Hardy-space/L^2 norm inequality, simple version}
\end{equation}

\vphantom{}\textbf{Proposition} \textbf{14}: $\left\Vert \cdot\right\Vert _{2,\frac{1}{y}}$
is a semi-norm on $\mathcal{A}\left(\mathbb{H}_{+i}\right)$ which
is \emph{not }a norm on $\mathcal{A}\left(\mathbb{H}_{+i}\right)$.
In particular, for any set $V\subseteq\mathbb{N}_{0}$ with $d\left(V\right)=0$,
the Fourier set-series $\psi_{V}\left(z\right)$ satisfies:
\[
\left\Vert \psi_{V}\right\Vert _{2,\frac{1}{y}}=0
\]

Proof:

I. (Absolute homogeneity): Let $\alpha\in\mathbb{C}$. Then:
\[
\left\Vert \alpha\psi\right\Vert _{2,\frac{1}{y}}\overset{\textrm{def}}{=}\sqrt{\lim_{y\downarrow0}y\int_{0}^{1}\left|\alpha\psi\left(x+iy\right)\right|^{2}dx}=\left|\alpha\right|\sqrt{\lim_{y\downarrow0}y\int_{0}^{1}\left|\psi\left(x+iy\right)\right|^{2}dx}=\left|\alpha\right|\left\Vert \psi\right\Vert _{2,\frac{1}{y}}
\]
and so, $\left\Vert \cdot\right\Vert _{2,\frac{1}{y}}$ is absolutely
homogeneous.

II. (Triangle Inequality): Let $f,g\in\mathcal{A}\left(\mathbb{H}_{+i}\right)$.
Then:

\begin{eqnarray*}
\left\Vert f-g\right\Vert _{2,\frac{1}{y}} & = & \sqrt{\lim_{y\downarrow0}y\int_{0}^{1}\left|f\left(x+iy\right)-g\left(x+iy\right)\right|^{2}dx}\\
\left(\sqrt{t}\textrm{ is continuous }\forall t\in\mathbb{R}\geq0\right); & = & \lim_{y\downarrow0}\sqrt{y}\sqrt{\int_{0}^{1}\left|f\left(x+iy\right)-g\left(x+iy\right)\right|^{2}dx}\\
 & = & \lim_{y\downarrow0}\sqrt{y}\left\Vert f\left(x+iy\right)-g\left(x+iy\right)\right\Vert _{L^{2}:x\in\left[0,1\right]}\\
\left(L^{2}\textrm{ norms satisfy }\Delta\textrm{-ineq.}\right); & \leq & \lim_{y\downarrow0}\sqrt{y}\left\Vert f\left(x+iy\right)\right\Vert _{L^{2}:x\in\left[0,1\right]}\\
 &  & +\lim_{y\downarrow0}\sqrt{y}\left\Vert g\left(x+iy\right)\right\Vert _{L^{2}:x\in\left[0,1\right]}\\
\left(\sqrt{t}\textrm{ is continuous }\forall t\in\mathbb{R}\geq0\right); & = & \sqrt{\lim_{y\downarrow0}y\int_{0}^{1}\left|f\left(x+iy\right)\right|^{2}dx}\\
 &  & +\sqrt{\lim_{y\downarrow0}y\int_{0}^{1}\left|g\left(x+iy\right)\right|^{2}dx}\\
 & = & \left\Vert f\right\Vert _{2,\frac{1}{y}}+\left\Vert g\right\Vert _{2,\frac{1}{y}}
\end{eqnarray*}

III. Let $V\subseteq\mathbb{N}_{0}$ be a set of zero natural density.
Then, for every $y>0$, the \textbf{Parseval--Gutzmer formula }(\textbf{PGF})\textbf{
}yields:
\begin{eqnarray*}
\int_{0}^{1}\left|\psi_{V}\left(x+iy\right)\right|^{2}dx & = & \int_{0}^{1}\left|\sum_{n=0}^{\infty}\mathbf{1}_{V}\left(n\right)e^{2\pi in\left(x+iy\right)}\right|^{2}dx\\
\left(\textrm{\textbf{PGF}}\right); & = & \sum_{n=0}^{\infty}\left|\mathbf{1}_{V}\left(n\right)\right|^{2}e^{-2n\pi y}\\
 & = & \sum_{n=0}^{\infty}\mathbf{1}_{V}\left(n\right)e^{-2n\pi y}\\
 & = & \psi_{V}\left(iy\right)\\
 & = & \sigma_{V}\left(2\pi y\right)
\end{eqnarray*}
where, recall, $\sigma_{V}$ is the exponential set-series for $V$.
As such, by the \textbf{Hardy-Littlewood Tauberian Theorem} (\textbf{HLTT}),
since $V$ was given to have zero density:
\begin{eqnarray*}
\left\Vert \psi_{V}\right\Vert _{2,\frac{1}{y}} & = & \sqrt{\lim_{y\downarrow0}y\int_{0}^{1}\left|\psi_{V}\left(x+iy\right)\right|^{2}dx}\\
 & = & \sqrt{\lim_{y\downarrow0}y\sigma_{V}\left(2\pi y\right)}\\
 & = & \sqrt{\frac{1}{2\pi}\lim_{y\downarrow0}2\pi y\sigma_{V}\left(2\pi y\right)}\\
 & = & \sqrt{\frac{1}{2\pi}\lim_{x\downarrow0}x\sigma_{V}\left(x\right)}\\
\left(\textrm{\textbf{HLTT}}\right); & = & \sqrt{\frac{d\left(V\right)}{2\pi}}\\
\left(d\left(V\right)=0\right); & = & 0
\end{eqnarray*}
and so, $d\left(V\right)=0$ implies $\left\Vert \psi_{V}\right\Vert _{2,\frac{1}{y}}=0$.

Q.E.D.

\vphantom{}\textbf{Definition 29}:

I. Let: 
\[
\textrm{Dom}\left(\left\Vert \cdot\right\Vert _{2,\frac{1}{y}}\right)\overset{\textrm{def}}{=}\left\{ \psi\in\mathcal{A}\left(\mathbb{H}_{+i}\right):\left\Vert \psi\right\Vert _{2,\frac{1}{y}}<\infty\right\} 
\]

\emph{Remark}: Note that, since $\left\Vert \cdot\right\Vert _{2,\frac{1}{y}}$
is a semi-norm, $\textrm{Dom}\left(\left\Vert \cdot\right\Vert _{2,\frac{1}{y}}\right)$
is then a $\mathbb{C}$-linear space.

II. Writing: 
\[
\textrm{Ker}\left(\left\Vert \cdot\right\Vert _{2,\frac{1}{y}}\right)\overset{\textrm{def}}{=}\left\{ \psi\in\textrm{Dom}\left(\left\Vert \cdot\right\Vert _{2,\frac{1}{y}}\right):\left\Vert \psi\right\Vert _{2,\frac{1}{y}}=0\right\} 
\]
define the \textbf{$\frac{1}{y}$-Semi-Hardy space }$\mathcal{H}_{1/y}^{2}\left(\mathbb{H}_{+i}\right)$
as the normed $\mathbb{C}$-linear space obtained by equipping the
quotient $\mathbb{C}$-linear space $\textrm{Dom}\left(\left\Vert \cdot\right\Vert _{2,\frac{1}{y}}\right)/\textrm{Ker}\left(\left\Vert \cdot\right\Vert _{2,\frac{1}{y}}\right)$
with the $\left\Vert \cdot\right\Vert _{\frac{1}{2},y}$ semi-norm:
\[
\mathcal{H}_{1/y}^{2}\left(\mathbb{H}_{+i}\right)\overset{\textrm{def}}{=}\left(\textrm{Dom}\left(\left\Vert \cdot\right\Vert _{2,\frac{1}{y}}\right)/\textrm{Ker}\left(\left\Vert \cdot\right\Vert _{2,\frac{1}{y}}\right),\left\Vert \cdot\right\Vert _{2,\frac{1}{y}}\right)
\]

Let $\mathcal{H}_{1/y}^{2}\left(\mathbb{H}_{+i}\right)$ denote the
linear subspace of $\mathcal{A}\left(\mathbb{H}_{+i}\right)$ of functions
whose \textbf{$\mathcal{H}_{1/y}^{2}\left(\mathbb{H}_{+i}\right)$-norm}:is
finite. The author calls this the \textbf{$\frac{1}{y}$-asymptotic
$2$-Hardy space} on the upper half plane.

\vphantom{}\textbf{Theorem} \textbf{13}: Suppose $\left\{ \psi_{n}\right\} _{n\in\mathbb{N}_{1}}\subseteq\mathcal{H}_{1/y}^{2}\left(\mathbb{H}_{+i}\right)$,
and, for each $n$, let:
\[
R_{n}\left(x\right)\overset{\textrm{def}}{=}\begin{cases}
\textrm{VRes}_{x}\left[\psi_{n}\right] & \textrm{if }x\in\mathbb{Q}\cap\textrm{P}\left(\psi_{n}\right)\\
0 & \textrm{if }x\in\mathbb{Q}\backslash\textrm{P}\left(\psi_{n}\right)
\end{cases}
\]
If the $\psi_{n}$s converge in $\mathcal{H}_{1/y}^{2}\left(\mathbb{H}_{+i}\right)$
to a limit $\psi\in\mathcal{H}_{1/y}^{2}\left(\mathbb{H}_{+i}\right)$,
then the $R_{n}$s converge in $L^{2}\left(\mathbb{Q}/\mathbb{Z}\right)$
to a limit $R\in L^{2}\left(\mathbb{Q}/\mathbb{Z}\right)$.

Proof: Write:
\[
\psi_{n}\left(z\right)=\sum_{k=0}^{\infty}c_{n,k}e^{2k\pi iz}
\]
\[
\psi\left(z\right)=\sum_{k=0}^{\infty}c_{k}e^{2k\pi iz}
\]
If the $\psi_{n}$s converge to $\psi$ in $\mathcal{H}_{1/y}^{2}\left(\mathbb{H}_{+i}\right)$,
then they form a Cauchy sequence in $\mathcal{H}_{1/y}^{2}\left(\mathbb{H}_{+i}\right)$.
By \textbf{Lemma 7}, it follows that:
\[
\frac{1}{\sqrt{2\pi}}\left\Vert \psi_{m}-\psi_{n}\right\Vert _{2,\frac{1}{y}}\geq\left\Vert R_{m}-R_{n}\right\Vert _{\mathbb{Q}/\mathbb{Z}}
\]
and thus, that the $R_{n}$s are Cauchy in $L^{2}\left(\mathbb{Q}/\mathbb{Z}\right)$.
Since $L^{2}\left(\mathbb{Q}/\mathbb{Z}\right)$ is complete, this
guarantees that the $R_{n}$s converge to some limit $R\in L^{2}\left(\mathbb{Q}/\mathbb{Z}\right)$.

Q.E.D.

\vphantom{}\textbf{Corollary} \textbf{5}: Suppose the $\psi_{n}$s'
limit, $\psi$, is fixed by $\mathscr{Q}_{H}$. Then, the $R_{n}$s'
limit, $R$, is fixed by $\mathfrak{Q}_{H}$.

Proof: Since $\mathfrak{Q}_{H}$ is continuous on $L^{2}\left(\mathbb{Q}/\mathbb{Z}\right)$,
if the $R_{n}$s satisfy $\left\Vert \left(\mathfrak{Q}_{H}-1\right)\left\{ R_{n}\right\} \right\Vert _{\mathbb{Q}/\mathbb{Z}}\rightarrow0$
and if $R_{n}\rightarrow R$ in $L^{2}\left(\mathbb{Q}/\mathbb{Z}\right)$,
then $\left(\mathfrak{Q}_{H}-1\right)\left\{ R\right\} =0$, and so,
$R\in\textrm{Ker}\left(1-\mathfrak{Q}_{H}\right)$.

Q.E.D.

\vphantom{}\textbf{Conjecture} \textbf{1}: Let everything be as in
\textbf{Theorem} \textbf{13}. Then:
\[
R\left(x\right)=\textrm{VRes}_{x}\left[\psi\right]
\]
that is to say:
\[
\lim_{n\rightarrow\infty}\textrm{VRes}_{x}\left[\psi_{n}\right]=\textrm{VRes}_{x}\left[\lim_{n\rightarrow\infty}\psi_{n}\right]
\]
where the limit on the left is taken in $L^{2}\left(\mathbb{Q}/\mathbb{Z}\right)$
and the limit on the right is taken in $\mathcal{H}_{1/y}^{2}\left(\mathbb{H}_{+i}\right)$;
in words: \emph{the limit of the dreamcatchers is the dreamcatcher
of the limit}.\emph{ }(Might it be sufficient merely to prove that
$\mathcal{H}_{1/y}^{2}\left(\mathbb{H}_{+i}\right)$ is a complete
metric space?)

\vphantom{}\textbf{Question} \textbf{1}: Let $\left\{ f_{n}\left(z\right)\right\} _{n\in\mathbb{N}_{1}}$
be a sequence of non-polynomial rational functions---with $f_{n}\left(z\right)=\frac{P_{n}\left(z\right)}{Q_{n}\left(z\right)}$
for co-prime polynomials $P_{n},Q_{n}$ for all $n\in\mathbb{N}_{1}$---and
let $V\subseteq\mathbb{N}_{0}$ so that the trigonometric polynomials
$F_{n}\left(z\right)\overset{\textrm{def}}{=}f_{n}\left(e^{2\pi iz}\right)$
converge in $\mathcal{H}_{1/y}^{2}\left(\mathbb{H}_{+i}\right)$ to
the Fourier set-series $\psi_{V}\left(z\right)$. What can be said
about such a $V$? For example---using the notation of Chapter 4---is
it true that for every set $W\in\mathcal{M}_{0}$ with $d\left(W\right)\neq0$,
there exists a $V\subseteq\mathbb{N}_{0}$ with $V\approx W$ and
a sequence of $F_{n}$s so that $\lim_{n\rightarrow\infty}\left\Vert F_{n}-\psi_{V}\right\Vert _{2,\frac{1}{y}}=0$?
If this is not always true, then, for what sort of $W$ does it hold
true?

\pagebreak{}

\end{document}